\documentclass[11pt]{article}

\usepackage[T1]{fontenc}
\usepackage{lmodern}
\usepackage{microtype}
\usepackage{amsmath,amssymb,amsthm,mathtools}
\usepackage{tikz-cd}
\usepackage{xcolor}
\usepackage[margin=1in]{geometry}
\usepackage{enumitem}
\usepackage[
  colorlinks=true,
  linkcolor=magenta,
  citecolor=cyan,
  urlcolor=magenta,
  pdftitle={Virtual Two-Step Nilpotency of the Fundamental Groups of Special Quasi-Projective Varieties},
  pdfauthor={Tianshu Jiang and Jiayu Li}
]{hyperref}

\setlist[enumerate]{leftmargin=2.2em,itemsep=2pt,topsep=4pt}
\setlist[itemize]{leftmargin=2em,itemsep=2pt,topsep=4pt}
\allowdisplaybreaks
\numberwithin{equation}{section}
\theoremstyle{plain}
\newtheorem{theorem}{Theorem}[section]            % label kind: thm
\newtheorem{proposition}[theorem]{Proposition}    % label kind: prop
\newtheorem{lemma}[theorem]{Lemma}                % label kind: lem
\theoremstyle{definition}
\newtheorem{definition}[theorem]{Definition}      % label kind: def
\theoremstyle{remark}
\title{Virtual Two-Step Nilpotency of the Fundamental Groups of Special Quasi-Projective Varieties}
\author{Tianshu Jiang \and Jiayu Li}
\date{}

\providecommand{\CC}{\mathbb C} % ground field
\providecommand{\QQ}{\mathbb Q} % rational numbers
\providecommand{\RR}{\mathbb R} % real numbers
\providecommand{\ZZ}{\mathbb Z} % integers
\providecommand{\Orbgrp}[2]{\Pi_{\mathrm{orb}}(#1,#2)} % regular-complement meridian-power quotient
\providecommand{\Supp}{\operatorname{Supp}} % support
\providecommand{\Tor}{\operatorname{Tor}} % torsion subgroup of an abelian group
\providecommand{\im}{\operatorname{im}} % image
\providecommand{\rk}{\operatorname{rk}} % rank
\providecommand{\vb}{\operatorname{vb}_1} % supremal first Betti number of connected finite etale covers

\begin{document}
\maketitle
\begin{abstract}
We prove that the fundamental group of every smooth connected complex
quasi-projective variety that is special in the sense of Campana is
virtually nilpotent of class at most two. This establishes the two-step
nilpotency conjecture of Cao--Deng--Hacon--P\u{a}un, a quasi-projective
analogue of Campana's almost abelian conjecture for special compact
K\"ahler manifolds. Examples of Cadorel--Deng--Yamanoi show that the
step bound is sharp. On a suitable finite etale cover, we construct a
canonical second Albanese map to a torus torsor over an abelian variety.
It suffices to prove that the induced homomorphism on fundamental
groups has finite kernel. Using good minimal models and meridian
groups of projective pairs, we show that
an infinite kernel would yield a saturated rank-one subsheaf of
logarithmic differentials whose Iitaka dimension contradicts specialness.
\end{abstract}
\medskip
{\small
\par
}
\tableofcontents

\section{Introduction}\label{sec:01-introduction}

A smooth complex quasi-projective variety $X$ is \emph{special} in
the sense of Campana~\cite{campana2004,campana2009} if, on a smooth
projective compactification $\overline X$ with reduced simple normal
crossing boundary $D=\overline X\setminus X$, every saturated rank-one
subsheaf $\mathcal L\subseteq\Omega^p_{\overline X}(\log D)$ satisfies
$\kappa(\overline X,\mathcal L)<p$ for $1\le p\le\dim X$, where
$\kappa(\overline X,\mathcal L)$ denotes the Iitaka dimension of
$\mathcal L$.
We prove the following theorem on its fundamental group.

\begin{theorem}\label{thm:introduction-main}\label{thm:main}
Let $X$ be a smooth connected special complex quasi-projective
variety. Then $\pi_1(X)$ has a finite-index subgroup $\Gamma_0$ such that
\[
 [\Gamma_0,[\Gamma_0,\Gamma_0]]=\{1\}.
\]
\end{theorem}

\medskip\noindent\textbf{Previous work.}
Campana conjectured that the fundamental group of a special compact
K\"ahler manifold is virtually abelian
\cite[Conjecture~7.1]{campana2004}. He proved the corresponding assertion
for finite-dimensional complex linear images
\cite[Theorem~7.8]{campana2004}. The compact conjecture was proved in
all dimensions by OpenAI~\cite[Theorem~1.1]{openai2026d}.

For quasi-projective varieties, Cadorel--Deng--Yamanoi constructed a
special surface whose fundamental group is nilpotent but not virtually
abelian~\cite[Example~4.27]{cadorel2026}. They proposed virtual
nilpotency as the replacement for virtual abelianity
\cite[Conjecture~2.6]{cadorel2026} and proved that every
finite-dimensional complex linear representation of the fundamental
group of a special quasi-projective variety has virtually nilpotent
image~\cite[Theorem~A]{cadorel2026}.

Cao--Deng--Hacon--P\u{a}un sharpened this result by proving that these
linear images are virtually nilpotent of class at most two
\cite[Theorem~B]{cao2026}. They formulated the same conclusion for the
entire fundamental group of a special quasi-projective variety as the
two-step nilpotency conjecture~\cite[Conjecture~2.7]{cao2026}.
Theorem~\ref{thm:introduction-main} establishes this conjecture. The
examples of Cadorel--Deng--Yamanoi show that the bound on the nilpotency
class is sharp.

\medskip\noindent\textbf{Sketch of proof.}\par

\smallskip\noindent\emph{Reduction to a finite kernel.}
Choose a connected finite etale cover $Y_0\to X$ on which the
first Betti number is maximal. Let $a_0:Y_0\to A_0$ be the quasi-Albanese map
and $F_0$ a smooth general fiber. Put
\[
 \begin{gathered}
 G_0=\pi_1(Y_0),\qquad
 K_0=\im\bigl(\pi_1(F_0)\to G_0\bigr),\\
 P_0=\ker\bigl(K_0\to (K_0)_{\rm ab}/\Tor((K_0)_{\rm ab})\bigr).
 \end{gathered}
\]
Theorem~\ref{thm:canonical-albanese} constructs the canonical second
Albanese map $f_0:Y_0\to Q_0$, where $Q_0$ is a torus torsor over an
abelian variety and
\[
 1\longrightarrow P_0\longrightarrow G_0
 \xrightarrow{\ f_{0*}\ }E_0=\pi_1(Q_0)\longrightarrow1.
\]
The group $E_0$ is torsion-free and nilpotent of class at most two.
The map $f_0$ is dominant, and $P_0$ is the image of the fundamental
group of a smooth general $f_0$-fiber, hence finitely generated.
A group is called \emph{FAb} if every finite-index subgroup has finite
abelianization. Proposition~\ref{prop:canonical-kernel-covers} proves
that $P_0$ is FAb.

It remains to prove that $P_0$ is finite. If so, the same proposition
gives a finite-index subgroup $H\le G_0$ with $H\cap P_0=1$.
The embedding of $H$ in $E_0$ then proves
Theorem~\ref{thm:introduction-main}. Suppose, for a
contradiction, that $P_0$ is infinite.

\smallskip\noindent\emph{A quotient of minimal relative dimension.}
For a good fiber $S_0$ of $f_0$ as in
Theorem~\ref{thm:ordinary-source-minimum}, consider a connected rational
fibration $h:S_0\dashrightarrow W_0$ with smooth projective base.
Let $C_h\triangleleft P_0$ be the image
of the fundamental group of a whole resolved general $h$-fiber.
Set
\begin{equation}\label{eq:sketch-minimal-base-dimension}
 \ell=\min\{\dim W_0:P_0/C_h\text{ is infinite}\},
\end{equation}
where the minimum is taken over all such fibers and fibrations.
A smooth projective compactification of $S_0$ supplies a candidate. Thus
$1\le\ell\le\dim Y_0-\dim Q_0$.

For a minimizing $h$, put $C=C_h$. A Hilbert scheme argument shows that
the normalizer $G=N_{G_0}(C)$ of $C$ in $G_0$ has finite index in $G_0$
and contains $P_0$.
Let $\rho:Q\to Q_0$ be the finite etale cover corresponding to the
subgroup $G/P_0\le\pi_1(Q_0)$. On the connected pullback
$Y=Y_0\times_{Q_0}Q$, the kernel is still $P_0$, and the fibers of
$f:Y\to Q$ are unchanged. Applying the Shafarevich map to $G/C$,
followed by a relative projective compactification, contracts the
chosen leaves and gives a rational factorization
\[
\begin{tikzcd}[column sep=large,row sep=large]
Y \arrow[r,dashed,"c"] \arrow[d,"\rho_Y"']
  \arrow[rr,bend left=22,"f"] &
Z \arrow[r,"g"] & Q \arrow[d,"\rho"] \\
Y_0 \arrow[rr,"f_0"'] && Q_0 .
\end{tikzcd}
\]
The outer square is Cartesian, and $g$ is projective and surjective;
both $c$ and $g$ have connected general fibers.
If $B_c\triangleleft G$ is the image of the fundamental group of a
resolved general $c$-fiber, then
\[
 C\subseteq B_c\subseteq P_0,\qquad B_c/C\text{ is finite},
 \qquad \Pi=P_0/B_c\text{ is infinite}.
\]
As a quotient of $P_0$, the group $\Pi$ is also FAb
(Theorem~\ref{thm:ordinary-source-minimum}(i)).
The contraction of the chosen leaves and the original minimality give
$\dim Z-\dim Q=\ell$. In particular, for any good $f$-fiber $S$ and
any connected rational fibration $S\dashrightarrow T_0$ with smooth
projective $T_0$ and $0<\dim T_0<\ell$, the image of its whole resolved general fiber in
$\Pi$ has finite index.

\smallskip\noindent\emph{Weighted projective fibers.}
Keep $Y,Q,c$ and the quotient $R=G/B_c$ fixed, so that
\[
 1\longrightarrow\Pi\longrightarrow R
 \longrightarrow\pi_1(Q)\longrightarrow1.
\]
Choose a smooth model $Z_m$ projective over $Q$ and a birational map
$\beta_m:Z_m\dashrightarrow Z$ over $Q$.
Resolving $c^{(m)}=\beta_m^{-1}\circ c$ gives a proper birational morphism
$\mu_m:M_m\to Y$, with $M_m$ smooth, and a morphism
$c_m:M_m\to Z_m$. We retain every exceptional divisor of $\mu_m$.
The maps satisfy $g_mc_m=f\mu_m$ and fit into the diagram
\[
\begin{tikzcd}[column sep=large,row sep=large]
M_m \arrow[r,"c_m"] \arrow[d,"\mu_m"'] &
Z_m \arrow[r,"g_m"] \arrow[d,dashed,"\beta_m"] & Q \arrow[d,equal] \\
Y \arrow[r,dashed,"c"'] & Z \arrow[r,"g"'] & Q .
\end{tikzcd}
\]
We choose a dense open
$Z_m^\circ=Z_m\setminus(D_m\cup\Sigma_m)$ over which $c_m$ is a smooth
surjective topological fiber bundle with connected fibers. Here $D_m$
is the reduced divisorial part of the excluded locus, arranged to have
simple normal crossings, and $\operatorname{codim}\Sigma_m\ge2$.
This is a prepared model for $R$ in the sense of
Definition~\ref{def:prepared-ordinary-model}.

Over a very general $u\in Q$ we have a commutative diagram
\[
\begin{tikzcd}[column sep=large,row sep=large]
S \arrow[r,dashed,"a_m"] \arrow[d,hook] &
V \arrow[r] \arrow[d,hook] & \{u\} \arrow[d,hook] \\
Y \arrow[r,dashed,"c^{(m)}"'] & Z_m \arrow[r,"g_m"'] & Q ,
\end{tikzcd}
\]
where $S=f^{-1}(u)$ is smooth quasi-projective and
$V=g_m^{-1}(u)$ is smooth projective of dimension $\ell$.
The restriction $D_m|_V$ has simple normal crossings.
For a prime component $D$ of $D_m$, let $m_D$ be the order in $R$
of the image of a meridian around $D$, after imposing the relations
from all retained source divisors, allowing $m_D=\infty$. Set
\[
 \Delta_m=\sum_D(1-1/m_D)D,\qquad \Delta=\Delta_m|_V.
\]
A boundary with coefficients $1-1/m$, where
$m\in\ZZ_{\ge2}\cup\{\infty\}$, is called \emph{standard};
we set $1/\infty=0$, and unmarked prime divisors have weight one.
For a normal projective variety $B$ with such a boundary
$C=\sum_i(1-1/m_i)D_i$, we use the meridian group
\begin{equation}\label{eq:meridian-group}
 \Orbgrp{B}{C}=\pi_1(B_{\rm reg}\setminus\Supp C)
 /\langle\!\langle\mu_i^{m_i}:m_i<\infty\rangle\!\rangle,
\end{equation}
where $\mu_i$ winds once positively around a general point of $D_i$.
Divisors of infinite weight remain deleted and impose no power relation.
The pair $(V,\Delta)$ has SNC standard boundary, and there are quotient maps
\[
 \pi_1(S)\twoheadrightarrow\Pi,
 \qquad \Orbgrp V\Delta\twoheadrightarrow\Pi.
\]
These maps are compatible with $a_m$ on a dense open subset of $S$.
By transversality, a meridian in $V\cap Z_m^\circ$ around a component
of $D\cap V$ maps to a conjugate of a meridian around $D$ in
$\pi_1(Z_m^\circ)$. Its image in $R$ lies in the subgroup $\Pi$.
Since $\Pi\hookrightarrow R$ is injective, its order in $\Pi$ is
exactly $m_D$, including when $m_D=\infty$.

The construction applies to every prepared model
$g_m:Z_m\to Q$ of Definition~\ref{def:prepared-ordinary-model}.

\smallskip\noindent\emph{Minimality forces bigness.}
If $K_V+\Delta$ is pseudoeffective, we use a $(K_V+\Delta)$-minimal
model program ending at a good minimal model $(M,\Delta_M)$, with
connected semiample contraction $p:M\to T$ and
$\dim T=\kappa(V,K_V+\Delta)$. This form of the minimal model program
follows from OpenAI~\cite[Theorem~11.1]{openai2026b} and the reduction
in the proof of Theorem~\ref{thm:source-minimal-bigness}.
Otherwise, the program with ample scaling terminates at a Mori
fiber space $p:(M,\Delta_M)\to T$
\cite[Section~2.20]{fujino2014b}.
In either case, set $d=\dim T$. The endpoint pair is log canonical,
and the birational steps extract no divisors and preserve the coefficients of the
surviving boundary components. Proposition~\ref{prop:ordinary-meridian-comparison}
therefore gives the meridian comparison
\[
 \Orbgrp M{\Delta_M}\twoheadrightarrow\Orbgrp V\Delta
 \twoheadrightarrow\Pi.
\]
We use the FAb property of $\Pi$ to rule out $d<\ell$.
Theorem~\ref{thm:antinef-group} implies that every FAb quotient of
the group of a projective lc standard pair with antinef adjoint is
finite. If $d=0$, then $K_M+\Delta_M$ is rationally trivial or
anti-ample. The theorem therefore implies that the FAb quotient
$\Pi$ of $\Orbgrp M{\Delta_M}$ is finite, a contradiction.

Suppose $0<d<\ell$. A very general $p$-fiber $(F,\Delta_F)$ is a
projective lc standard pair. Its adjoint is rationally trivial in the
semiample case and anti-ample in the Mori case. Passing to a smooth
projective resolution $\widehat T\to T$, the composite through
$S\dashrightarrow V\dashrightarrow M\to T$ gives a connected rational
fibration $S\dashrightarrow\widehat T$ with $\dim\widehat T=d$.
The group
\[
 J=\im\bigl(\Orbgrp F{\Delta_F}\to\Pi\bigr)
\]
contains, by compatibility on dense opens, the image $I$ of a whole
resolved general fiber of this composite. The minimality of $\ell$
in \eqref{eq:sketch-minimal-base-dimension} gives $[\Pi:I]<\infty$
(Theorem~\ref{thm:ordinary-source-minimum}(iii)). Hence $J$ has finite
index in $\Pi$. Since the FAb property passes to
finite-index subgroups, $J$ is FAb. As a quotient of
$\Orbgrp F{\Delta_F}$, it is therefore finite by
Theorem~\ref{thm:antinef-group}, again contradicting the infinitude of
$\Pi$. The Mori case, for which $d<\ell$, is excluded. In the
semiample case we must have $d=\ell$, so
\[
 \kappa(V,K_V+\Delta)=\ell.
\]
Thus the weighted adjoint is big on a very general projective fiber of
\emph{every} such model $Z_m$ (Theorem~\ref{thm:source-minimal-bigness}).

\smallskip\noindent\emph{The contradiction with specialness.}
We use this relative bigness to construct logarithmic differential
forms on a finite cover of the source.
For each model $Z_m$, choose a smooth projective compactification
$\overline Z_m$ such that the reduced boundary
$D_{\infty,m}=(\overline Z_m\setminus Z_m)_{\rm red}$ and the closure
$\overline{\Delta_m}$ have combined SNC support. Following
Definition~\ref{def:canonical-kernel-logarithmic-pair}, set
\[
 \kappa_{\log}(Z_m,\Delta_m)
 =\kappa\!\left(\overline Z_m,
 K_{\overline Z_m}+D_{\infty,m}+\overline{\Delta_m}\right).
\]
This Iitaka dimension is independent of the chosen
compactification (Lemma~\ref{lem:canonical-kernel-logarithmic-comparison}).
Over all such models of the fixed map $g:Z\to Q$, set
\begin{equation}\label{eq:sketch-minimal-log-dimension}
 b=\min_m\kappa_{\log}(Z_m,\Delta_m)\ge\ell>0,
\end{equation}
where the lower bound follows from
Lemma~\ref{lem:ordinary-minimum-mixed-leaves}(i).
Lemma~\ref{lem:ordinary-minimum-mixed-leaves} and the construction in
the proof of Theorem~\ref{thm:cofinal-specialness-obstruction} allow us
to choose a smooth prepared model $U$ for the same quotient $R$, with
actual meridian boundary $\Delta_U$ and
$\kappa_{\log}(U,\Delta_U)=b$. Its logarithmic Iitaka map is a
dominant morphism $\varphi:U\to B$ to a smooth projective $b$-fold.
The map $g_U:U\to Q$ is projective, and the source map lifts to
$c_U:Y\dashrightarrow U$ with the same general-fiber group image $B_c$:
\begin{equation}\label{eq:sketch-iitaka-resolution}
\begin{tikzcd}[column sep=large,row sep=large]
Y \arrow[r,dashed,"c_U"] \arrow[dr,"f"'] &
U \arrow[r,"\varphi"] \arrow[d,"g_U"'] & B \\
& Q &
\end{tikzcd}
\end{equation}
The very general fibers $H$ of $\varphi$ are smooth and connected;
we call them the \emph{Iitaka leaves}. The boundary $\Delta_U$ is
vertical over $B$, so these leaves avoid its support. Each whole leaf
$H$ maps properly to $Q$, and $\pi_1(H)\to\pi_1(Q)$ is injective.
The image of $\pi_1(H)$ in $R$ satisfies
\[
 N=\im\bigl(\pi_1(H)\to R\bigr)\triangleleft R,
 \qquad N\cap\Pi=1.
\]

As in the proof of Theorem~\ref{thm:cofinal-specialness-obstruction},
Lemmas~\ref{lem:ordinary-minimum-finite-leaves}(iii) and
\ref{lem:canonical-kernel-family-splitting} give a
single finite etale cover
$\rho':Q'\to Q$ and a smooth surjection $r:Q'\to Q_1$ such that
each very general connected lifted leaf maps properly birationally
onto a whole fiber of $r$. Thus $Q_1$ parametrizes the images of
the lifted leaves in $Q'$; it is again a torus torsor over an abelian
variety. The construction also factors $r$ through a principal torus
torsor $Q'\to A\times Q_1$, where $A$ is an abelian variety.
We make the same base change on the source and on $U$:
\begin{equation}\label{eq:sketch-iitaka-cover}
\begin{tikzcd}[column sep=large,row sep=large]
Y' \arrow[r,dashed,"c'_U"] \arrow[d]
  & U' \arrow[r,"g'"] \arrow[d] & Q' \arrow[d,"\rho'"] \\
Y \arrow[r,dashed,"c_U"'] & U \arrow[r,"g_U"'] & Q .
\end{tikzcd}
\end{equation}
Here $Y'=Y\times_QQ'$ and $U'=U\times_QQ'$, and the vertical maps
are connected finite etale covers. The general $c$-fibers lift unchanged,
so their group image is still $B_c$. Put
$R'=\pi_1(Y')/B_c\le R$; it contains $\Pi$. The meridian boundary
$\Delta_{U'}$ is the pullback of $\Delta_U$, with unchanged orders.

The logarithmic Iitaka fibration of $(U',\Delta_{U'})$ has the connected
lifted leaves as its very general fibers, and its base still has dimension $b$.
Since $rg'$ is constant on these leaves, it induces a rational map
from this base to $Q_1$. The construction in the proof of
Theorem~\ref{thm:cofinal-specialness-obstruction} gives a smooth
birational model $W$ of this base, with a projective morphism
$s:W\to Q_1$.
Thus $W$ parametrizes the lifted leaves, while $s$ records their
images in $Q'$. Set $E=Q'\times_{Q_1}W$. The birational maps of the
leaves onto the $r$-fibers give a birational map $U'\dashrightarrow E$:
\begin{equation}\label{eq:sketch-iitaka-pullback}
\begin{tikzcd}[column sep=large,row sep=large]
Y' \arrow[r,dashed,"c'_U"] &
U' \arrow[r,dashed,"\sim"] \arrow[dr,"g'"'] &
E \arrow[r,"p"] \arrow[d,"q"] &
W \arrow[d,"s"] \\
&& Q' \arrow[r,"r"'] & Q_1 .
\end{tikzcd}
\end{equation}
The right square is Cartesian, and the composite is a
dominant rational map $\psi:Y'\dashrightarrow W$ with $\dim W=b>0$.

We define the meridian boundary on $W$ using $\psi$.
Let $K'$ be the image in $\pi_1(Y')$ of a whole resolved general
$\psi$-fiber. The image of a connected lifted Iitaka leaf in $R'$ is
$N'=N\cap R'$, which is normal in $R'$ and satisfies $N'\cap\Pi=1$.
Lemma~\ref{lem:ordinary-minimum-mixed-leaves}(ii) gives
$K'/B_c=N'$, hence $\pi_1(Y')/K'\simeq R'/N'$. In particular,
\begin{equation}\label{eq:sketch-leaf-quotient}
\begin{tikzcd}[column sep=large,row sep=large]
\Pi \arrow[r,hook] \arrow[d,equal] & R' \arrow[d,two heads] \\
\Pi \arrow[r,hook] & R'/N'\simeq\pi_1(Y')/K' .
\end{tikzcd}
\end{equation}
Thus passing to the quotient associated with $\psi$ preserves the
meridian orders in $\Pi$. Let $\Delta_E$ be the meridian boundary on
$E$, with weights computed in $R'$. The verticality of $\Delta_U$,
together with the birational comparison in the proof of
Theorem~\ref{thm:cofinal-specialness-obstruction}, gives
\mbox{$\Delta_E=p^*\Theta_W$}. The coefficients of $\Theta_W$ are
defined by the actual meridian orders in $\pi_1(Y')/K'$.

The same construction works after any further projective birational
morphism $\eta:W_1\to W$ with $W_1$ smooth. Taking a further
resolution when needed to make the boundary SNC, put
$E_1=Q'\times_{Q_1}W_1$. The corresponding boundaries satisfy
$\Delta_{E_1}=p_1^*\Theta_{W_1}$, with weights computed in the same
quotient $\pi_1(Y')/K'$:
\begin{equation}\label{eq:sketch-iitaka-refinement}
\begin{tikzcd}[column sep=large,row sep=large]
E_1 \arrow[r] \arrow[d,"p_1"'] & E \arrow[d,"p"] \\
W_1 \arrow[r,"\eta"'] & W .
\end{tikzcd}
\end{equation}
This square is Cartesian, and $E_1\to A\times W_1$ is again a
principal torus torsor. The section comparison in the proof of
Theorem~\ref{thm:cofinal-specialness-obstruction}, the minimality of
$b$ in \eqref{eq:sketch-minimal-log-dimension}, and
Lemma~\ref{lem:canonical-kernel-logarithmic-comparison} give
\[
 \kappa_{\log}(W_1,\Theta_{W_1})
 =\kappa_{\log}(E_1,p_1^*\Theta_{W_1})
 \ge b=\dim W_1.
\]
Hence the full logarithmic adjoint on a smooth projective
compactification of every such $W_1$ is big. We may therefore make
further base modifications to pull its weighted logarithmic sections
back to $Y'$.

Fix a smooth projective compactification $(\overline Y',D_Y)$ with
full reduced SNC boundary. Using the base preparation in the proof of
Theorem~\ref{thm:cofinal-specialness-obstruction}, we choose a higher
model $W_1$ with a smooth projective compactification $\overline W_1$
such that every prime divisor of the fixed $\overline Y'$ either
dominates $\overline W_1$ or dominates one of its prime divisors.
Every prime divisor in $Y'$ has its generic image in $W_1$.
Write $D_{W_1}$ for the reduced boundary at infinity; we also arrange
that $D_{W_1}+\overline\Theta_{W_1}$ has SNC support.
Resolving the graph of the induced rational map gives
\begin{equation}\label{eq:sketch-special-source-pullback}
\begin{tikzcd}[column sep=large,row sep=large]
& \widehat Y' \arrow[dl,"\mu"'] \arrow[dr,"\widehat\psi"] & \\
\overline Y' \arrow[rr,dashed,"\overline\psi_1"'] && \overline W_1 .
\end{tikzcd}
\end{equation}
Here $\mu$ is proper birational and $\widehat\psi$ is a morphism
resolving $\overline\psi_1$. We use $\widehat\psi$ to compute
pullbacks and identify them, through $\mu$, with rational
differential tensors on the fixed $\overline Y'$.

On a dense open in $Y'$ where $\overline\psi_1$ is smooth, the
pullbacks of local top forms on the $b$-dimensional base span a
rank-one subbundle of $\Omega^b_{\overline Y'}$. Let
\[
 \mathcal L\subseteq\Omega^b_{\overline Y'}(\log D_Y)
\]
be the saturation of this generic line: its local sections are the
logarithmic $b$-forms whose generic values lie in that line.
Pullbacks of weighted
logarithmic sections are initially rational sections of
$\mathcal L^{\otimes k}$; we check their poles along the prime
divisors of $\overline Y'$.

For a prime $D\subset W_1$, let $m_D$ be its
meridian order in $\pi_1(Y')/K'$, so its coefficient in
$\Theta_{W_1}$ is $1-1/m_D$. Suppose that a prime divisor
$F\subset Y'$ dominates $D$. At its general point, write
$\overline\psi_1^*t=x^{e_F}u$, where $t=0$ defines $D$, $x=0$
defines $F$, and $u$ is a unit. Thus $e_F$ is the vanishing order
along $F$ of the pulled-back defining equation.
A transverse circle around $F$ bounds a disk in $Y'$ and maps
to the $e_F$-th power of a meridian of $D$. Hence the image of this
meridian in $\pi_1(Y')/K'$ has $e_F$-th power one, so $m_D$ is
finite and divides $e_F$.
In particular, a prime of infinite meridian order has no such
interior divisor above it. The inequality
\[
 e_F-1\ge e_F\left(1-\frac1{m_D}\right)
\]
shows that the zeros introduced by differential pullback cancel
the allowed weighted poles at $F$.

At a component of $D_Y$ lying over a base divisor, the identity
$d(x^eu)=x^e(eu\,dx/x+du)$ gives vanishing of order at least $e$
in a logarithmic frame, enough for every coefficient at most one,
including the boundary at infinity. A source prime dominating
the base introduces no pole from a base divisor. By the base
preparation, these checks cover every prime of the fixed
$\overline Y'$. Normality therefore extends the pullbacks to
global sections of $\mathcal L^{\otimes k}$, giving injections
preserving section ratios for sufficiently divisible $k$:
\[
 H^0\!\left(\overline W_1,
 k(K_{\overline W_1}+D_{W_1}+\overline\Theta_{W_1})\right)
 \hookrightarrow H^0(\overline Y',\mathcal L^{\otimes k}).
\]
Bigness on $\overline W_1$ implies
$\kappa(\overline Y',\mathcal L)\ge b$.
Since $Y'$ is a finite etale cover of the special variety $Y$,
it is special
\cite[Proposition~10.11]{campana2009}, a contradiction.

\section{The canonical second Albanese map and its kernel}\label{sec:09-canonical-kernel}

We construct the canonical second Albanese map to a torus torsor over
an abelian variety and identify the image of the fundamental group of
a general fiber with the kernel on fundamental groups. The FAb and
finite-cover properties of this kernel reduce the main theorem to
its finiteness.

\subsection{Canonical quotient and Albanese realization}\label{subsec:canonical-quotient-realization}

\begin{definition}\label{def:canonical-kernel}
Let $X$ be smooth, connected, complex quasi-projective and special.
Write $\vb(X)$ for the supremum of $b_1(X')$ over connected finite
etale covers $X'\to X$.
Let $Y\to X$ be a connected finite etale cover satisfying
$b_1(Y)=\vb(X)$, and let $a:Y\to A$ be its based quasi-Albanese morphism,
the universal based algebraic morphism to a semiabelian variety.
Choose the base point $y$ in a smooth connected general fiber $F$ of $a$.
Set
\[
 G=\pi_1(Y),\qquad K=\im(\pi_1(F)\to G),\qquad
 P=\ker\bigl(K\to K_{\rm ab}/\Tor K_{\rm ab}\bigr),\qquad E=G/P.
\]
Here $K_{\rm ab}=K/[K,K]$. We call these groups, together with the
chosen cover and morphism, the canonical kernel data.
\end{definition}

Thus $K/P$ is the maximal torsion-free abelian quotient of the fiber
image $K$.

For a based smooth connected quasi-projective variety $(V,v)$, let
$U_2(V,v)$ be the rational two-step unipotent (Malcev) completion of
$\pi_1(V,v)$, and let $\Gamma_2(V,v)$ be the image of $\pi_1(V,v)$ in
$U_2(V,v)(\QQ)$. Write $\mathfrak u_2(V,v)=\operatorname{Lie}U_2(V,v)$,
equipped with the mixed Hodge structure of Morgan
\cite[Theorem~9.2]{morgan1978} and Hain
\cite[Theorem~5.1]{hain1987geometry}.
The Hodge filtration makes $F^0\mathfrak u_2(V,v)_{\CC}$ a Lie
subalgebra \cite[Corollary~3.26]{hainzucker1987}. Following Hain--Zucker
\cite[Definition~5.15]{hainzucker1987}, the \emph{based second Albanese
manifold} is
\begin{equation}\label{eq:based-second-albanese}
 \operatorname{Alb}_2(V,v)
 =\Gamma_2(V,v)\backslash U_2(V,v)(\CC)
   /\exp\bigl(F^0\mathfrak u_2(V,v)_{\CC}\bigr).
\end{equation}
Their construction gives a based holomorphic map
$a_{2,v}:V^{\rm an}\to\operatorname{Alb}_2(V,v)$, sending $v$ to
the identity double coset and recovering the usual based quasi-Albanese
map after projection to level one
\cite[(5.17), Corollaries~5.20 and~5.40]{hainzucker1987}.
Hain gives an account of these manifolds and maps in
\cite{hain1987albanese}.

\begin{theorem}\label{thm:canonical-albanese}
Let $X$ be a smooth connected special complex quasi-projective
variety. There is a connected finite etale cover $Y\to X$ with
$b_1(Y)=\vb(X)$. For any canonical kernel data as in
Definition~\ref{def:canonical-kernel}, the following hold.
\begin{enumerate}[label=(\roman*)]
\item The group $K$ is finitely generated, $P$ is normal in $G$, and
there is a central extension
\[
 1\longrightarrow K_{\rm ab}/\Tor K_{\rm ab}
 \longrightarrow E\longrightarrow\pi_1(A)\longrightarrow1
\]
with finitely generated free abelian kernel and quotient.
In particular $E$ is torsion-free of nilpotency class at most two.
\item Let $\overline Y$ be a smooth projective compactification and
$B=\operatorname{Alb}(\overline Y)$. The based second Albanese
$\operatorname{Alb}_2(Y,y)$ admits an algebraic realization as a smooth
quasi-projective variety $Q$, with an algebraic principal torus torsor
$q:Q\to B$ and algebraic maps $f:Y\to Q$ and $p:Q\to A$ satisfying
$pf=a$. The analytic space $Q^{\rm an}$ is a $K(E,1)$, meaning that
$\pi_1(Q^{\rm an})\simeq E$ and $\pi_i(Q^{\rm an})=0$ for every $i\ge2$.
The homomorphism $f_*:G\to\pi_1(Q)=E$ is the canonical quotient.
This realization uses the given cover $Y$, abelian variety $B$, and lattice $E$.
\item The map $f$ is dominant with smooth geometrically connected general
fiber $S$. The image of its fundamental group is
\[
 \im\bigl(\pi_1(S)\to G\bigr)=P.
\]
In particular $P$ is finitely generated.
\end{enumerate}
\end{theorem}

We prove parts (i) and (ii) in this subsection. The proof of part (iii)
is completed in Subsection~\ref{subsec:kernel-exactness}, after the
required geometry of torus torsors and special images has been established.

\begin{proof}[Proof of existence and part (i) of Theorem~\ref{thm:canonical-albanese}]
\emph{The integral first Albanese.}
The quotient in (i) depends on the integral period lattice of the
quasi-Albanese, which we now construct.
For a smooth connected quasi-projective $V$, choose a smooth projective
compactification $\overline V$ with full SNC boundary $D$, and put
$H=H_1(V,\ZZ)$, $\overline H=H_1(\overline V,\ZZ)$ and
$H_{\mathrm f}=H/\Tor H$, with the analogous notation for
$\overline H$. Set
\[
 \mathcal L=H^0(\overline V,\Omega^1_{\overline V}(\log D)),\quad
 \mathcal H=H^0(\overline V,\Omega^1_{\overline V}),\quad
 q=\dim\mathcal H,\quad d=\dim\mathcal L-q.
\]
Logarithmic Hodge degeneration and strictness
\cite[Theorems~2.3.5 and~3.2.5, Corollaries~3.2.13--3.2.15]{deligne1971}
show that these forms are closed and that $b_1(V)=2q+d$. They also
identify $\mathcal L$ with $F^1H^1(V,\CC)$, independently of
compactification.
Loops in $\overline V$ avoid $D$, so $H\to\overline H$ is onto.
The kernel
\[
 K_{\mathrm s}=\ker(H_{\mathrm f}\to\overline H_{\mathrm f})
\]
is a saturated direct summand of rank $d$. Choose representatives
$\eta_1,\ldots,\eta_d\in H$ of its basis, and lifts
$\xi_1,\ldots,\xi_{2q}$ of a basis of $\overline H_{\mathrm f}$.
The images of the $\eta_j$ in $\overline H$ are torsion;
their periods against $\mathcal H$ therefore vanish. Every positive
meridian $\delta_a$ of a component $D_a$ has the form
\[
 \delta_a=\sum_jm_{ja}\eta_j+\tau_a,
 \qquad m_{ja}\in\ZZ,\quad \tau_a\in\Tor H.
\]
Integration identifies $\mathcal L/\mathcal H$ with
$\operatorname{Hom}(K_{\mathrm s},\CC)$: a form with zero periods
on $K_{\mathrm s}$ has zero residues on every $D_a$, hence is
holomorphic, and both spaces have dimension $d$.
Choose $\varphi_j\in\mathcal L$ with
$\int_{\eta_i}\varphi_j=\delta_{ij}$ and a basis of $\mathcal H$.
In the dual coordinates the entire period lattice is
\[
 \Lambda=\sum_{i=1}^{2q}\ZZ(a_i,b_i)
               +\sum_{j=1}^d\ZZ(0,e_j)
       \subset\CC^q\times\CC^d,
\]
where the $a_i$ are the ordinary periods and the $e_j$ are the
standard basis. Projective Hodge decomposition identifies
$H_1(\overline V,\RR)$ with $\mathcal H^\vee$ as real vector spaces
by integration, so the $a_i$ are a real basis. An ample class
polarizes this Hodge structure; rescaling the dual rational
polarization makes it integral on the same lattice. The Riemann-form
criterion \cite[Chapter~I, Theorems~2.8--2.9]{milne2008course}
identifies the quotient as the abelian variety
$B_V=\operatorname{Alb}(\overline V)$.
The holomorphic integration map $\overline V\to B_V$ has an
algebraic graph by projective GAGA
\cite[Section~12, Theorems~2--3]{serre1956}, so it is algebraic.
In particular the displayed generators of $\Lambda$ are
real-linearly independent. Thus $\Lambda$ is discrete
and the kernel of the period map on $H$ is exactly $\Tor H$.
Projection gives a principal $(\CC^*)^d$-bundle
$\mathcal L^\vee/\Lambda\to B_V$ with its full integral kernel.

To algebraize this bundle, exponentiate the last $d$ coordinates.
Its $j$th factor is the frame bundle of the flat holomorphic line
bundle on $B_V$ with character $a_i\mapsto\exp(2\pi i b_{ij})$.
Projective GAGA \cite[Section~12, Theorems~2--3]{serre1956}
algebraizes these line bundles. Their character identities give
holomorphic isomorphisms
$\operatorname{pr}_1^*L_j\otimes\operatorname{pr}_2^*L_j
\simeq m_{B_V}^*L_j$ on $B_V\times B_V$, and the analogous
inverse isomorphisms. By GAGA on these projective bases, the group operations
are algebraic. Hence the original quotient is a semiabelian
variety $A_V$, with connected torus kernel over $B_V$.
Based integration defines $\alpha_V:V\to A_V$ on all of $V$.
We apply the algebraicity argument of
\cite[Lemma~3.10]{fujino2024b} with the saturated basis chosen above.
Its local calculation only needs the integral meridian expansion:
the frame coordinates are constants times
$\prod_a z_a^{m_{ja}}\exp(2\pi i g_j)$, with $g_j$ holomorphic.
Thus they have finite divisorial poles on $\overline V$. After
twisting by those poles, projective GAGA shows that they are algebraic.
Their restrictions and inverses are regular on $V$.
Continuation of integration along loops gives
\[
 H_{\mathrm f}\xrightarrow[\alpha_{V*}]{\simeq}
 H_1(A_V,\ZZ)=\pi_1(A_V)=\Lambda.
\]
Strictness gives the corresponding isomorphism of global
logarithmic one-forms, as in
\cite[Lemmas~3.11--3.12]{fujino2024b}.

To prove the universal property, we use the following criterion.
Every based algebraic morphism between semiabelian varieties is a
group homomorphism: apply
\cite[Lemma~5.4.8]{brion2017structure} to its composition with
addition. Smooth connected subgroups and quotients are semiabelian
\cite[Corollary~5.4.6, Propositions~2.7.4 and~2.8.4]{brion2017structure}.
For a semiabelian variety of abelian dimension $g$ and torus rank $t$,
the bundle homotopy sequence gives $\pi_2=0$ and
$\pi_1=H_1\simeq\ZZ^{2g+t}$.
If a based algebraic map of semiabelian varieties induces an
isomorphism on integral $H_1$, its quotient by the closed image
has zero $\pi_1$ and is a point. The map is therefore onto.
The homotopy sequence then shows that the identity component of the
kernel has zero $\pi_1$ and that the whole kernel is connected.
The kernel is trivial, so the map is an algebraic isomorphism.

Serre's universal-morphism theorem
\cite[Theorem~7]{serre1959albanese} gives a based algebraic map
$a_V:V\to S_V$, universal for maps to semiabelian varieties.
Here his category consists of extensions of abelian varieties by
tori; translating the base value to zero turns his unique affine
homomorphism into a group homomorphism. Thus $\alpha_V$ factors as
$F\circ a_V$ for an algebraic group
homomorphism $F:S_V\to A_V$.
The map $a_{V*}$ is surjective on ordinary integral first homology
\cite[Theorem~10 and the following discussion, p.~19]{serre1959albanese}.
Since $H_1(S_V,\ZZ)$ is torsion-free, it induces a surjection
$u:H_{\mathrm f}\to H_1(S_V,\ZZ)$. The composite
$F_*u=\alpha_{V*}$ is the integral isomorphism proved above,
so $u$ and $F_*$ are both isomorphisms. By the algebraic homology
criterion just proved, $F$ is an isomorphism. Transporting
Serre's universal property along $F$ gives universality and uniqueness
for $\alpha_V$. Translation gives the unbased version; for projective
$V$ this is the ordinary Albanese. For the algebraic homology
criterion and its role in universality, see also
\cite[Lemmas~1--3 and Remark~4]{fujino2025correction}.

\emph{The maximal cover and the central quotient.}
Every connected finite etale cover $Z\to X$ is special by Campana's
etale-invariance theorem \cite[Proposition~10.11]{campana2009}; its
quasi-Albanese map is dominant by \cite[Lemma~4.4]{cadorel2026}.
If its target has abelian dimension $g$ and torus rank $t$, then
\[
 b_1(Z)=2g+t\le2(g+t)\le2\dim X
\]
by the integral construction just proved.
Choose a cover $Y$ with maximal first Betti number. The exact sequence
of fundamental groups for $a:Y\to A$ is
$1\to K\to G\to\Lambda_A:=\pi_1(A)\to1$
\cite[Proposition~4.12]{cadorel2026}. Compatible finite triangulation
\cite[Theorems~3.4.4 and~9.2.1]{bochnak1998}, applied to a compactification
and its deleted locus, gives a finite CW model for $F$: after
barycentric subdivision, normalize the weights on vertices outside
the deleted subcomplex to retract its complement to a finite
subcomplex. Hence $\pi_1(F)$ and $K$ are finitely generated.
Thus $\Tor K_{\rm ab}$ is finite;
$P$ is characteristic in $K$ and normal in $G$.

Put $V_K=K_{\rm ab}$ and $R_K=G/[K,K]$. Specialness implies the
weak specialness in Rogov's statement: finite etale covers and proper modifications
remain special, while logarithmic pullback excludes a dominant map
to a positive-dimensional variety of general type. By
\cite[Corollary~7.2]{rogov2025}, this finitely generated
solvable quotient of a special source is virtually nilpotent.
Classical results on finitely generated nilpotent groups give a
torsion-free subgroup of finite index and a faithful integral
unitriangular representation of that subgroup; see
\cite[Proposition~2.7(1)]{pan2025} and
\cite[Theorem~2.10]{dekimpe2017}, respectively.
Its inverse image in $G$ is the group of a special finite etale
cover of $Y$. The resulting complex linear representation has image
virtually nilpotent of class at most two by \cite[Theorem~B]{cao2026}.
Taking the normal core in $R_K$ gives a normal finite-index
torsion-free class-two subgroup $R_0$.

Let $V_0=R_0\cap V_K$, let $\Lambda_0$ be its image in $\Lambda_A$,
and put $r=b_1(Y)=\rk\Lambda_0$. Maximality gives
$r\le b_1(R_0)\le r$, so
$(R_0)_{\rm ab}\otimes\QQ\simeq\Lambda_0\otimes\QQ$.
Every $v\in V_0$ therefore has a positive power in
$C_0=[R_0,R_0]\subset Z(R_0)$; torsion-freeness and
$[v,h]^a=[v^a,h]$ give $V_0\subset Z(R_0)$.
Conjugation $\alpha_z$ by any $z\in R_K$ fixes $\Lambda_0$ and hence
the rational abelianization of $R_0$. Choosing $a,b>0$ with
$\alpha_z(u)^au^{-a},\alpha_z(w)^bw^{-b}\in C_0$ gives
\[
 \alpha_z([u,w])^{ab}
 =[\alpha_z(u)^a,\alpha_z(w)^b]=[u^a,w^b]=[u,w]^{ab}.
\]
Thus $\alpha_z$ fixes $C_0$ and its isolator $V_0$ pointwise.
Since $[V_K:V_0]<\infty$, its action on $V_K/\Tor V_K$ is trivial.
This proves the central extension. Its kernel and quotient are free
abelian, so $E$ is torsion-free of class at most two.

With $W_K=K/P$, the central extension gives
$\bigwedge^2\Lambda_A\to W_K$ with image $[E,E]$.
The surjections $G\twoheadrightarrow E\twoheadrightarrow\Lambda_A$
and equality of the outer first Betti numbers imply
$W_K/[E,E]$ finite. This pairing is therefore rationally surjective.
\end{proof}

\begin{lemma}\label{lem:canonical-kernel-unipotent-quotient}
Let $X$ be smooth connected special complex quasi-projective,
and let $Y,G,K,P,E$ be canonical kernel data as in
Definition~\ref{def:canonical-kernel}. Then every homomorphism from $G$
to a torsion-free nilpotent group has image of class at most two and
kills $P$. Moreover $E$ admits a faithful rational unitriangular
representation, and the canonical map from $G$ to its rational
two-step Malcev group has image canonically isomorphic to $E$.
\end{lemma}

\begin{proof}
The image $N$ of such a homomorphism is finitely generated and has a
faithful unitriangular representation \cite[Theorem~2.10]{dekimpe2017}.
By \cite[Theorem~B]{cao2026}, it has a class-two finite-index subgroup.
An algebraic subgroup of a characteristic-zero unipotent group is
connected, so this subgroup and $N$ have the same Zariski closure.
The polynomial identity $[x,[y,z]]=1$ extends to that closure;
hence $N$ itself has class at most two.

Every element of $K$ has torsion image in $G_{\rm ab}$ by the
quasi-Albanese homology isomorphism. Its image in $N$ has a power in
$[N,N]$, so is central by the class-two power identity and
torsion-freeness. Thus the image of $K$ is torsion-free abelian and
kills $P$. Conversely $E$ admits a faithful unitriangular representation.
Write $\gamma_1(G)=G$ and $\gamma_{j+1}(G)=[G,\gamma_j(G)]$ for the
lower central series.
The quotient $(G/\gamma_3G)/\Tor(G/\gamma_3G)$ consequently factors
through $E$, and $E$ factors through it; the quotient maps are inverse.
By the description of the Malcev completion of a nilpotent group
\cite[Appendix~A, Corollary~3.8]{quillen1969}, this identifies the image
of $G$ in its rational two-step Malcev group with $E$.
\end{proof}

\begin{proof}[Proof of part (ii) of Theorem~\ref{thm:canonical-albanese}]
\emph{The Hodge structure and integral lattice.}
Put $U_2=U_2(Y,y)$ and $\mathfrak u=\operatorname{Lie}U_2$.
Lemma~\ref{lem:canonical-kernel-unipotent-quotient}
identifies $\Gamma_2(Y,y)$ with $E$, and the pairing
in (i) identifies $\mathfrak d=[\mathfrak u,\mathfrak u]$
with $(K/P)\otimes\QQ$. By the functoriality of the mixed Hodge
structures constructed by Morgan
\cite[Theorem~9.2]{morgan1978} and Hain
\cite[Theorem~5.1]{hain1987geometry},
the inclusion $F\hookrightarrow Y$ induces a morphism of mixed Hodge
structures on Malcev Lie algebras. Its image is $\mathfrak d$, since
the image of $\pi_1(F)$ in $E$ is $K/P$. As $\mathfrak d$ is abelian,
this morphism factors through a mixed Hodge surjection
$H_1(F,\QQ)\twoheadrightarrow\mathfrak d$.
Its source has weights $-1,-2$ and Tate weight-$-2$ part, while
brackets in $\mathfrak u$ have weights at most $-2$. Strictness
\cite[Theorems~2.3.5 and~3.2.5]{deligne1971} implies that
$\mathfrak d$ is pure of type $(-1,-1)$. Hence $\mathfrak u$ has only weights $-1,-2$;
$\mathfrak c=W_{-2}\mathfrak u$ is central, pure Tate, and fits into
\[
 0\to\mathfrak c\to\mathfrak u
   \xrightarrow{\beta}H_1(B,\QQ)\to0,
 \qquad F^0\mathfrak c_{\CC}=0.
\]
Strictness identifies $F^0\mathfrak u_{\CC}$ with $F^0H_1(B,\CC)$.

Loops avoid the compactification boundary, so
$G\to\pi_1(\overline Y)\to\Lambda_B:=\pi_1(B)$ is surjective.
The homomorphism $A\to B$ has connected torus kernel: in the
first-Albanese Hodge realization its kernel is
$W_{-2}H_1(Y,\CC)/(\pi_1(A)\cap W_{-2}H_1(Y,\QQ))$,
a pure Tate vector space modulo its full integral lattice.
Writing $t=\dim W_{-2}H_1(Y,\QQ)$ gives
$0\to\ZZ^t\to\pi_1(A)\to\Lambda_B\to0$.
Consequently
\[
 J=\ker(E\to\Lambda_B)=E\cap\exp(\mathfrak c)(\QQ),
 \qquad 0\to K/P\to J\to\ZZ^t\to0.
\]
The group $J$ is central free abelian, $\log J$ is a full lattice in
$\mathfrak c_{\RR}$, and both $E/J$ and $J/(K/P)$ are torsion-free.

\emph{Algebraization.}
The second Albanese from \eqref{eq:based-second-albanese}, with
$\Gamma_2(Y,y)$ identified with $E$, is
$N=E\backslash U_2(\CC)/\exp(F^0\mathfrak u_{\CC})$. It
projects to $B^{\rm an}$ as a holomorphic torsor under
$T_J^{\rm an}=\mathfrak c_{\CC}/\log J$.
Indeed, lifting base representatives through $E\to\Lambda_B$ and
$F^0\mathfrak u_{\CC}\simeq F^0H_1(B,\CC)$ leaves a central
difference; the stabilizer is trivial since
$\Lambda_B\cap F^0H_1(B,\CC)=0$ and $F^0\mathfrak c_{\CC}=0$.
The projection is a submersion and therefore locally a torsor.

Central multiplication is a morphism of mixed Hodge groups preserving
$J\times E\to E$. Rogov's definability theorem for morphisms of
nil-Jacobians makes the action and projection definable
\cite[Definition~7 and Theorem~5.4(i)]{rogov2026}.
An integral basis of $J$ identifies the acting pure Tate torus with
the algebraic $(\CC^*)^{\rk J}$, compatibly with its definable
structure \cite[Theorem~6.3(iv)]{rogov2026}.
Definable transverse slices give local holomorphic sections.
Since $B$ is compact, finitely many such charts suffice. Their
transition functions define $\rk J$ definable line bundles.
Bakker--Mullane's proper definable GAGA identifies them with
analytic line bundles on the compact base $B$
\cite[Proposition~2.1, Corollary~2.2 and Lemmas~2.3--2.4]{bakker2022}.
Projective GAGA then algebraizes them to
$L_1,\ldots,L_{\rk J}$ on $B$
\cite[Section~12, Theorems~2--3]{serre1956}. Thus
$Q=\prod_B L_i^\times$ is smooth quasi-projective and definably
biholomorphic to $N$. The Hain--Zucker map $a_{2,y}$ and the projection
to level one are definable by Rogov's theorem
\cite[Theorem~6.3(i), (iii)]{rogov2026}.
Applying Peterzil--Starchenko's definable Chow theorem
\cite[Theorem~5.1]{peterzil2008} to their graphs gives algebraic maps
$f:Y\to Q$ and $p:Q\to A$, with $pf=a$; see also
\cite[Corollary~3.7(i)]{rogov2026} for this consequence.

\emph{The fundamental group.}
Choose a complement $L$ to $F^0H_1(B,\CC)$ and a splitting
$\sigma:H_1(B,\CC)\to\mathfrak u_{\CC}$ preserving $F^0$.
The map
\[
 (c,l)\longmapsto\exp(c)\exp(\sigma l)\exp(F^0\mathfrak u_{\CC})
\]
identifies $\mathfrak c_{\CC}\times L$ with
$U_2(\CC)/\exp(F^0\mathfrak u_{\CC})$: projection fixes $l$,
then the remaining central coordinate fixes $c$.
The free $E$-action is properly discontinuous, since its projected
translations form $\Lambda_B$ and transformations with a fixed
projection differ by $J$. This affine space is the universal cover,
so $Q^{\rm an}$ is $K(E,1)$. Continuation of the based higher Albanese
map along a loop gives its canonical Malcev lattice image. Thus
$f_*$ is the canonical map $G\twoheadrightarrow E$.
\end{proof}

\subsection{Logarithmic geometry of torus torsors}\label{subsec:torsor-logarithmic-geometry}

\begin{proposition}\label{prop:finite-cover-torus-torsor}
Let $q:Q\to B$ be a principal torus torsor over an abelian variety,
with fiber $T$, and let $j:Q'\to Q$ be a connected finite etale cover.
Put $M=\pi_1(T)$, $\Lambda=\pi_1(B)$ and $H=j_*\pi_1(Q')$.
Then $1\to M\to\pi_1(Q)\to\Lambda\to1$ is central.
There are an abelian isogeny $i:B'\to B$, a torus isogeny
$\phi:T'\to T$, and a principal $T'$-torsor $q':Q'\to B'$ such that
$qj=iq'$ and $j$ is $\phi$-equivariant. The maps on fundamental groups
identify $\pi_1(B')$ with the image $\Lambda'$ of $H$ in $\Lambda$
and $\pi_1(T')$ with $M'=H\cap M$.
\end{proposition}
\begin{proof}
The homotopy sequence of the torus bundle is exact on the left
because $\pi_2(B)=0$. Its fiber subgroup is central, since the two
factors of $\pi_1(T\times Q)$ commute under the action map.
Choose $d$ with $d\Lambda\subseteq\Lambda'$. Quotienting $B$ by
$(d^{-1}\Lambda')/\Lambda\subset B[d]$ and descending multiplication
by $d$ gives an abelian isogeny $B'\to B$ with image lattice
$\Lambda'$. The map $qj$ lifts to $q':Q'\to B'$.
Such a lift is algebraic: its graph selects a degree-one component of
the finite etale pullback cover, by Riemann existence
\cite[Expos\'e~XII, Theorem~5.1]{grothendieck2004}.

Set $\widetilde Q=Q\times_BB'$. The finite etale map
$\widetilde j=(j,q'):Q'\to\widetilde Q$ is surjective and
$\pi_1(\widetilde Q)=HM$. Since $M$ is central, $H$ is normal in $HM$.
The components over a torus fiber correspond to
$H\backslash HM/M$, a singleton. Hence every $q'$-fiber is connected
and is the finite torus cover with lattice $M'$. Let
$\phi:T'\to T$ be the corresponding torus isogeny.

The map $T'\times Q'\to\widetilde Q$ given by
$(t',z)\mapsto\phi(t')\widetilde j(z)$ has fundamental-group image
in $H$, so it has a unique based lift to $Q'$. This lift is algebraic by the
preceding Riemann-existence argument. Uniqueness of
based lifts proves the identity and associativity laws. On each fiber,
the orbit map identifies two connected covers of $T$ with the same
lattice $M'$, hence is an isomorphism. Thus
$T'\times Q'\to Q'\times_{B'}Q'$ is an isomorphism (it is bijective
and etale), and the smooth surjection $q'$ is a principal $T'$-torsor.
A character basis of the split torus expresses it as a product of
frame bundles of line bundles.
\end{proof}

\begin{definition}\label{def:canonical-kernel-logarithmic-pair}
Let $V$ be smooth connected quasi-projective and let $\Gamma$ be an
effective rational SNC divisor on $V$ with coefficients in $[0,1]$.
Choose a smooth projective compactification $\overline V$ such that
the full reduced boundary $D_V=\overline V\setminus V$ and the closure
$\overline\Gamma$ have combined SNC support. We write
\[
 \kappa(V,\Gamma)=
 \kappa(\overline V,K_{\overline V}+D_V+\overline\Gamma),
 \qquad \overline\kappa(V)=\kappa(V,0).
\]
A proper smooth resolution of a normal quasi-projective variety retains
all divisors over that variety. Its deleted boundary consists only of
the complement in a projective compactification.
\end{definition}

Components of $\Gamma$ belong to the open variety, including those
with coefficient one.

\begin{lemma}\label{lem:canonical-kernel-logarithmic-comparison}
The number in Definition~\ref{def:canonical-kernel-logarithmic-pair}
is independent of the chosen compactification. Proper birational maps
between smooth quasi-projective varieties preserve logarithmic Kodaira
dimension and specialness. Specialness descends under dominant
morphisms between smooth quasi-projective varieties.

Let $Q\to B$ be a principal split-torus torsor over an abelian variety.
If a smooth connected quasi-projective variety $H$ has a morphism to
$Q$ generically finite onto its image, then $\overline\kappa(H)\ge0$.
Let $A_0$ be an abelian variety, let $W$ be smooth connected
quasi-projective, and let $p:E_0\to W$ be the composite of a principal
split-torus torsor $E_0\to A_0\times W$ with projection. For every
effective rational SNC divisor $\Theta$ on $W$ with coefficients in
$[0,1]$, one has
\[
 \kappa(E_0,p^*\Theta)=\kappa(W,\Theta).
\]
After fixing an invariant volume on the abelian and torus fibers,
the compactifications below, with their full deleted boundaries, have
isomorphisms compatible with multiplication:
\[
 H^0\!\left(\overline E_0,mk(K_{\overline E_0}+D_{E_0}
                    +\overline{p^*\Theta})\right)
 \simeq
 H^0\!\left(\overline W,mk(K_{\overline W}+D_W
                    +\overline\Theta)\right)
 \quad(m\ge0)
\]
for sufficiently divisible $k$. These identifications persist on
compatible compactifications. Via $p^*$ they identify the full fields
generated by ratios of sections.
\end{lemma}

\begin{proof}
\emph{Model comparison.}
At a blowup center of codimension $c$, contained in boundary components
of coefficients $\beta_j$, the exceptional correction to the adjoint
is $c-1+\delta-\sum\beta_j$, where $\delta$ is the new coefficient.
It is nonnegative at infinity, where $\delta=1$, and for retained
interior blowups of an unweighted open, where it is $c-1$.
At arbitrary valuations the same bound follows from
$x_j=t^{a_j}u_j$ and the Jacobian order $\sum a_j-1$.
More precisely, at a deleted exceptional prime the correction is at
least $\sum_j(1-\beta_j)a_j\ge0$, where unmarked center parameters
have coefficient zero. Thus, for compatible canonical divisors, the
new adjoint is $\mu^*L+E$ with $E\ge0$ exceptional. After clearing
denominators, the projection formula gives
\[
 \mu_*\mathcal O(m\mu^*L+mE)=\mathcal O(mL).
\]
The sheaf of an effective exceptional divisor pushes forward to
$\mathcal O$: normality extends rational functions with no
codimension-one poles. Common resolutions therefore identify the
divisible section rings on compactifications of the same weighted
open and on proper smooth modifications of an unweighted open.

Logarithmic pullback followed by saturation preserves independent
section ratios. Conversely, a line on a proper smooth modification
restricts to a generic line downstairs, whose saturation and sections
extend across codimension two. Indeed its saturation in the locally
free logarithmic differential sheaf is reflexive of rank one, hence
a line bundle on the smooth target: its double dual embeds in that
sheaf, and saturation excludes a nonzero torsion quotient. At every
target prime the birational map is an isomorphism generically, so
each upstairs section has no divisorial pole downstairs. Normality
then extends the sections without changing their ratios, proving
birational invariance of specialness. The pullback argument also proves
its descent under a dominant morphism: every preimage of deleted target boundary is
deleted in the source, and locally
$d\log(t^au)=a\,dt/t+d\log u$.

For a dominant generically finite map of smooth opens, logarithmic
pullback induces a nonzero map of logarithmic canonical bundles on compatible
compactifications, since the differential has full rank at the generic point. Its tensor powers
inject all pluricanonical section spaces and preserve their ratios.
For a connected fibration, a nonzero logarithmic pluricanonical
section on the total space restricts nontrivially to a whole
sufficiently general resolved fiber. To check this, compactify
the map properly, make the source boundary SNC, and avoid the images
of vertical strata and the zero divisor's vertical components.
Adjunction on a smooth compact fiber identifies the restricted
adjoint with its logarithmic canonical bundle, since the normal
bundle is constant and the restricted boundary is exactly the
deleted boundary. A horizontal component of the section's zero
divisor cannot contain that whole fiber. Thus the fiber has
nonnegative logarithmic Kodaira dimension.

\emph{Torsors.}
Write $Q=\prod_B L_i^\times$ and compactify it by
$\overline Q=\prod_B\mathbb P(\mathcal O_B\oplus L_i)$,
with toric boundary $J_Q$. Local coordinates give
\[
 0\to\overline q^*\Omega_B^1
 \to\Omega^1_{\overline Q}(\log J_Q)
 \to\mathcal O_{\overline Q}^{\dim T}\to0.
\]
Both outer terms are trivial. Resolve the graph of the map from $H$
outside $H$ to extend it to a compactification. Put $a=\dim H$.
The $a$th exterior power of logarithmic pullback is a nonzero map to
$\omega_{\overline H}(D_H)$ from a bundle filtered by trivial
quotients. Its first nonzero filtration term yields a nonzero
logarithmic top form. Thus $\overline\kappa(H)\ge0$.

For $E_0$, extend its character line bundles to
$A_0\times\overline W$ by closure of Cartier divisors. Their
projective-line bundles give $\overline p:\overline E_0\to\overline W$
with connected integral fibers and full boundary
$D_{E_0}=J_0+\overline p^*D_W$. Wedge an invariant top form on
$A_0$ with the relative $d\log$ coordinates. Coordinate changes
add differentials from $A_0\times\overline W$; the $A_0$ terms
disappear on wedging with its invariant volume, and the
$\overline W$ terms disappear in the relative determinant. Hence
\[
 \omega_{\overline E_0}(D_{E_0})
   \simeq\overline p^*\omega_{\overline W}(D_W),
 \qquad\overline p_*\mathcal O_{\overline E_0}=\mathcal O_{\overline W}.
\]
Since $\overline{p^*\Theta}=\overline p^*\overline\Theta$, the
projection formula gives the stated multiplicative isomorphisms
in all divisible degrees. Model comparison extends these isomorphisms to
other compactifications. The corresponding powers of the same fiber
volume cancel in equal-degree ratios, identifying the full fields
of section ratios by pullback.
\end{proof}

\begin{lemma}\label{lem:canonical-kernel-divisor-positivity}
Let $Q\to B$ be a principal split-torus torsor over an abelian variety,
and let $D^\circ$ be a prime divisor of $Q$. There is a smooth
projective equivariant compactification $\overline Q$, with full
reduced SNC boundary $J_Q$, such that
$K_{\overline Q}+J_Q\sim0$, the closure $D$ contains no boundary
stratum, and $\kappa(\overline Q,D)\ge1$.
\end{lemma}

\begin{proof}
Start with the product of projective-line bundles used in
Lemma~\ref{lem:canonical-kernel-logarithmic-comparison}.
Its nonempty boundary intersections are connected toric bundles.
The combinatorial desingularization of
\cite[Corollary~2 and Section~5.1]{molinasamper2018}, applied to the
boundary orders of the local equations of $D^\circ$, gives a finite
sequence of blowups of whole boundary strata after which the closure
contains no stratum. The geometric application is described in
\cite[Section~3]{fujino2024a}. It applies here in the toric fiber
charts: changing a character-bundle trivialization multiplies each
toric coordinate by a unit, preserving these boundary orders and
the strata chosen as centers. Thus the blowups glue over $B$.
They are equivariant and preserve the logarithmic volume.

If $D$ is not invariant under the torus, a translate distinct from $D$
is linearly equivalent to $D$ by
\cite[Theorem~2.14 and Lemmas~2.12--2.13]{brion2014}:
the acting group is a split torus, for which the stable exponent
in that result is one. Their section ratio is
nonconstant, proving positive Iitaka dimension. If $D$ is invariant,
it is the pullback of a prime $D_B$ on $B$. On a trivializing chart
this follows because an invariant ideal in a Laurent polynomial ring
is homogeneous and all monomials are units. The theorem of the square
\cite[Theorem~6.7 and Remark~6.9]{milne1986} gives
$t_b^*D_B+t_{-b}^*D_B\sim2D_B$.
Choose a translate of $D_B$ different from $D_B$.
For such a translate the displayed effective divisor differs from
$2D_B$, giving a nonconstant section ratio. Pulling it back proves
the assertion.
\end{proof}

\begin{lemma}\label{lem:canonical-kernel-zero-covers}
Let $Q\to B$ be a principal split-torus torsor over an abelian variety.
Then the following hold.
\begin{enumerate}[label=(\roman*)]
\item If $v:V\to Q$ is finite surjective with $V$ connected and normal,
and a proper smooth resolution of $V$ has logarithmic Kodaira dimension
zero, then $v$ is etale.
\item If $a:H\to Q$ is a birational morphism with $H$ smooth connected
quasi-projective and $\overline\kappa(H)=0$, then there is a closed
$Z\subset Q$ of codimension at least two such that $a$ is an isomorphism
over $Q\setminus Z$. Moreover $a_*:\pi_1(H)\to\pi_1(Q)$ is an
isomorphism.
\end{enumerate}
\end{lemma}

\begin{proof}
\smallskip\noindent\emph{The logarithmic valuation calculation.}
Choose the compactification in
Lemma~\ref{lem:canonical-kernel-divisor-positivity} for a specified prime
$D^\circ$, and choose a logarithmic volume $\omega$ with divisor
$-J_Q$. For a generically finite morphism
$g:\overline H\to\overline Q$, and a prime $E_H$ whose center has
codimension $c$ and lies in $k$ boundary components, write
$g^*x_j=t^{a_j}u_j$ for center coordinates, $1\le j\le c$.
In the top wedge at most one factor contributes $dt/t$, so
\begin{equation}\label{eq:canonical-kernel-valuation}
 \operatorname{ord}_{E_H}(g^*\omega)
       \ge\sum_{j=k+1}^{c}a_j-1.
\end{equation}
Adding coefficient one at a deleted prime gives a nonnegative
coefficient, strictly positive when $c>k$. At an interior prime the
coefficient is at least $c-1$, and is strictly positive for $c\ge2$.
At a nonexceptional interior prime it is exactly the ramification
index minus one. Therefore a finite Galois cover ramified at
$D^\circ$ has, on its proper resolution with full boundary,
\[
 K_{\overline H}+D_H\ge\epsilon g^*D
\]
for some rational $\epsilon>0$. All sheets over $D^\circ$ have the
same ramification index, and a center in $D\cap J_Q$ cannot be an
entire boundary stratum. Pullback of sections gives logarithmic Kodaira
dimension at least one. Thus a Galois cover of logarithmic Kodaira
dimension zero is unramified along divisors. Zariski--Nagata purity
\cite[Lemma~58.21.4, Tag~0BMB]{stacks2026} then implies that the cover is etale.

For a birational $a$ as in (ii), suppose no prime of $H$ dominates
$D^\circ$. Compactify the morphism; its strict transform of
$D$ is deleted and has coefficient one in $K_{\overline H}+D_H$.
Formula~\eqref{eq:canonical-kernel-valuation} gives the same effective
inequality at all other components of $g^*D$. This again contradicts
$\overline\kappa(H)=0$. Every target prime is therefore dominated by
a source prime. The two corresponding discrete valuation rings agree
because the map is birational and both varieties are smooth. Near
their generic points the map is an isomorphism. The complement of
the largest such target open has codimension at least two; take its
closure as $Z$. The common open surjects on $\pi_1(H)$, since loops
avoid a proper complex subvariety, and its map to $Q$ is an isomorphism
on fundamental groups by transversality in codimension two.

\smallskip\noindent\emph{Reduction to one torus coordinate.}
For rank zero, $V$ is complete because it is finite over $B$.
The Kawamata--Viehweg characterization theorem
\cite[Main Theorem]{kawamata1980}, in the finite-cover formulation of
\cite[Theorem~9.1]{fujino2024b}, implies that $V\to B$ is etale.
For positive rank write $Q\to Q_1$ as the last $\mathbb G_m$-bundle
in a product of character bundles. Normalize its projective-line
compactification over $Q_1$ in $\CC(V)$ and take the proper Stein
factor. Its base $V_1$ is finite normal over $Q_1$; the generic
connected fiber is a smooth projective curve over $\CC(V_1)$, finite
over $\mathbb P^1$, and its open part is the inverse image of
$\mathbb G_m$. On a general complex fiber the punctures include
the preimages of zero and infinity. Riemann--Hurwitz gives
\[
 \deg(K_{\overline C}+D_C)
    =\sum_{x\notin D_C}(e_x-1)\ge0.
\]
Resolve the source and base while retaining their open parts.
\cite[Theorem~1.1]{fujino2005}, together with the
nonnegativity assertion of
Lemma~\ref{lem:canonical-kernel-logarithmic-comparison}, gives
\[
 0\ge\overline\kappa(V_1)+\overline\kappa(C)\ge0.
\]
By induction, $V_1\to Q_1$ is etale.
Proposition~\ref{prop:finite-cover-torus-torsor} identifies $V_1$ as a
torus torsor over the source of an abelian isogeny $B_1\to B$. Hence the pullback
$Q_2=Q\times_{Q_1}V_1$ is also a torus torsor over $B_1$, because
the last character bundle is pulled back from $B_1$.

The zero logarithmic degree of $C$ means that its geometric completion
has genus zero and exactly two punctures. Each lies over a different
rational endpoint and is therefore rational over $\CC(V_1)$.
With a coordinate $u$ having these as zero and pole, the original
coordinate is $t=\lambda u^d$ for $\lambda\in\CC(V_1)^*$.
The residual extension $\CC(V)/\CC(Q_2)$ is cyclic, since this field
contains all $d$th roots of unity. The finite morphism
$V\to Q_2$ is thus Galois on function fields and has the same zero
logarithmic Kodaira dimension. It is etale by the preceding valuation argument
and purity. Composing with the etale map $Q_2\to Q$ proves
(i) for every rank.
\end{proof}

\subsection{Special images and families}\label{subsec:special-images-families}

\begin{lemma}\label{lem:canonical-kernel-zero-images}
Let $q:Q\to B$ be a principal torus torsor with torus $T$ over an
abelian variety. Let $H$ be smooth connected quasi-projective with
$\overline\kappa(H)=0$, and let $h:H\to Q$ be generically finite onto
its reduced closed image closure $S$. Then there are a connected
subtorus $U\subseteq T$, an abelian subvariety $B_0\subseteq B$,
a point $b\in B$, and an abelian variety $C$ finite etale over
$b+B_0$, with a finite morphism $C\to Q/U$, such that
\[
 S^\nu=Q\times_{Q/U}C
\]
as varieties over $Q$. The induced map
$h_*:\pi_1(H)\to\pi_1(Q)$ is injective. If $h$ is proper, its finite
normal factor over $S^\nu$ is etale and its degree is the index of
this subgroup in $\pi_1(S^\nu)$.
\end{lemma}

\begin{proof}
\smallskip\noindent\emph{Closed images and rank induction.}
A proper smooth resolution $X_S\to S$ has logarithmic Kodaira
dimension zero: resolve the rational lift from $H$, pull back
logarithmic pluriforms, and use the nonnegativity in
Lemma~\ref{lem:canonical-kernel-logarithmic-comparison}.

We use Kawamata's quasi-Albanese theorem in logarithmic dimension zero
\cite[Theorem~28 and Corollary~29]{kawamata1981}.
For a smooth quasi-projective $V$ with $\overline\kappa(V)=0$,
its quasi-Albanese map $\alpha_V:V\to A_V$ is dominant with
irreducible general fibers, and
\[
 \dim A_V\le\dim V,
 \qquad \dim A_V=\dim V\ \Longleftrightarrow\
 \alpha_V\text{ is birational}.
\]
Kawamata's term \emph{open algebraic fiber space} means dominance
and irreducibility of the general fibers. In characteristic zero,
generic smoothness makes the whole general fiber geometrically
integral. Based universality identifies his quasi-Albanese with the
one constructed in Theorem~\ref{thm:canonical-albanese}(i).

We prove the normalization assertion for such a closed $S$ by induction
on torus rank. In rank zero, invariant forms on $B$ pull back to at least
$\dim S$ independent logarithmic one-forms on $X_S$.
Here $X_S$ is projective, so its quasi-Albanese is its ordinary
Albanese. The preceding argument shows that its Albanese map is birational and
that the target has dimension $\dim S$. The factor homomorphism from
that abelian variety to $B$ has closed image equal, after translation,
to $S$, and is finite onto it by equality of dimensions.
The degree of $X_S\to S$ and of the Albanese map are both one;
the homomorphism therefore has degree one and trivial kernel.
It is an isomorphism onto that abelian translate, giving the
asserted normalization.

In positive rank choose a coordinate quotient $\rho:Q\to Q_1$
with kernel $\mathbb G_m$, and let $S_1$ be the reduced closure of
$\rho(S)$. If $\dim S=\dim S_1+1$, the closed integral $S$ equals
$\rho^{-1}(S_1)$ by dimension. Its generic fiber is $\mathbb G_m$.
The addition theorem for curves and nonnegativity show that a proper
resolution of $S_1$ has logarithmic Kodaira dimension zero. Induction classifies
$S_1^\nu$; pulling it back along $\rho$ gives the normalization of
$S$ and the asserted subtorus, the inverse image of the subtorus
obtained in rank one less.

Suppose instead that $S\to S_1$ is generically finite. Pullback of
logarithmic sections again shows that a proper resolution of $S_1$
has logarithmic Kodaira dimension zero. By induction, $S_1^\nu$ is a
torus torsor over an abelian variety. Normalize $S_1^\nu$ in $\CC(S)$ and
call the resulting finite normal variety $V_1$. Normality gives a lift
$X_S\to V_1$. A proper resolution of $V_1$ has logarithmic Kodaira dimension
zero by pullback and nonnegativity. By Lemma~\ref{lem:canonical-kernel-zero-covers}(i),
$V_1\to S_1^\nu$ is etale. Proposition~\ref{prop:finite-cover-torus-torsor}
then identifies $V_1$ as a torus torsor over an abelian variety $C$
finite etale over the abelian base of $S_1^\nu$.

The birational morphism $X_S\to V_1$ is an isomorphism outside a
target subset of codimension at least two, by
Lemma~\ref{lem:canonical-kernel-zero-covers}(ii). On that open its
map to $Q$ is a nowhere-zero section of the pullback of the last
character line bundle. The section and its inverse extend over
$V_1$ by normality, giving a morphism $V_1\to Q$.
It is finite, because $V_1\to Q_1$ is finite and its graph over the
separated $Q\to Q_1$ is closed. Its image is $S$, and it is
birational onto $S$, so $V_1=S^\nu$.

On a local trivialization over $C$, a nowhere-zero function on a
torus is a constant times a character. The character is constant
on the connected base. The map $V_1\to Q$ is therefore equivariant
for a homomorphism from its acting torus to $T$, lifting the known
isogeny onto the lower-rank subtorus. Its kernel is finite and is
trivial because $V_1\to S$ has generic degree one. Its image is a
connected subtorus $U\subset T$. Quotienting by $U$ gives the map
$C\to Q/U$. This map is finite, since $C\to b+B_0$ is finite and
$Q/U$ is separated over $B$. The equivariant map identifies
$V_1$ with $Q\times_{Q/U}C$ fiber by fiber. In both cases
$q(S)=b+B_0$, completing the induction.

\smallskip\noindent\emph{The fundamental group of the source.}
For the original $H$, let $V_H$ be the normalization of $S^\nu$ in
$\CC(H)$. The morphism lifts to $V_H$ by normality.
A proper resolution of $V_H$ has logarithmic
Kodaira dimension zero, so its finite map to the torsor $S^\nu$
is etale by Lemma~\ref{lem:canonical-kernel-zero-covers}(i).
The variety $V_H$ is a torus torsor by
Proposition~\ref{prop:finite-cover-torus-torsor}.
Lemma~\ref{lem:canonical-kernel-zero-covers}(ii) identifies
$\pi_1(H)$ with $\pi_1(V_H)$.
The map $\pi_1(S^\nu)\to\pi_1(Q)$ is injective: in the two torus
bundle sequences the left map is the inclusion of cocharacter
lattices of $U\subseteq T$ and the right map is the injection of
the finite etale abelian cover $C\to b+B_0\subseteq B$.
This proves injectivity for $h_*$ and the stated covering index.
When $h$ is proper, $H\to V_H$ is proper birational, so this is
also its finite factor in the proper Stein factorization.
\end{proof}

\begin{lemma}\label{lem:canonical-kernel-family-splitting}
Let $q:Q\to B$ be a principal torus torsor with torus $T$ over an
abelian variety. Fix a connected subtorus $U\subseteq T$ and an
abelian subvariety $B_0\subseteq B$. Let $\mathcal I$ and $H_0$ be
integral complex varieties of finite type with dominant morphisms
$\pi:\mathcal I\to H_0$ and $e:\mathcal I\to Q$. Suppose that, for
every $h$ outside a countable union of proper closed subsets of $H_0$,
the whole image $e(\pi^{-1}(h))$ is a closed irreducible subvariety
$S_h$ with a proper smooth resolution of logarithmic Kodaira dimension
zero, that $U$ preserves $S_h$, that $q(S_h)$ is a translate of $B_0$,
and that $\dim S_h=\dim U+\dim B_0$.
Then there is an isogeny $\alpha:B_0\times B_1\to B$ such that the
pullback $Q'=Q\times_B(B_0\times B_1)$ admits a smooth
surjection $r:Q'\to Q_1$, where $Q_1$ is a principal $(T/U)$-torsor
over $B_1$ and $Q'\to B_0\times Q_1$ is a principal $U$-torsor.
The reduced inverse image of each such $S_h$ is a finite disjoint
union of whole fibers of $r$.
\end{lemma}

\begin{proof}
Let $L_\chi$ be a character line bundle of $Q/U\to B$.
Lemma~\ref{lem:canonical-kernel-zero-images} identifies the normalization
of $S_h$ with a torus torsor over an abelian variety $C_h$ finite
etale over $b_h+B_0$. Its stabilizing torus equals $U$: the product
of it with $U$ acts freely in a fiber, whose dimension is $\dim U$.
The map $C_h\to Q/U$ trivializes $L_\chi$ after this finite etale
pullback. If its degree is $d_h$, the projection formula gives
$(L_\chi|_{b_h+B_0})^{\otimes d_h}\simeq\mathcal O$:
the norm is
$\det(\lambda_*\lambda^*L_\chi)\otimes
\det(\lambda_*\mathcal O)^{-1}=L_\chi^{\otimes d_h}$.

Consider the restriction classes
$\Phi_\chi(b)=[L_\chi|_{b+B_0}]$ in $\operatorname{Pic}^0(B_0)$,
using translation to identify the coset with $B_0$.
They all have zero Chern class, since one is torsion and translation
is homotopic to the identity. This is a holomorphic map of complex
tori: locally in $b$, holomorphic transition functions and their
logarithms give its class in
$H^1(B_0,\mathcal O)/H^1(B_0,\ZZ)$ via the exponential sequence.
By translation invariance of degree-zero line bundles
\cite[Proposition~9.2]{milne1986}, this class is constant along each
coset $b_h+B_0$ arising from the family.
If $\Phi_\chi$ were nonconstant, each inverse image of a torsion point
under $\Phi_\chi q e$ would be a proper closed analytic subset of
$\mathcal I$. There are countably many torsion points. Together
with the inverse images of the bad parameter sets these subsets
would cover $\mathcal I$, contrary to the Baire theorem on a smooth
analytic coordinate neighborhood. Dominance of both maps makes all
these subsets proper. Thus $\Phi_\chi$ is constant with a torsion
value $P_\chi$.

Poincare complete reducibility
\cite[Proposition~12.1]{milne1986} gives an abelian complement $B_1$
such that the addition map $\alpha_0:B_0\times B_1\to B$ is an isogeny.
On $B_0\times B_1$, the
pullback of $L_\chi$, tensored with the inverse of $P_\chi$ from
$B_0$, is trivial on every $B_0$-fiber.
The seesaw theorem \cite[Theorem~5.1]{milne1986} gives
\[
 \alpha_0^*L_\chi\simeq
 \operatorname{pr}_0^*P_\chi\otimes\operatorname{pr}_1^*M_\chi,
 \qquad \alpha_0(a,c)=a+c.
\]
Choose one integer $n$ killing the finitely many torsion bundles for
a character basis. Since $P_\chi\in\operatorname{Pic}^0(B_0)$,
translation invariance and the same seesaw theorem give the identity
$m^*P_\chi\simeq\operatorname{pr}_1^*P_\chi\otimes
\operatorname{pr}_2^*P_\chi$ for the addition map $m$ of $B_0$.
It follows that $[n]^*P_\chi=P_\chi^{\otimes n}$.
Under $\alpha(a,c)=na+c$, all quotient character bundles therefore
pull back from $B_1$. Their frame bundles define $Q_1$, and give the
torsor $Q'\to B_0\times Q_1$.

The quotient image of $S_h$ in $Q/U$ is the proper image of $C_h$.
After the isogeny, each irreducible component of this quotient image
lies over one point of $B_1$. Its projection to the corresponding
affine torus fiber of $Q_1$ is constant, because the component is
proper. Its dimension is $\dim B_0$, so it is exactly
$B_0\times\{z\}$ for that constant $z\in Q_1$.
Distinct components give distinct constants. Pulling back along
the principal $U$-torsor gives the claimed disjoint union of whole
$r$-fibers.
\end{proof}

\begin{lemma}\label{lem:canonical-kernel-special-proper-image}
Let $Q\to B$ be a principal torus torsor over an abelian variety.
Let $Z$ be smooth connected special quasi-projective and let
$g:Z\to Q$ be a projective morphism generically finite onto its
image. Then $\overline\kappa(Z)=0$.
\end{lemma}

\begin{proof}
\smallskip\noindent\emph{The full leaves and their period subgroups.}
Nonnegativity follows from
Lemma~\ref{lem:canonical-kernel-logarithmic-comparison}. Put
$k=\overline\kappa(Z)$ and suppose $k>0$.
The quasi-Albanese of $Z$ is dominant by
\cite[Lemma~4.4]{cadorel2026}; universality shows that
the closure of $qg(Z)$ is an abelian translate $B^\flat$.
Work in the torsor $Q^\flat=q^{-1}(B^\flat)$, identifying this
translate with an abelian variety by a point.

Choose a smooth projective compactification $\overline Z$ of $Z$
with full SNC boundary $D_Z$, and resolve its full logarithmic
Iitaka map by $\mu:\overline Z_m\to\overline Z$.
Write $Z_m=\mu^{-1}(Z)$ and
$h:\overline Z_m\to H$ for the map with connected fibers, where $H$ is
smooth projective of dimension $k$. The boundary
$D_m=(\overline Z_m\setminus Z_m)_{\rm red}$ is SNC;
every divisor over $Z$ is retained.
For a very general $w\in H$, the whole open leaf
$F_w=h^{-1}(w)\cap Z_m$ has logarithmic Kodaira dimension zero.
Indeed the compactification change gives
\[
 K_{\overline Z_m}+D_m
   =\mu^*(K_{\overline Z}+D_Z)+E,\qquad E\ge0
   \text{ exceptional},
\]
by Lemma~\ref{lem:canonical-kernel-logarithmic-comparison}.
The standard Iitaka-fiber statement in
\cite[Notation~4.1 and Lemma~4.2]{taji2021} says that the pullback
of a line bundle defining its full Iitaka fibration, plus an
effective resolution-exceptional divisor, has Iitaka dimension zero
on a very general compact fiber. Apply it to
$K_{\overline Z}+D_Z$ and $E$. Adjunction identifies the restriction
with the canonical divisor of the compact fiber plus its full
deleted boundary, which is the logarithmic adjoint of $F_w$.

The leaf maps to $Q^\flat$ are projective and generically finite onto
their closed images $S_w$, of dimension $\dim Z-k$. By
Lemma~\ref{lem:canonical-kernel-zero-images}, their normal images
are torus torsors and their finite normal factors are etale.
Put $G^\flat=\pi_1(Q^\flat)$, and let
$P_w\le N_w\le G^\flat$ be the images of $\pi_1(F_w)$ and
$\pi_1(S_w^\nu)$. Both maps are injective, and $[N_w:P_w]$ is the degree of the finite
normal factor.

These pairs of subgroups are conjugate as $w$ varies in a connected
dense parameter open. To check this for the normal images as well
as the leaves, take the reduced image of
$(h,g):Z_m\to H\times Q^\flat$, and normalize it both in its own
function field and in that of $Z_m$. After restricting to a parameter open, generic flatness and geometric
normality identify the fibers of these two finite models with
$S_w^\nu$ and its finite normal cover. The birational
maps are isomorphisms on opens meeting every such fiber; hence they
remain birational on the integral fibers. Shrink again so the normal
families are smooth and the finite map is etale. This is possible
because their bad loci have constructible images missing very general
parameters. Thom's first isotopy lemma and compatible Whitney
stratifications \cite[Sections~1.5 and~1.7]{goresky1988} give the local
triviality statement in \cite[Theorem~2.2]{achinger2019}.
Apply this statement to this single family of
normalized images. The theorem gives a finite stratification of the base on
which the whole inverse-image families are locally topologically
trivial. Choose a smooth connected dense stratum. Over contractible analytic
neighborhoods in this stratum, the whole
fibers of the family of normalized images are therefore trivialized.

Over a contractible such neighborhood, lift the resulting paths
through the one finite etale normal cover. Unique path lifting
gives compatible trivializations of that cover and all its
connected fibers.
If $\Phi_t$ is the transport in the family of normalized images and
$\gamma$ a loop in its starting fiber, then
$(s,t)\mapsto e(\Phi_t(\gamma(s)))$ is a homotopy in $Q^\flat$.
The moving base point follows the single path
$t\mapsto e(\Phi_t(\gamma(0)))$. The same path is used for the
cover, because its evaluation factors through this family.
Consequently the two subgroup images are conjugated by the same
element of $G^\flat$. Chaining such neighborhoods proves
conjugacy of the pairs $(N_w,P_w)$ throughout the connected
parameter open. The proper birational leaf maps give exactly the
images of these finite normal covers by
Lemma~\ref{lem:canonical-kernel-zero-images}.

\smallskip\noindent\emph{Birational leaf maps after a finite cover.}
By Proposition~\ref{prop:finite-cover-torus-torsor}, $G^\flat$ is
finitely generated and nilpotent. We use Mal'cev's classical
subgroup-separability theorem in the following form: if $G_*$ is
finitely generated nilpotent,
$H_*\le G_*$, and $x\notin H_*$, there is a finite-index subgroup
containing $H_*$ and omitting $x$.
This also follows directly from
\cite[Theorem~3.1]{wong2009separability}. In its notation take
both ambient groups equal to $G_*$, put $K_1=H_*$,
$K_2=\{1\}$ and $y=1$. The theorem gives a normal finite-index
$J_x\triangleleft G_*$ with $x\notin H_*J_x$; the subgroup
$H_*J_x$ has the required properties.

Apply this argument with $G_*=G^\flat$ to one transported pair
$P\le N$. For each of the finitely many nontrivial $P$-cosets
in $N$, choose a finite-index subgroup containing $P$ and
excluding its representative. Their intersection meets $N$ in exactly $P$.
Its normal core $J$ has finite index and satisfies
$J\cap N_w\subseteq P_w$ for every transported pair.
Normality also gives
\[
 J\cap N_w=n(J\cap N_w)n^{-1}\subseteq nP_wn^{-1}
 \qquad(n\in N_w).
\]

After pullback to the connected finite etale cover $Q^e\to Q^\flat$
corresponding to $J$, every component of every very general lifted leaf
has degree one onto its normal image. The remaining covering degree is
$[J\cap N_w:J\cap nP_wn^{-1}]=1$.
The cover is algebraic by \cite[Theorems~3.4 and 21.3]{milne2013}.
Proposition~\ref{prop:finite-cover-torus-torsor} then identifies
$Q^e$ as a torus torsor over an abelian variety $B^e$.

Choose a source component over this cover. It is special by
\cite[Proposition~10.11]{campana2009} and has logarithmic
Kodaira dimension at least $k$ by pullback. Take the proper Stein
factor of its compact map to $H$. Its base has dimension $k$ and
its full open fibers are the connected lifted leaves. Their maps
are proper birational onto their smooth normal images.
Their period subgroups determine a fixed connected subtorus
$U\subseteq T^e$ and a fixed abelian subvariety $B_0\subseteq B^e$:
intersection with the central torus lattice fixes $U$, and the
rational span of the projection to $\pi_1(B^e)$ fixes $B_0$.
These directions satisfy $\dim U+\dim B_0=\dim Z-k$.

Form the incidence by taking the product of this full leaf family
with $T^e$ and evaluating by torus translation in $Q^e$.
Both maps are dominant. For evaluation, the source projects densely
to $B^e$, and the torus fills each fiber over this dense image.
Lemma~\ref{lem:canonical-kernel-family-splitting} gives one further
isogeny pullback $Q'\to Q^e$ and a smooth map $r:Q'\to Q_1$, with
$Q'\to B_0\times Q_1$ a principal $U$-torsor. The very general
connected lifted leaves are now proper birational onto whole
$r$-fibers. Choose a corresponding source component $Z'$; it
remains special and has logarithmic Kodaira dimension at least $k$.

\smallskip\noindent\emph{Descent to the quotient.}
Let $S'$ be its closed image in $Q'$ and let $W$ be the closure of
$r(S')$ in $Q_1$, both reduced. Dense very general leaf images are
whole $r$-fibers contained in $S'$. Fiber dimension therefore gives
$\dim W=k$. The integral variety $r^{-1}(W)$ has dimension
$\dim W+\dim U+\dim B_0=\dim Z$, and contains the closed integral
$S'$ of that dimension. Hence $S'=r^{-1}(W)$.

Let $V'$ be the finite normal factor of the proper generically finite
map $Z'\to S'$, so $Z'\to V'$ is proper birational.
Let $K_H$ be the function field of the compact Stein base refining
$H$ after both covers, so that the map to this base has connected fibers. It is algebraically closed in
$\CC(Z')$, has transcendence degree $k$, and contains $\CC(W)$:
the map to $W$ is constant on the very general connected leaves.
Thus $K_H/\CC(W)$ is finite. The regular extension
$\CC(S')/\CC(W)$ is linearly disjoint from this finite extension.
The map on a very general leaf has degree one, so
\[
 \CC(Z')=\CC(S')K_H.
\]
Indeed the degree of the left side over the compositum is the generic
degree of the map to the fiber product with that Stein base, and
restricts to the degree on the general connected leaf.

Normalize $W$ in $K_H$, obtaining a finite normal $W'\to W$.
Local integrality gives a lift $V'\to W'$. The induced map
$V'\to Q'\times_{Q_1}W'$ is finite, since both sides are finite over
$S'$, and birational by
the field equality. Its target is integral and normal: it is a
principal $U$-torsor over $B_0\times W'$. Thus it is an isomorphism.
For a proper smooth resolution $\widetilde W\to W'$, the
torsor $Q'\times_{Q_1}\widetilde W$ is a proper smooth resolution
of $V'$. It is special and has logarithmic Kodaira dimension at
least $k$, by comparison with $Z'$. Lemma~\ref{lem:canonical-kernel-logarithmic-comparison}
gives
\[
 \overline\kappa(\widetilde W)\ge k=\dim\widetilde W.
\]
But specialness descends to $\widetilde W$ by that same lemma,
contradicting its logarithmic top-form condition when $k>0$.
This proves $\overline\kappa(Z)=0$.
\end{proof}

\begin{lemma}\label{lem:canonical-kernel-dominance}
Let $Y$ be smooth connected special quasi-projective, let
$Q\to B$ be a principal torus torsor over an abelian variety, and
let $f:Y\to Q$ induce a surjection of fundamental groups.
Then $f$ is dominant and has smooth geometrically connected general
fiber.
\end{lemma}

\begin{proof}
Let $S$ be the reduced closure of the image. A proper smooth
resolution $X_S\to S$ is dominated by a proper smooth modification
of $Y$, obtained from the graph. Lemma~\ref{lem:canonical-kernel-logarithmic-comparison}
implies that $X_S$ is special; by
Lemma~\ref{lem:canonical-kernel-special-proper-image}, its logarithmic
Kodaira dimension is zero.
Lemma~\ref{lem:canonical-kernel-zero-images} therefore describes
$S^\nu$ as a principal $U$-torsor over an abelian variety $C$ finite
etale over an abelian translate $b+B_0\subseteq B$.
Since $Y$ is normal, its morphism to $S$ lifts to $S^\nu$ by
integrality. Hence $\pi_1(S^\nu)\to\pi_1(Q)$ is surjective.
In the torsor sequences, its right map is the injection
$\pi_1(C)\to\pi_1(B)$, and its left map is the cocharacter
inclusion for $U\subseteq T$. Surjectivity first forces $B_0=B$
and $C\to B$ to have degree one, then forces $U=T$.
Thus $\dim S=\dim Q$ and the closed $S$ equals $Q$.

Let $L$ be the relative algebraic closure of $\CC(Q)$ in $\CC(Y)$.
It is finite over $\CC(Q)$. Normalize $Q$
in $L$ to obtain a finite normal $v:V\to Q$. Normality of $Y$ gives
a lift $Y\to V$. Every proper smooth
resolution of $V$ is special by the graph and dominance argument
above. Its logarithmic Kodaira dimension is zero by
Lemma~\ref{lem:canonical-kernel-special-proper-image}, so
Lemma~\ref{lem:canonical-kernel-zero-covers}(i) implies that $v$ is etale.
Since $f_*$ factors through $v_*$ and is surjective, the covering
index is one. Thus $L=\CC(Q)$ and $V=Q$.

The geometric generic fiber is integral, since the function-field
extension is regular in characteristic zero. Compactifying the
map and using its proper Stein factor shows that its
general fibers are nonempty open subsets of integral compact
fibers, hence connected. After shrinking the target to avoid the image of the nonsmooth locus,
generic smoothness gives smooth fibers, which are then geometrically
integral over $\CC$.
\end{proof}

\subsection{Exactness of the fiber sequence}\label{subsec:kernel-exactness}

\begin{lemma}\label{lem:canonical-kernel-ordinary-exactness}
Let $Y$ be smooth connected special quasi-projective, let
$Q\to B$ be a principal torus torsor over an abelian variety, and
let $f:Y\to Q$ be dominant with smooth connected general fiber $F$.
Assume that $f_*:\pi_1(Y)\to\pi_1(Q)$ is surjective. Then
\[
 \im(\pi_1(F)\to\pi_1(Y))=\ker f_*.
\]
This subgroup is finitely generated.
\end{lemma}

\begin{proof}
\smallskip\noindent\emph{The quotient and its meridian orders.}
Put $G_Y=\pi_1(Y)$ and $N=\im(\pi_1(F)\to G_Y)$.
Apply \cite[Theorem~2.2]{achinger2019} to $f$ and choose a smooth
connected dense stratum $Q^\circ$ on which
$f^{-1}(Q^\circ)\to Q^\circ$ is a topological bundle with whole
fiber $F$. The bundle sequence, followed by the dense-open surjection,
shows that $N$ is normal and induces
\[
 \eta:\pi_1(Q^\circ)\twoheadrightarrow R=G_Y/N,
 \qquad \varphi:R\twoheadrightarrow\pi_1(Q).
\]

For smooth $V$ projective birational over $Q$, resolve the main
component of $Y\times_QV$ retaining every divisor over $Y$.
Proper smooth birational modification preserves the fundamental group:
the common open omits only codimension two downstairs and surjects
on the upstairs group. Thus $G_Y$ and $N$ remain unchanged, and
the bundle argument gives compatible epimorphisms
$\eta_V:\pi_1(V^\circ)\twoheadrightarrow R$.
For every interior prime $D\subset V$, let $m_D$ be the
order of its meridian in this fixed group $R$, allowing infinity. Put
$\Delta_V=\sum_{m_D>1}(1-1/m_D)D$. Its support is finite, since orders are one on the bundle locus.
Relations from retained source divisors already hold in $R$.
Common resolutions preserve nonexceptional orders. Recompute
exceptional orders in $R$ and resolve the support of $\Delta_V$;
this gives a nonempty cofinal collection $\mathcal C$ with SNC support.
Sections on upper models push down: interior coefficients agree at
common primes, the deleted boundary has coefficient one, and
normality extends sections across codimension two. Thus the minimum
over $\mathcal C$ is a lower bound for the Iitaka dimension on every
model with these meridian weights, even if exceptional coefficients
decrease.

\smallskip\noindent\emph{A positive minimum and finite auxiliary weights.}
Suppose $\ker\varphi\ne1$. A transverse null-homotopy disk shows
that this kernel is normally generated by omitted-divisor meridians.
Hence some prime $D^\circ\subset Q$ has $m_D>1$.
Choose the compactification of
Lemma~\ref{lem:canonical-kernel-divisor-positivity} for this prime.
For any model with these meridian weights, compactify its projective
morphism to $Q$ over
this compactification, with inverse image of $Q$ exactly the model.
Its logarithmic canonical divisor $B_V=K_{\overline V}+D_V$ is
effective by \eqref{eq:canonical-kernel-valuation}. At the strict
transform of $D^\circ$, $\overline\Delta_V$ has positive coefficient.
At every other prime over its closure the same formula has positive
coefficient, because an interior exceptional center has codimension
at least two and no boundary stratum lies in that closure.
Consequently
\[
 B_V+\overline\Delta_V\ge\epsilon_V g_V^*D
\]
for some rational $\epsilon_V>0$. It follows that
\[
 b=\min_{V\in\mathcal C}\kappa(V,\Delta_V)
 \quad\text{satisfies}\quad 1\le b\le\dim Q.
\]

Choose a minimizing compactified model with effective logarithmic
canonical divisor $B_V$, and put $C_V=\overline\Delta_V$ there.
The inequalities
\[
 \tfrac12(B_V+C_V)\le B_V+\tfrac12C_V\le B_V+C_V
\]
show that the middle divisor and the original adjoint generate the
same full field of section ratios. Inclusions of effective divisors
preserve ratios, and positive rational scaling preserves the full field by
passing to divisible degrees. Resolve this Iitaka map to
% All ratios, rather than only their transcendence degree, are retained.
$h:\overline V_1\to H$, with connected fibers and $\dim H=b$.
Explicitly, a ratio $s/t$ survives a Veronese of index $a$ as
$(s\,t^{a-1})/t^a$; multiplication by the same complementary
effective section preserves it under either divisor inclusion.
Retain all interior exceptional divisors in $V_1$.
If an interior exceptional prime has discrepancy $a_E>-1$ for the
lower pair and meridian coefficient $\beta_E\in[0,1]$, choose its
auxiliary coefficient $t_E$ with
\[
 \max(0,-a_E,\beta_E/2)\le t_E<1.
\]
Keep half the original coefficients on the strict transforms, and
assign coefficient one only to deleted primes. The resulting finite auxiliary divisor
$C_{\rm aux}$ satisfies
\[
 K_{\overline V_1}+D_{V_1}+C_{\rm aux}
   =\mu^*(B_V+\tfrac12C_V)+E_{\rm aux},
 \qquad E_{\rm aux}\ge0\ \text{exceptional}.
\]
Thus its section spaces and full Iitaka field are unchanged. Every
nonzero interior coefficient of $\Delta_{V_1}$ has positive auxiliary
coefficient, including those equal to one.

\smallskip\noindent\emph{The leaves and a common quotient torsor.}
For a very general $w\in H$, put $H_w=h^{-1}(w)\cap V_1$.
Apply \cite[Notation~4.1 and Lemma~4.2]{taji2021} to a Cartier
multiple of $B_V+\tfrac12C_V$ and the corresponding multiple of
$E_{\rm aux}$. Adjunction on the very general compact fiber,
with its full deleted boundary, gives weighted Iitaka dimension
zero on $H_w$, as in
Lemma~\ref{lem:canonical-kernel-special-proper-image}.
Its unweighted logarithmic
dimension is also zero: it is nonnegative by
Lemma~\ref{lem:canonical-kernel-logarithmic-comparison} and at most
the weighted dimension. The map from $H_w$
to $Q$ is projective and birational onto its closed image $S_w$,
because a very general leaf meets the isomorphism locus of
$V_1\to Q$ densely. By
Lemma~\ref{lem:canonical-kernel-zero-images}, $S_w^\nu$ is a smooth
torus torsor and $H_w\to S_w^\nu$ is proper birational.

Every positive component of $C_{\rm aux}|_{H_w}$ is exceptional
for this map. If such a component dominated a prime on $S_w^\nu$,
adapt the compactification of that torsor to the prime using
Lemma~\ref{lem:canonical-kernel-divisor-positivity}. Formula
\eqref{eq:canonical-kernel-valuation} and its positive coefficient
on the nonexceptional component would bound the full weighted
adjoint below by a positive rational multiple of the pullback of that prime.
All other primes above it have positive logarithmic discrepancy
as in the preceding minimum argument. The weighted Iitaka dimension
would then be positive, a contradiction.

Apply \cite[Theorem~2.2]{achinger2019} to the single morphism
$h|_{V_1}$. Over a connected dense parameter stratum, the resulting
topological trivializations transport the whole leaves $H_w$.
Evaluating this transport in $Q$, as in
Lemma~\ref{lem:canonical-kernel-special-proper-image}, shows that their
fundamental-group images in $\pi_1(Q)$ are conjugate.
Their intersections with the central fiber lattice and
their rational projections to $\pi_1(B)$ determine one connected
subtorus $U$ and one abelian subvariety $B_0$. Their dimensions sum
to $\dim Q-b$. The inverse image of this parameter open in $V_1$
defines the incidence variety, with maps to the parameter space and $Q$.
Both are dominant, and each whole fiber image is the closed $S_w$.
Lemma~\ref{lem:canonical-kernel-family-splitting} gives a single
abelian isogeny pullback $Q'\to Q$ and a smooth map
\[
 r:Q'\to Q_1,\qquad
 Q'\to B_0\times Q_1\text{ a principal }U\text{-torsor},
 \qquad\dim Q_1=b.
\]
Very general connected lifted leaves map by proper birational morphisms
onto whole $r$-fibers.

The pullback $Y'=Y\times_QQ'$ is connected by surjectivity of $f_*$,
and is special finite etale over $Y$. Put $G'=\pi_1(Y')$ and
$R'=G'/N\le R$; the fiber image $N$ is unchanged. Interior meridians
map trivially to $\pi_1(Q)$, hence lift with winding one and retain
their orders in $R'$, including infinity.

\smallskip\noindent\emph{Horizontal meridians and the quotient by the leaves.}
Every prime of $Q'$ dominating $Q_1$ has meridian order one. Otherwise
its strict transform on the lifted $V_1$ has positive auxiliary
coefficient and meets a very general leaf in a noncontracted
divisor, contradicting the fact that this support is exceptional.
A vertical prime of $r$ is the irreducible reduced inverse image of
a prime of $Q_1$: the map is smooth and each fiber is a connected
$U$-torsor over $B_0$. Thus the divisor defined by the meridian
orders on $Q'$ is $r^*\Theta$.
For every smooth projective birational $W\to Q_1$, the pullback
$E_W=Q'\times_{Q_1}W\to B_0\times W$ is a principal $U$-torsor,
projective birational over $Q'$.
Every horizontal prime is nonexceptional over $Q'$ and has order
one; every vertical prime is the reduced inverse image of a
prime of $W$. Hence recomputing meridian orders gives the divisor
$p_W^*\Theta_W$. Further such models make its support SNC.

Fill all order-one primes and codimension-two subsets. On the open
$E^\circ=r^{-1}(Q_1\setminus\Supp\Theta)$ this gives an epimorphism
$\eta^\circ:\pi_1(E^\circ)\to R'$.
The map $r$ is a topological bundle: the projective-bundle
compactification over a small base neighborhood has smooth compact
fibers and relative SNC boundary, so local trivializations of the
compact pair restrict to its open complement. Let $L_r$ be the image
in $R'$ of the group of a general whole $r$-fiber $Z_z$. It is normal
by the bundle sequence.
Its composite to $\pi_1(Q')$ is injective. In fact the fiber and
total torsor sequences have injective outer maps
\[
 \pi_1(U)\hookrightarrow\pi_1(T),\qquad
 \pi_1(B_0)\hookrightarrow\pi_1(B_0\times B_1),
\]
so their middle map is injective as well.

The morphism $rf':Y'\to Q_1$ has connected geometric generic
fiber, because both successive function-field extensions are
regular. A smooth connected general fiber $Y'_z$ maps dominantly
to $Z_z$ with connected general fiber $F$. The argument on the open set where the map is a bundle
shows that $\pi_1(Y'_z)\to\pi_1(Z_z)$ is surjective. Let $K'$ be
its actual image in $G'$. It contains $N$, and
\[
 K'/N=L_r,\qquad G'/K'=R'/L_r.
\]
To prove the equality, represent a loop in $Y'_z$ on its dense
$f'$-bundle open. Modulo $N$ its image is the image under
$\eta^\circ$ of its projected loop in $Z_z$, hence belongs to
$L_r$. Conversely the connected-fiber bundle sequence and the
dense-open surjection make $\pi_1(Y'_z)\to\pi_1(Z_z)$ onto,
so every element of $L_r$ occurs. The equality therefore holds for
the whole-fiber image, including relations from retained primes over
base codimension two.
In particular $K'$ is normal. For $\varphi':R'\to\pi_1(Q')$,
injectivity on $L_r$ gives $L_r\cap\ker\varphi'=1$, so quotienting
by $L_r$ preserves every order in this kernel, including infinity.

Each interior meridian on every $E_W$ lies in this kernel, since
its transverse disk maps to $Q'$. The bundle sequence over the
complement of $\Theta_W$ consequently gives compatible epimorphisms
\[
 \beta_W:\pi_1(W\setminus\Supp\Theta_W)\twoheadrightarrow G'/K'
\]
whose meridian orders are exactly those computed in
$R'$. The fiber subgroup $L_r$ is unchanged on every model.

\smallskip\noindent\emph{The lower bound on quotient models.}
For any such $W$, resolve equivariantly the simultaneous graph of all deck
translates of $E_W\to Q'$. The isogeny cover is Galois and its
action is free, so the resulting smooth $J$ descends to a smooth
projective birational $J_0\to Q$ with $J=J_0\times_QQ'$.
Projectivity descends by tensoring translates of a relatively ample
bundle. A cofinal $V_0\in\mathcal C$ dominating $J_0$ therefore
has pullback $V'_0$ dominating $E_W$.

The meridian coefficients on $V'_0$ are the etale pullbacks of those
on $V_0$, because $R'\hookrightarrow R$ preserves the orders.
Divisible logarithmic sections pull back, including at infinity:
for coefficients $\gamma_j\le1$, the Jacobian bound and coefficient
one on the full deleted boundary upstairs give the nonnegative correction
$\sum_j(1-\gamma_j)a_j\ge0$ at every deleted valuation.
Pushing these sections down along $V'_0\to E_W$ then gives
\[
 \kappa(E_W,p_W^*\Theta_W)\ge
 \kappa(V'_0,\Delta_{V'_0})\ge
 \kappa(V_0,\Delta_{V_0})\ge b.
\]
Lemma~\ref{lem:canonical-kernel-logarithmic-comparison} identifies
the left side with $\kappa(W,\Theta_W)$. Since $\dim W=b$, the
adjoint with full deleted boundary on every such $W$ is big.

\smallskip\noindent\emph{Pullback multiplicities and specialness.}
Fix a smooth projective compactification $\overline Y'$ with full
boundary $D_{Y'}$.
We use the following preparation for rational maps from a fixed source.
Let $X_*$ be a fixed smooth projective variety and
$X_*\dashrightarrow Z_*$ a dominant rational map to an integral
projective variety. Put $d=\dim X_*-\dim Z_*$.
Take its integral projective graph and a generic-flatness open of
$Z_*$. Apply Raynaud--Gruson flattening in the form of
\cite[Lemma~76.39.1, Tag~\texttt{087E}]{stacks2026}: a blowup of
the base, with center outside this open, makes the strict main
transform flat and of finite presentation. The strict transform
is the closure of the part over that open and is still proper
birational over $X_*$. It is projective and flat over the modified
base; constancy of the Hilbert polynomial gives fiber dimension
$d$ everywhere. Normalize the base and this main component.
Normality gives the lift to the normalized base, and the
normalization of the main component is finite. Its fibers over
the normalized base therefore have dimension at most $d$.

Every prime $A$ of the original $X_*$ has a divisorial center on
this normalized main component, dominating $A$: a proper
birational morphism to the smooth $X_*$ is an isomorphism at its
generic point. By the dimension bound, the image of this center on the normalized
base has dimension at least $\dim Z_*-1$.
It thus dominates the base or a prime $D$ of the base. In the
second case the restriction of its valuation to $\CC(Z_*)$ is
$e\,\operatorname{ord}_D$ for a positive integer $e$, since the
local ring at $D$ is a DVR. On every higher proper birational
normal base the center of this valuation is the generic point of
the strict transform of $D$. In the first case the restriction
is trivial and the center continues to dominate the base.
Consequently further base resolutions and blowups making any
specified finite divisor support SNC preserve this conclusion
for every prime of the fixed $X_*$.
Resolve a dominating source graph after these base changes.
Every source prime mapped into base codimension at least two
is exceptional over $X_*$.

Apply this construction with $X_*=\overline Y'$ and $Z_*$ a
compactification of $Q_1$, then resolve the recomputed meridian
support and the full base boundary. Over $Q_1$ we obtain a smooth
variety $W$ projective birational over $Q_1$, with SNC weighted
support, and a source model retaining all divisors over $Y'$
and mapping to $W$.

No retained source prime maps to compactified infinity: the
map to $Q_1$ is regular, and properness of $W\to Q_1$ extends its
generic lift at that prime's discrete valuation ring into $W$.
For a retained source prime over a prime $D\subset W$, write
$e$ for its pullback multiplicity. Its small transverse loop is
null-homotopic in the whole source. Compatibility with $\beta_W$
gives $\beta_W(\mu_D)^e=1$. Thus a finite meridian order $m_D$ divides
$e$, and an infinite-order prime has no retained source prime
dominating it.

Let $\mathcal L$ be the saturation of the generic differential
line of this $b$-dimensional base map in
$\Omega^b_{\overline Y'}(\log D_{Y'})$.
Pull back a divisible section of
$K_{\overline W}+D_W+\overline\Theta_W$ to a rational section of
a power of $\mathcal L$. At a retained prime over a meridian of
finite order $m_D$, the differential vanishes to order at least
\[
 e-1\ge e(1-1/m_D),
\]
because $m_D\mid e$. A component of coefficient one has no
retained prime above it. At a deleted source prime,
$d(x^e u)=x^e(eu\,dx/x+du)$ gives logarithmic vanishing at least
$e$, enough for every coefficient at most one. A prime dominating
the base gives no pole. These checks cover all primes on the
fixed $\overline Y'$, since the remaining contracted primes on
the prepared source are exceptional over it. Since $\mathcal L$ is a
line bundle by
Lemma~\ref{lem:canonical-kernel-logarithmic-comparison}, normality
extends these sections to global sections of $\mathcal L^{\otimes m}$.
All degrees use the same generic differential map. Its powers
give multiplicative injections
\[
 H^0\bigl(\overline W,
   m(K_{\overline W}+D_W+\overline\Theta_W)\bigr)
 \hookrightarrow H^0(\overline Y',\mathcal L^{\otimes m})
\]
for every sufficiently divisible $m$. The differential factor
cancels in equal-degree ratios, so these maps identify each base
ratio with its usual pullback in $\CC(Y')$.
Dominance preserves their independent ratios, giving
\[
 \kappa(\overline Y',\mathcal L)\ge b>0.
\]
This contradicts specialness of the finite etale cover $Y'$, so
$\ker\varphi=1$ and $N=\ker f_*$.
The group $\pi_1(F)$ is finitely generated by compatible finite
triangulation, as in the proof of
Theorem~\ref{thm:canonical-albanese}(i); hence so is its image $N$.
\end{proof}

\begin{proof}[Proof of part (iii) of Theorem~\ref{thm:canonical-albanese}]
The source $Y$ is special by
\cite[Proposition~10.11]{campana2009}. Part (ii) identifies
the map $f_*$ with the quotient $G\twoheadrightarrow E$.
By Lemma~\ref{lem:canonical-kernel-dominance}, $f$ is dominant with
smooth geometrically connected general fiber.
Lemma~\ref{lem:canonical-kernel-ordinary-exactness} then identifies
its fiber image with $\ker f_*=P$ and proves finite generation.

The construction is compatible with every connected finite etale target cover
$j:Q'\to Q$. If $E'=j_*\pi_1(Q')$, the surjectivity of $f_*$ implies that
$Y'=Y\times_QQ'$ is connected, with
$G'=\pi_1(Y')=f_*^{-1}(E')$. The lift $f':Y'\to Q'$ has
the same whole general fiber $S$, because the cover restricts
trivially over that fiber. Therefore, inside the original $G$,
\[
 P\subseteq G',\qquad
 1\longrightarrow P\longrightarrow G'
   \xrightarrow{\,f'_*\,}E'\longrightarrow1,\qquad
 \im(\pi_1(S)\to G')=P.
\]
The cover $Y'$ is special, and $Q'$ is a torus torsor over an
abelian variety by Proposition~\ref{prop:finite-cover-torus-torsor}.
\end{proof}

\subsection{The FAb kernel and finite covers}\label{subsec:fab-kernel-covers}

\begin{lemma}\label{lem:canonical-kernel-solvable-image}
Let $X$ be smooth connected special complex quasi-projective,
and let $Y,G,K,P,E$ be canonical kernel data as in
Definition~\ref{def:canonical-kernel}. If a homomorphism $\psi:G\to R$
has virtually solvable image, then $\psi(P)$ is finite.
\end{lemma}

\begin{proof}
Replace $R$ by $\psi(G)$. Choose a solvable
finite-index subgroup $R_s\le R$. Its inverse image under $\psi$ is
the group of a special finite etale cover of $Y$.
By \cite[Corollary~7.2]{rogov2025}, $R_s$ is virtually nilpotent. Then
\cite[Proposition~2.7(1)]{pan2025} gives a torsion-free nilpotent
subgroup of finite index, and \cite[Theorem~2.10]{dekimpe2017}
embeds it in an integral unitriangular group. Apply
\cite[Theorem~B]{cao2026} to this linear representation of the
corresponding special cover. Taking the normal core in $R$ gives
a normal finite-index torsion-free class-two subgroup $N$.
Put $G_N=\psi^{-1}(N)$. Maximality of $b_1(Y)$ and injectivity of
restriction give
\[
 \operatorname{Hom}(G,\QQ)\xrightarrow{\sim}
 \operatorname{Hom}(G_N,\QQ).
\]
Every rational character of $N$, pulled back to $G_N$, therefore
extends to $G$. It follows both that conjugation by $R$ is trivial on
$N_{\rm ab}\otimes\QQ$ and that $\psi(K\cap G_N)$ has zero image there.
For the second assertion we used that rational characters of $G$ vanish
on $K$, by the full integral period identification in the proof of
Theorem~\ref{thm:canonical-albanese}(i); see also
\cite[Definition~2.31 and Lemma~2.34]{elduque2024}.

The class-two calculation in the proof of
Theorem~\ref{thm:canonical-albanese}(i) now applies to $N$: conjugation
by $R$ fixes $[N,N]$ pointwise, and it fixes any element of $N$ having
a positive power in $[N,N]$. Thus
$M_0=\psi(K\cap G_N)$ is torsion-free abelian and central in $R$.
It has finite index in the finitely generated group $M=\psi(K)$.
Since $M_0\le Z(M)$ and $[M:M_0]<\infty$, Schur's theorem
\cite[Theorem~1.1]{dixon2015schur} implies that $[M,M]$ is finite.
Here $[M,M]=\psi([K,K])\le\psi(P)$.
Finally $\psi(P)/[M,M]$ is an image of the finite group
$P/[K,K]=\Tor K_{\rm ab}$. Both this quotient and $[M,M]$ are finite,
so $\psi(P)$ is finite.
\end{proof}

\begin{proposition}\label{prop:canonical-kernel-covers}
For canonical kernel data, $P$ is FAb. For every finite-index subgroup
$L\le P$, there are a characteristic finite-index subgroup $P_0\le P$
contained in $L$ and a normal finite-index subgroup $G_0\le G$ such
that $G_0\cap P=P_0$. Thus $G_0/P_0\to E$ is injective with finite-index
image.
\end{proposition}
\begin{proof}
\smallskip\noindent\emph{A prescribed normal finite quotient of $P$.}
Let $P_0\le P$ have finite index and be normal in $G$. Put
$R_0=G/P_0$ and $D_0=P/P_0$. The centralizer $C_0=C_{R_0}(D_0)$ has
finite index because $D_0$ is finite. Its intersection with $D_0$ is
central, and its quotient by that intersection embeds into $E$.
Consequently
\[
 \gamma_3(C_0)\subseteq C_0\cap D_0\subseteq Z(C_0),
 \qquad\gamma_4(C_0)=1.
\]
The inverse image of $C_0$ in $G/[P_0,P_0]$ is an extension of this
nilpotent group by the abelian group $P_0/[P_0,P_0]$. Thus
$G/[P_0,P_0]$ is virtually solvable. Applying
Lemma~\ref{lem:canonical-kernel-solvable-image} gives
\[
 |P/[P_0,P_0]|<\infty.
\]
The triangulation of the compactification used in
Theorem~\ref{thm:canonical-albanese}(i) shows that $G$ is finitely
generated. Its quotient $R_0$ is finitely generated, and so is
the finite-index subgroup $C_0$. By
\cite[Proposition~2.7(1)]{pan2025}, it has a torsion-free
subgroup of finite index. Its normal core $N_0$ in $R_0$ is still
torsion-free and has finite index, so $N_0\cap D_0=1$.
The inverse image $G_0$ of $N_0$ in $G$ is normal of finite index and
satisfies $G_0\cap P=P_0$. Hence $G_0/P_0\to E$ is injective with
finite-index image.

\smallskip\noindent\emph{All finite-index subgroups of $P$.}
By Theorem~\ref{thm:canonical-albanese}(iii), $P$ is finitely generated.
Given $L\le P$ of index $a$, let $P_0$ be the intersection of the
finitely many subgroups of $P$ of index at most $a$. It is
characteristic of finite index in $P$, lies in $L$, and is normal in $G$.
The preceding argument constructs the required $G_0$ and proves that
$P/[P_0,P_0]$ is finite.
Since $[P_0,P_0]\subseteq[L,L]$, the abelianization $L/[L,L]$ is finite.
\end{proof}

\begin{proposition}\label{prop:finite-kernel-criterion}
For canonical kernel data, if $P$ is finite, then $G$ has a
finite-index subgroup of nilpotency class at most two.
\end{proposition}
\begin{proof}
Apply Proposition~\ref{prop:canonical-kernel-covers} to
$\{1\}\le P$. It gives a finite-index subgroup $G_0\le G$ with
$G_0\cap P=1$. Thus $G_0$ embeds in $E=G/P$, which has class at
most two by Theorem~\ref{thm:canonical-albanese}(i).
\end{proof}

\section{Meridian groups}\label{sec:10-ordinary-endpoints}

We compare the meridian groups defined in \eqref{eq:meridian-group}
under birational maps and study their boundary collar images.
The birational comparisons are used in
Section~\ref{sec:13-ordinary-minimum}; the results on collar images
enter the proof of virtual solvability in
Section~\ref{sec:12-bounded-special}.

\subsection{Birational comparison}\label{subsec:meridian-birational-comparison}

\begin{proposition}\label{prop:ordinary-meridian-comparison}
Let $d:T\to T'$ be a divisorial contraction between normal projective
varieties in a negative minimal model program for a standard pair
$(T,C)$, and put $C'=d_*C$. Comparison on the common open in the
regular loci, with the boundaries removed, gives an epimorphism
\[
 \Orbgrp{T'}{C'}\twoheadrightarrow\Orbgrp TC.
\]
Its kernel is the normal closure of $\lambda_E^{m_E}$ for the
contracted primes $E$ of finite original weight $m_E$, including
weight one for an unmarked prime. Here $\lambda_E$ is the transported
generic normal meridian. An infinite-weight prime contributes no
relation. A flip with the unchanged surviving weights induces an
isomorphism of these groups.

If $S$ is smooth projective with SNC standard boundary $\Delta$ and
$r:T\to S$ is a smooth projective birational morphism, let $\Delta_T$
retain every strict-transform weight and assign weight infinity to
every exceptional prime. Then the common-open map induces
\[
 \Orbgrp T{\Delta_T}\simeq\Orbgrp S\Delta.
\]
\end{proposition}

\begin{proof}
Removing complex codimension at least two from a smooth variety does
not change $\pi_1$: loops and homotopy disks avoid a stratification
by transversality. Restoring a divisor instead kills exactly the
normal closure of its generic meridians, as a transverse
null-homotopy disk, punctured at its intersections, expresses its
boundary as a product of their conjugates.

Remove the contracted primes and all boundary primes to obtain a
common regular open $O$. All remaining omissions from the relevant
regular loci have codimension at least two. Let $\mathcal H$ be
$\pi_1(O)$ modulo the common finite boundary powers. Since no prime
is restored on $T'$, its meridian group is $\mathcal H$; on $T$,
an unmarked contracted prime kills its first meridian, a finite
weight $m_E$ kills its $m_E$-th power, and an infinite weight adds
no relation. Hence
\[
 \Orbgrp TC=\mathcal H/
 \langle\!\langle\lambda_E^{m_E}:m_E<\infty\rangle\!\rangle.
\]
For a flip the common presentations have identical prime relations,
giving an isomorphism.
For the smooth birational map $r$, its exceptional image downstairs
has codimension at least two and every exceptional prime upstairs
has infinite weight. Again the presentations are identical.
\end{proof}

\subsection{Boundary collars and adjunction}\label{subsec:boundary-collars-adjunction}

\begin{definition}\label{def:whole-boundary-collar}
Let $V$ be normal projective with a standard boundary $B$, and let
$S$ be a component of $\lfloor B\rfloor$. Put
\[
 R_S=(S_{\rm reg}\cap V_{\rm reg})\setminus
       \bigcup_{B_i\ne S}B_i.
\]
A punctured normal tube over the whole component is a punctured
complex normal tube $P_S$ over all of this smooth embedded hypersurface.
Its positive radius may decrease towards the omitted locus. Its
collar group is
$J_S=\im(\pi_1(P_S)\to\Orbgrp VB)$ after a chosen base-point path.
A positive circle fiber defines $h_S\in J_S$.
\end{definition}

The stratum $R_S$ is a closed smooth hypersurface in the regular
open obtained by removing the other boundary components and the
singular locus of $S$. Tubular neighborhoods with decreasing positive
radii exist there. Two choices admit a common smaller tube and give
the same image by radial deformation and tubular isotopy; changing
the base-point path conjugates that image. The complex normal line
orients its circle globally.

\begin{theorem}\label{thm:reduced-collar-detection}
Let $(V,A)$ be a connected normal projective lc pair of dimension at
least one, with $A\ne0$ reduced and $K_V+A\sim_{\QQ}0$.
If $G=\pi_1(V_{\rm reg}\setminus\Supp A)$ is virtually solvable,
then the collar image $J_S$ of every component has finite index in $G$.
\end{theorem}

\begin{proof}
Put $U=V_{\rm reg}\setminus\Supp A$ and $n=\dim V$. If $n=1$,
the degree equation gives $(V,A)=(\mathbb P^1,\{0,\infty\})$;
every boundary tube, also after a finite cover, generates
$\pi_1(\CC^*)$. Assume $n\ge2$.

\emph{Detection by one tube.}
Take a projective log resolution $f:X\to V$, isomorphic over $U$,
with full SNC deleted boundary $B$. Then
\[
 K_X+B\sim_{\QQ}E
   =\sum_{F\text{ exceptional}}(a(F,V,A)+1)F\ge0.
\]
For the strict transform $C$ of a boundary component and
$L=f^*H$ with $H$ ample on $V$, one has
$E\cdot L^{n-1}=0$ and $C\cdot L^{n-1}>0$.
A nonzero logarithmic one-form vanishing along $C$ has saturated
image $\mathcal O_X(Z)$, $Z\ge C$, with torsion-free quotient
$\mathcal Q$ of rank $n-1$. Hence
\[
 c_1(\mathcal Q)\cdot L^{n-1}
 =(K_X+B-Z)\cdot L^{n-1}<0.
\]
The inequality persists for a small ample perturbation of $L$.
Since the pair is SNC and $K_X+B\sim_{\QQ}E\ge0$ is pseudoeffective,
Campana--P\u{a}un's generic semipositivity theorem
\cite[Theorem~2.1]{campana2015} gives the opposite inequality for
every torsion-free quotient of $\Omega_X^1(\log B)$.
Thus no such form exists.

For $\alpha\in H^1(U,\RR)$, let $v_j$ be its periods on positive
meridians of the components $B_j$. The residue exact sequence gives
$\sum_jv_j[B_j]=0$ in $H^2(X,\RR)$.
By the residue criterion \cite[Theorem~1.5 and Corollary~4.1]{fang2018},
there is a closed logarithmic form with residues $v_j/(2\pi i)$;
denote it by $\omega_0$.
The class $\alpha-[\operatorname{Re}\omega_0]$ comes from
$H^1(X,\RR)$, hence is the real part of a holomorphic form $\theta$.
Thus $\alpha=[\operatorname{Re}\omega]$ for
$\omega=\omega_0+\theta$ with the same residues.

Suppose $\alpha$ vanishes on the punctured tube over
$C^o=C\setminus\bigcup_{B_j\ne C}B_j$.
Vanishing on the normal circle forces the residue at $C$ to be zero;
small loops near $C\cap B_j$ give zero residues along every component
meeting $C$. Thus $\omega$ is holomorphic near $C$. The fundamental
group of the tube surjects onto $\pi_1(C^o)$; shrinking lifted loops
in the normal direction shows that $\omega|_C$ has zero real periods.
Since $\pi_1(C^o)\twoheadrightarrow\pi_1(C)$,
compact Hodge theory gives $\omega|_C=0$. Together with its zero
residue, this means that all logarithmic coefficients vanish along
$C$, so the preceding argument gives $\omega=0$ and $\alpha=0$.

Restricting a tube to a dense algebraic open of its smooth base
preserves its group image. Indeed, the base inclusion is surjective on
$\pi_1$, so the projected class of any loop in the whole tube lifts
to a loop in the restricted circle bundle. By the bundle homotopy
sequence, the two classes differ by an element in the image of the
normal circle, which also lies in the restricted bundle. Such a
restriction of the resolution tube maps into the prescribed $P_S$
after decreasing radii. Consequently
$H^1(U,\RR)\to H^1(P_S,\RR)$ is injective.

\emph{All covers and finite index.}
For any connected finite etale $U'\to U$, normalize $V$ in
$\CC(U')$ to obtain $g:V'\to V$, with
$A'=(g^{-1}A)_{\rm red}$. Hurwitz and absence of ramification off
$A$ in codimension one give
\[
 K_{V'}+A'=g^*(K_V+A),\qquad
 a(F',V',A')+1=e(F'/F)(a(F,V,A)+1)\ge0.
\]
Resolve also $V'\setminus U'$; all exceptional primes remain deleted.
The difference between $U'$ and $V'_{\rm reg}\setminus A'$
has codimension at least two, since it lies over the singular locus
of $V$ under a finite morphism. Resolving that deleted set preserves
the same effective exceptional representative used above.
The same form argument applies to each strict boundary component.
Locally a cover over $R_S$ is $(z,w)\mapsto(z^e,w)$ after etale
base change. Each connected lifted tube therefore contains, over a
dense open, the tube of one prime over $S$. For every such lift $P'$,
\[
 H_1(P',\QQ)\twoheadrightarrow H_1(U',\QQ).
\]

The group $G$ is finitely generated by compatible triangulation
\cite[Theorems~3.4.4 and~9.2.1]{bochnak1998}.
Since $G$ is virtually solvable, it has a finite-index nilpotent
subgroup $G_1$ by \cite[Corollary~6.6]{rogov2025}, realized by a
finite etale cover
\cite[Expos\'e~XII, Theorem~5.1]{grothendieck2004}.
A lifted collar has image $M=J_S\cap G_1$ spanning
$H_1(G_1,\QQ)$, so $M$ has finite-index image in the integral
abelianization. Choose generators $s_j$ of $G_1$ and $m>0$ with
$t_j\in M$ congruent to $s_j^m$ modulo $[G_1,G_1]$.
In $\gamma_a(G_1)/\gamma_{a+1}(G_1)$, the corresponding length-$a$
commutators in the $t_j$ are the $m^a$-th powers of the generating
commutators in the $s_j$. Thus $M\cap\gamma_a(G_1)$ has finite-index
image in every successive quotient. Descending the finite lower
central series proves $[G_1:M]<\infty$, hence $[G:J_S]<\infty$.
\end{proof}

\begin{theorem}\label{thm:positive-standard-collar}
Let $(V,B)$ be a connected projective $\QQ$-factorial dlt standard
pair of dimension $n\ge2$. Let $S$ be a component of $\lfloor B\rfloor$
and let $A$ be an ample Cartier divisor with $S^2A^{n-2}>0$.
Let $\Theta=\operatorname{Diff}_S(B-S)$ be the adjunction different,
so that $K_S+\Theta\sim_{\QQ}(K_V+B)|_S$, and put $G=\Orbgrp VB$.
Then $S$ is normal, $(S,\Theta)$ is dlt with standard coefficients,
and its collar over the whole component satisfies
\[
 [G:J_S]<\infty,\qquad h_S\in Z(J_S),\qquad
 \Orbgrp S\Theta\twoheadrightarrow J_S/\langle h_S\rangle.
\]
Both groups are finitely generated. In particular, virtual
solvability of $\Orbgrp S\Theta$ implies virtual solvability of $G$.
If $\rho(V)=1$, the positivity assumption holds for every floor
component and every ample Cartier $A$.
\end{theorem}

\begin{proof}
\emph{The circle quotient.}
By standard dlt adjunction, $S$ is normal and the adjunction pair is dlt
with standard coefficients
\cite[Propositions~6.12 and~6.18, Example~6.8]{graf2013}.
Outside ambient codimension three, transverse charts along $S$ are
$(\mathbb D_x\times\mathbb D_y)/\mu_r$ times a smooth factor, with
$\zeta(x,y)=(\zeta^b x,\zeta y)$, $S=(x=0)/\mu_r$ and
$\gcd(b,r)=1$ if $r>1$. The other branch $y=0$ has weight $m$,
including $m=1$ when absent; the different has weight $rm$.
The normal circle is central in the tube sequence
$\ZZ\to\pi_1(P_S)\to\pi_1(R_S)\to1$.
In the chart, let $h,\beta$ be full turns in $x,y$, and let $t$
have turns $(b/r,1/r)$. Then $t^r=h^b\beta$ and $t$ projects to
the positive meridian of $v=y^r$ on $S$. For finite $m$,
$\beta^m=1$ in $G$ (for $m=1$ it bounds a disk), so
$t^{rm}=h^{bm}$. Modulo $h$, this is precisely the different
relation. An infinite weight imposes no relation.
The remaining discrepancy between $R_S$ and
$S_{\rm reg}\setminus\Supp\Theta$ has codimension at least two.
Thus $h_S\in Z(J_S)$ and
$\Orbgrp S\Theta\twoheadrightarrow J_S/\langle h_S\rangle$.

\emph{Positivity and the whole collar.}
Replace $A$ by a high multiple and take a general intersection
surface $W$ of $n-2$ members, avoiding the preceding exceptional
set and transverse to all retained strata; set $W=V$ when $n=2$.
By Bertini \cite[Lemmas~2.23--2.25]{greb2011}, $W$ is normal and
$C=S|_W$ is smooth and irreducible, with
$C^2=S^2A^{n-2}>0$. The transverse quotient charts are unchanged.
For $U_W=W\cap(V_{\rm reg}\setminus\Supp B)$,
\cite[Theorem~1.1.3(ii)]{hamm1985} gives
\[
 \pi_1(U_W)\simeq\pi_1(V_{\rm reg}\setminus\Supp B).
\]
The singular locus stays deleted at each cut, and the variety being
cut has dimension at least three.

Resolve $W$ by $f:Y\to W$, isomorphically over $W_{\rm reg}$.
Writing $E_j$ for the exceptional curves over $C$, one has
$F=f^*C=\widetilde C+\sum a_jE_j$ with $a_j>0$.
Choose a general ample divisor $H_Y$ meeting all $E_j$. Its
pushforward is rational Cartier near these quotient points;
the exceptional part $E=\sum e_jE_j$ of its pullback has $e_j>0$
and $E\cdot E_j=-H_Y\cdot E_j<0$.
For small rational $\varepsilon>0$, $F-\varepsilon E$ is effective
with the same support and positive intersection with every component:
\[
 (F-\varepsilon E)\cdot E_j>0,
 \qquad (F-\varepsilon E)\cdot\widetilde C
       =C^2-\varepsilon E\cdot\widetilde C>0.
\]
Use $F$ if no exceptional curves occur. Clearing denominators gives
a connected effective Cartier thickening $L_Y$ of the whole
$f^{-1}C$, with $\mathcal O_Y(L_Y)|_{L_Y}$ ample by the componentwise
degree criterion for proper curves. Set
$Z_Y=Y\setminus f^{-1}(U_W)$, retaining every exceptional curve.
Nori's WLT(A) \cite[pp.~305--306]{nori1983} implies that
$\pi_1(N_C\cap U_W)$ has finite-index image in $\pi_1(U_W)$ for
every sufficiently small connected neighborhood $N_C$ of $C$:
apply it upstairs inside $f^{-1}(N_C)$ to $L_Y$ and delete $Z_Y$.
Here the locally invertible map required by WLT(A) is the open
inclusion $f^{-1}(N_C)\hookrightarrow Y$. The support of $L_Y$
is connected because $C$ and the fibers of $f$ are connected, and
the ampleness criterion applies to its possibly nonreduced Cartier
thickening.

To capture this image by the prescribed collar, form $N_C$ from a
normal disk tube $N_0$ away from the finitely many quotient and
crossing points, together with small quotient bidisks there.
Each punctured patch is
$(\mathbb D_x^*\times\mathbb D_y)/\mu_r$ or
$(\mathbb D_x^*\times\mathbb D_y^*)/\mu_r$.
Its overlap with $N_0$ has the homotopy type of the quotient two-torus.
The inclusion into the patch is surjective on fundamental groups,
by filling $y=0$ in the first case and by radial deformation in the
second. Shrinking in the nonzero $x$-direction
handles varying radii. Van Kampen gives
$\pi_1(N_0\cap U_W)\twoheadrightarrow\pi_1(N_C\cap U_W)$.
Transversality identifies this tube's normal line with the
restriction of that of $S$; tubular isotopy puts its group image
inside $J_S$. The hyperplane isomorphism and Nori's theorem now give
$[G:J_S]<\infty$.

Both groups are finitely generated by finite
triangulation \cite[Theorems~3.4.4 and~9.2.1]{bochnak1998}.
A solvable finite-index subgroup of $J_S/\langle h_S\rangle$
has solvable inverse image, since the kernel is central cyclic.
If $\rho(V)=1$, then
$[S]=c[A]$ with $c>0$ by $S\cdot A^{n-1}>0$; hence
$S^2A^{n-2}=c^2A^n>0$.
\end{proof}

\begin{lemma}\label{lem:ordinary-endpoints-principal-adjunction}
Let $V$ be a connected normal projective klt variety of dimension
$d\ge2$, and let $D=\sum_iD_i$ be a nonzero reduced divisor.
Assume that $(V,D)$ is log canonical and $K_V+D$ is principal
Cartier. For each component let $W_i\to D_i$ be its normalization,
let $A_i$ be its normalized different, and let $P_i$ be its
tube over the whole stratum $R_i$ of
Definition~\ref{def:whole-boundary-collar}. Then $A_i$ is effective
and reduced, $(W_i,A_i)$ is log canonical, and $K_{W_i}+A_i$ is
principal Cartier. The stratum $R_i$ identifies with an open of
$(W_i)_{\rm reg}\setminus\Supp A_i$ whose complement has complex
codimension at least two. There is an exact sequence
\[
 \ZZ\longrightarrow\pi_1(P_i)\longrightarrow
 \pi_1((W_i)_{\rm reg}\setminus\Supp A_i)\longrightarrow1
\]
in which the image of $\ZZ$ is central and generated by a positive
normal circle.
\end{lemma}

\begin{proof}
\emph{The normalized residue.}
Choose a rational $d$-form $\omega$ with
$\operatorname{div}_V(\omega)=-D$. Its residue along the generic
smooth point of $D_i$ is a rational $(d-1)$-form $\eta_i$ on
$W_i$. The normalized adjunction pair is log canonical by Kawakita's
inversion of adjunction, in the formulation of
\cite[Theorem~1.1]{fujino2023b}. For a sufficiently divisible tensor
power, divisorial adjunction agrees generically with
$\omega^{\otimes m}\mapsto\eta_i^{\otimes m}$, hence agrees
as a rational map on $W_i$. Since $\omega$ trivializes the full
log-canonical sheaf, this gives
$A_i=-\operatorname{div}_{W_i}(\eta_i)$. In particular $A_i$
is integral and $K_{W_i}+A_i$ is principal Cartier.

\emph{The codimension-two boundary.}
Outside ambient codimension three,
\cite[Proposition~9.3 and Section~9.C.4]{greb2011} provides a transverse
smooth surface cover with finite group having no quasi-reflections.
The pulled-back boundary is reduced and lc: the finite Hurwitz formula
multiplies
log discrepancies by the ramification index. Take a general
ample complete-intersection surface transverse to these strata.
It is lc by \cite[Lemmas~2.23--2.25]{greb2011}; on
the smooth cover its equations solve for the invariant coordinates,
giving an lc plane pair.
A reduced lc plane curve is smooth or nodal. Its first blowup
has discrepancy $1-\nu$ at a point of multiplicity $\nu$, so $\nu\le2$.
If a double tangent occurs, its transform is either still singular
on the coefficient-one exceptional curve or tangent to it; the
next one or two blowups produce total boundary multiplicity at
least three. This contradicts lc. The implicit function theorem
and the holomorphic Morse lemma with parameters then straighten
the smooth or nodal branches along the stratum.

With one branch $x=0$, equivariant coordinates give a diagonal
action. The pulled-back log form is a unit times
$dx/x\wedge dy\wedge dw_1\wedge\cdots\wedge dw_{d-2}$.
Invariance at the stratum forces the $y$-character to be trivial.
Since quasi-reflections are excluded, the group is trivial.
With two branches $xy=0$, the corresponding form is a unit times
$dx/x\wedge dy/y\wedge dw_1\wedge\cdots\wedge dw_{d-2}$.
An element exchanging the branches changes its value by a minus
sign at the stratum, so no such element occurs in the local
isotropy group. The group preserves both branches, is diagonal,
and has two faithful characters. It is therefore cyclic, acting
as $(x,y)\mapsto(\zeta^a x,\zeta y)$ with $\gcd(a,r)=1$.
Along a connected crossing stratum, the branch normalization is
finite etale of degree two and may be disconnected or connected.
Thus global branch monodromy can occur even though the local
isotropy group preserves each branch.

On the normalized branch the parameter is $t=y^r$; the residue
is a nonzero unit times
$r^{-1}dt/t\wedge dw_1\wedge\cdots\wedge dw_{d-2}$.
It has order $-1$ there and order zero along a smooth isolated
branch. These divisors include self-conductor primes, intersections
with other $D_j$, and generic singular strata of $V$ contained
in $D$. Hence $A_i$ is effective and reduced, and these are
exactly the divisorial parts omitted from $R_i$ in its normalization.
The remaining omitted sets have codimension at least two on $W_i$,
since normalization is finite.

Transversality identifies the two base fundamental groups. The
homotopy sequence of the complex normal line bundle gives the
displayed exact sequence, and complex transition functions preserve
the positive circle, so its image is central.
\end{proof}

\subsection{Groups with trivial adjoint}\label{subsec:trivial-adjoint-groups}

\begin{theorem}\label{thm:principal-fano-collars}
Let $V$ be a connected projective $\QQ$-factorial klt variety of
dimension $d\ge2$, with $\rho(V)=1$ and $-K_V$ ample. Let
$D=\sum_iD_i$ be nonzero and reduced, with $(V,D)$ log canonical
and $K_V+D$ principal Cartier. For each $i$, let $W_i\to D_i$ be
the normalization and let $A_i$ be its adjunction
different. Put $G=\pi_1(V_{\rm reg}\setminus\Supp D)$ and let
$J_i$ be the collar image of the whole component $D_i$.
Assume, for each $i$ with $A_i\ne0$, that
$\pi_1((W_i)_{\rm reg}\setminus\Supp A_i)$ is virtually solvable.
Then $G$ is virtually abelian.
\end{theorem}

\begin{proof}
By Lemma~\ref{lem:ordinary-endpoints-principal-adjunction}, the
adjunction pairs are reduced and lc with principal adjoints. The fundamental
group of each tube surjects onto $\Gamma_i=\pi_1((W_i)_{\rm reg}\setminus\Supp A_i)$ with
central cyclic kernel generated by its normal meridian $\mu_i$.

\emph{Incidences of whole tubes.}
Remove an algebraic $Z\subset V$ of codimension at least three
containing the exceptions to the lemma's charts and the nongeneric
loci of all crossing strata. The remaining crossing pieces $B_j^o$
are smooth connected with transverse type
$\CC^2/\mu_r(a,1)$ and boundary $xy=0$.
Their branch normalization is either two degree-one covers or one
connected degree-two etale cover; the latter allows branch monodromy.
Shrink the prescribed tubes over $R_i\setminus Z$ to $T_i$, without
changing their images $J_i$, and put $E_i=\pi_1(T_i)$.
Deleted crossing neighborhoods $K_j$ and these tubes cover
$H=O\setminus D$ for a neighborhood $O$ of $D\setminus Z$.
The $T_i$ and the $K_j$ are separately disjoint; overlaps are the
sectors where one nonzero normal coordinate is smaller than the
other. On a locally finite chart cover, choose smaller radii so that
the neighborhoods are compatible on overlaps and under branch permutations.

An overlap $C_e=T_i\cap K_j$ corresponds to a connected component
$Q_e^o$ of the branch cover. The transverse winding lattice is
$\Lambda_{r,a}=\ZZ^2+\ZZ(a/r,1/r)$.
The two sectors have respective bases
$(1,0),(a/r,1/r)$ and $(0,1),(1/r,b/r)$, with $ab\equiv1\pmod r$;
for $r=1$ use the coordinate basis. Both determinants have absolute
value $1/r$, so either sector supplies the whole transverse group.
Let $M_e\le\Gamma_i$ be the image of its projected punctured tube
along the corresponding component of $A_i$. Writing
$p_i:E_i\to\Gamma_i$, the incidence images are
\begin{align}
 \im(\pi_1(C_e)\to E_i)&=p_i^{-1}(M_e),
 \label{eq:ordinary-endpoints-tube-incidence}\\
 \im(\pi_1(C_e)\to\pi_1(K_j))
 &=\bigl(\pi_1(K_j)\to\pi_1(B_j^o)\bigr)^{-1}
       \im\bigl(\pi_1(Q_e^o)\to\pi_1(B_j^o)\bigr).
 \label{eq:ordinary-endpoints-crossing-incidence}
\end{align}
To see these equalities, lift the indicated base loops and correct
by the normal circle in the first case, or by the full transverse
fiber image in the second.
An incidence implies $A_i\ne0$. By its assumed virtual solvability
and Theorem~\ref{thm:reduced-collar-detection}, $M_e$ has finite
index in $\Gamma_i$; dense-open restrictions do not change tube
images. The other incidence image has index one or two, the degree
of the branch cover. Thus every edge image has finite index on both
sides, including self-crossings.

\emph{Generation and commensuration.}
Since $\rho(V)=1$, every $D_i$ is ample rational Cartier; distinct
components meet in codimension two, so the incidence graph is
connected. Take a general intersection surface $V_2$ avoiding $Z$
and transverse to the retained strata, with $V_2=V$ for $d=2$.
By Bertini and the product charts, the surface is normal and klt, with
$(V_2)_{\rm reg}=V_2\cap V_{\rm reg}$
\cite[Lemmas~2.23--2.25]{greb2011}. For $C=D|_{V_2}$,
\cite[Theorem~1.1.3(ii)]{hamm1985} gives
$\pi_1((V_2)_{\rm reg}\setminus C)=G$.
Resolve $f:X_2\to V_2$, and delete
$Z_2=f^{-1}(C\cup\operatorname{Sing}V_2)$ in full.
For $lC$ Cartier, the full pullback $F_2=f^*(lC)$ has a multiple
whose linear system has two-dimensional image. Nori's neighborhood
surjection \cite[Corollary~2.3]{nori1983}, applied inside
$f^{-1}(O\cap V_2)$, therefore factors through $H$ and proves
$\pi_1(H)\twoheadrightarrow G$.

By van Kampen for the incidence graph, $\pi_1(H)$ is generated by
vertex group images and stable letters outside a maximal tree. Since the
edge images have finite index on both sides, the vertex images are
commensurable along that tree.
Each vertex group normalizes itself and commensurates this common
class; each stable letter conjugates finite-index subgroups on its
two ends. Thus the tube images in $G$ are pairwise commensurable
and commensurated. This includes both overlaps at a split
self-crossing and the index-two overlap of a connected branch cover.

\emph{Finite roots.}
For $q\ge2$, the pair $(V,(1-1/q)D)$ is klt, since its log
discrepancy is the corresponding positive combination of those
of $(V,0)$ and $(V,D)$. Its negative adjoint is $D/q$, hence ample.
By \cite[Theorem~2]{braun2020},
\[
 Q_q=G/\langle\!\langle\mu_i^q:\text{all }i\rangle\!\rangle
 \quad\text{is finite}.
\]
Set $J=J_1$. Commensuration defines
$c:G\to\operatorname{Comm}(J)$, where partial isomorphisms of
finite-index subgroups are identified upon further restriction.
Each $\mu_i$ centralizes $J_i$, so $c(\mu_i)=1$ and $c$ factors
through $Q_2$. Thus $H_0=\ker c$ is normal, finite-index and
finitely generated. Intersect $H_0$, all $J_i$, and finite-index
subgroups of $J$ centralized by each generator of $H_0$.
The result $J_0$ has finite index in $J$ and lies in $Z(H_0)$.
Its full normal core $N$ is the intersection of its finitely many
$G/H_0$-conjugates. By commensuration, this core has finite index in
every $J_i$; it is abelian because $N\subset Z(H_0)$. Choose $q\ge2$ with
$\mu_i^q\in N$ for all $i$. Since $G/N$ is a quotient of $Q_q$,
$N$ has finite index in $G$.
\end{proof}

\begin{theorem}\label{thm:reduced-trivial-group}
Let $V$ be a connected normal projective variety and $A$ an effective
reduced Weil divisor, possibly zero. Suppose $(V,A)$ is log canonical and
$K_V+A$ is principal Cartier.
Then $\pi_1(U)$, where $U=V_{\rm reg}\setminus\Supp A$, is virtually
solvable.
\end{theorem}

\begin{proof}
We induct on $n=\dim V$. In dimension one the pair is an elliptic
curve with empty boundary or $\mathbb P^1$ with two marked points.
Assume $n\ge2$ and the assertion in smaller dimensions.
Take a crepant $\QQ$-factorial dlt modification by Hacon's dlt blow-up theorem
\cite[Theorem~10.4]{fujino2011},
\[
 h:X\to V,\qquad C=h_*^{-1}A+\sum_{E\text{ exceptional}}E,
 \qquad K_X+C=h^*(K_V+A).
\]
The boundary $C$ is reduced, and $X$ is klt. At a center on $C$,
the log discrepancy of $X$ is the pair log discrepancy plus the
positive valuation of $C$, so it is positive. Away from $C$,
the dlt condition already implies that $X$ is klt.
No exceptional prime lies over $U=V_{\rm reg}\setminus\Supp A$.
A projective small birational map to a smooth target is an
isomorphism: a relatively ample divisor pushed to that target is
Cartier and has the same pullback, precluding a contracted curve.
Thus $h$ is an isomorphism over $U$, and the remainder of
$X_{\rm reg}\setminus C$ over $V\setminus U$ has codimension at
least two. Consequently
$\pi_1(U)=\pi_1(X_{\rm reg}\setminus C)$.

If $C=0$, then $X$ is projective klt and $K_X\sim0$.
The ordinary group $\Gamma=\pi_1(X_{\rm reg})$ is virtually abelian
by Theorem~1.3 of \cite{fu2026}.
For the ordinary universal cover of $X_{\rm reg}$,
\cite[Proposition~4.1 and Theorem~6.1]{fu2026} give a proper
geodesic RCD completion $\mathcal M$ and a proper cocompact
action of $\Gamma$. The maximal metric splitting
used in Lemma~7.1 is $\mathcal M=\mathbb R^r\times C_0$,
where $C_0$ is proper geodesic and contains no line.
By the triangle inequality for the product metric, the coordinate
projections of a line are constant-speed geodesics. Its $C_0$ coordinate
is therefore constant.
Euclidean leaves are the classes of points that coincide
or lie on a line. An isometry has the form
$(x,c)\mapsto(A_cx+b_c,\sigma(c))$ on these leaves.
Comparing distances for $(x,c),(x,d)$ for all $x$ shows that $A_c=A_d$.
For arbitrary $(x,c),(y,d)$, the same comparison shows that $b_c=b_d$
and that $\sigma$ is an isometry.
Thus cocompactness descends to $C_0$. If $C_0$ were unbounded,
move the midpoints of minimizing segments of length $2j$ into
one compact set. Properness and Arzela--Ascoli give a subsequence
converging on every bounded parameter interval to a metric line,
a contradiction. Hence $C_0$ is compact.
The remaining argument of \cite[Lemma~7.1]{fu2026} now gives
a finite-index free abelian subgroup of the deck group $\Gamma$.

Suppose now $C\ne0$. Then
$K_X\cdot H^{n-1}=-C\cdot H^{n-1}<0$ for ample $H$.
There is a $K_X$-MMP with ample scaling ending in a Mori fiber
space $p:Y\to Z$ \cite[Corollary~1.3.2]{birkar2010}, preserving klt and
$\QQ$-factorial models \cite[Lemma~4.8.13 and Propositions~4.8.14,
4.8.16]{fujino2014a}. Push $C$ through the program. The equality
$K_X+C=\operatorname{div}(g)$ pushes to the same principal equality
on every model, and its pullbacks agree on common resolutions.
Thus the transformed reduced pairs remain lc. The endpoint boundary
$D$ is nonzero, since $-K_Y$ is ample on positive-dimensional fibers.
The common-open calculation of
Proposition~\ref{prop:ordinary-meridian-comparison}, which only uses
unchanged weights, gives
\begin{equation}\label{eq:ordinary-endpoints-mori-quotient}
 \pi_1(Y_{\rm reg}\setminus D)\twoheadrightarrow\pi_1(U).
\end{equation}

If $Z$ is a point, then $\rho(Y)=1$ and $-K_Y$ is ample.
By Lemma~\ref{lem:ordinary-endpoints-principal-adjunction}, the boundary
normalizations carry reduced lc pairs with principal adjoints.
Their groups are virtually solvable by induction. The endpoint group is
therefore virtually abelian by Theorem~\ref{thm:principal-fano-collars},
and the claim follows from
\eqref{eq:ordinary-endpoints-mori-quotient}.

\emph{A positive-dimensional Mori base.}
Suppose $0<b=\dim Z<n$. Resolve $(Y,D)$ and its full deleted set,
isomorphically over $U_Y=Y_{\rm reg}\setminus D$ and its smooth
one-branch boundary stratum. Keep every exceptional divisor deleted.
Over a connected smooth dense $Z^o\subset Z$, arrange that $Y\to Z$
is flat with normal fibers, and the resolved total space and every
horizontal boundary stratum are smooth over $Z^o$; omit vertical strata.
Use generic flatness for $Y\to Z$ and generic smoothness for the
resolved total space and all dominant components of the finitely
many closed SNC boundary intersections. These maps are proper,
so we can remove the images of their nonsmooth loci and all vertical
components to obtain one common dense base open. Restricting the
principal discrepancy formula to a general fiber gives
\[
 (F,D_F)\text{ lc},\quad F\text{ klt},\quad
 K_F+D_F\text{ principal Cartier},\quad D_F\ne0.
\]
Here $D_F$ is reduced, $-K_F$ is ample, and $\dim F=n-b<n$.
Shrink also so that $p$ is smooth on $Y_{\rm reg}$. A regular
fiber point of a flat map to a smooth base is regular in the total
space (lift local parameters), so $F_{\rm reg}=F\cap Y_{\rm reg}$.
The full SNC boundary, including all exceptional components,
defines a finite Whitney stratification of the resolution.
Every stratum is submersive over $Z^o$ by the preparation above.
The proper resolved family therefore admits local products
preserving every stratum by
\cite[Corollaries~10.2--10.3 and Proposition~11.1]{mather1970}.
Their restrictions to the entire complement, identified with
$U_Y^o$, give a bundle $U_Y^o\to Z^o$ with connected fiber
$U_F=F_{\rm reg}\setminus D_F$.

Choose a horizontal component $D_i$; it exists because $D_F\ne0$.
By adjunction and induction in dimension $n-1$, its whole tube group
$\pi_1(P_i)$ is virtually solvable. Restrict its smooth embedded base over $Z^o$
to a connected dense open $S_i$ and choose a fiberwise punctured
tube $P_i^o$ whose image lies in that of $P_i$. After further
stratification, $S_i\to Z^o$ is a locally trivial family with
$e_i>0$ connected components per fiber, each dense in a smooth
boundary stratum of $D_F$.
Put
\[
 E^o=\pi_1(U_Y^o),\quad L^o=\pi_1(Z^o),\quad
 K^o=\im\pi_1(U_F),\quad H_i^o=\im\pi_1(P_i^o).
\]
The bundle induces a surjection $q:E^o\twoheadrightarrow L^o$ with
kernel $K^o$.
Transitive monodromy on the $e_i$ components and the circle-bundle
sequence give $[L^o:q(H_i^o)]=e_i$.
By induction, $\pi_1(U_F)$ is virtually solvable; by
Theorem~\ref{thm:reduced-collar-detection}, every whole fiber
boundary tube has finite-index image. Dense-open restrictions
preserve this image. Its image in $E^o$ lies in
$K^o\cap H_i^o$, so
\[
 [E^o:H_i^o]=[L^o:q(H_i^o)]\,[K^o:K^o\cap H_i^o]<\infty.
\]
Since $E^o\twoheadrightarrow\pi_1(U_Y)$, the whole collar image
has finite index there and is virtually solvable. The endpoint group
and its quotient are therefore virtually solvable.
\end{proof}

\begin{theorem}\label{thm:numerically-trivial-group}
Let $S$ be a connected normal projective complex variety, let $\Delta$
be a standard boundary, and suppose that $(S,\Delta)$ is log canonical
and $L=K_S+\Delta$ is rational Cartier and numerically trivial.
Then $\Orbgrp S\Delta$ is virtually solvable.
\end{theorem}

\begin{proof}
Gongyo's abundance theorem for a numerically trivial lc adjoint
\cite[Theorem~3.1]{gongyo2010} gives $L\sim_{\QQ}0$.
Choose the least $N>0$ with $NL=\operatorname{div}(f)$ and normalize
$S$ in $\CC(S)(t)$, $t^N=f$, obtaining $p:Y\to S$.
Its degree is $N$: the extension is Galois with group in $\mu_N$;
if its degree were $e<N$, then $t^e\in\CC(S)$ and
$eL=\operatorname{div}(t^e)$, contradicting minimality.
At a prime $D$, putting $a_D=\operatorname{ord}_D(f)$, the local
Kummer equation gives ramification index
\[
 e_D=N/\gcd(N,a_D).
\]
At a finite weight $m_i$, one has $m_i\mid N$ and
$a_{D_i}\equiv(N/m_i)(m_i-1)\pmod N$, hence $e_{D_i}=m_i$.
At unmarked and infinite-weight primes the coefficient of $L$ is
integral, so $e_D=1$.

Let $D_\infty$ be the sum of infinite-weight components and put
$B_Y=(p^{-1}D_\infty)_{\rm red}$. The crepant upstairs coefficient
above a coefficient $c$ with ramification index $e$ is
$ec-(e-1)$: it is zero at finite weights and one at infinite weights.
Thus
\[
 K_Y+B_Y=p^*L=\operatorname{div}_Y(t).
\]
On normalized models carrying any divisorial valuation, Hurwitz gives
$a(F',Y,B_Y)+1=e(F'/F)(a(F,S,\Delta)+1)\ge0$.
Hence $(Y,B_Y)$ is lc with principal Cartier adjoint.

Put $U_S=S_{\rm reg}\setminus\Supp\Delta$ and $U_Y=p^{-1}(U_S)$.
Locally on $U_S$, write the Cartier divisor $L=K_S$ as
$\operatorname{div}(g)$. Then $t/g$ is an $N$-th root of
the unit $f/g^N$, so the cover is etale there. It is connected
because $Y$ is the normalization in a field and $U_Y$ is dense.
The coefficient calculation gives
$p^*\Delta=\operatorname{Ram}(p)+B_Y$. Both boundaries are standard,
including their coefficient-one components. Thus the finite map
between normal varieties is compatible in the sense of
\cite[Definition~2.6 and Remark~2.7]{gachet2025}; its orbifold Galois
sequence \cite[Proposition~2.9]{gachet2025} is
\[
 1\to\pi_1(Y_{\rm reg}\setminus\Supp B_Y)
   \to\Orbgrp S\Delta\to\ZZ/N\ZZ\to1.
\]
The first group is virtually solvable by
Theorem~\ref{thm:reduced-trivial-group}; its finite-index overgroup
has the same property.
\end{proof}

\section{Groups of generalized antinef pairs}\label{sec:12-bounded-special}

We prove that the meridian groups of effective standard subboundaries
of projective generalized log canonical pairs with antinef adjoint
are virtually solvable (Theorem~\ref{thm:antinef-group}). Their FAb
quotients are therefore finite. In Section~\ref{sec:13-ordinary-minimum},
this finiteness, combined with source minimality, excludes lower-dimensional
bases and forces the logarithmic adjoint to be big.

The proof is by induction on dimension. In the Mori fiber space step,
the generalized canonical bundle formula introduces a nef moduli
b-divisor on the base, so the induction requires generalized pairs.

\begin{definition}\label{def:generalized-pair}
A rational nef b-divisor $\mathbf M$ on a normal projective variety $T$
is the pullback closure of a nef rational Cartier divisor on a
projective birational model. Write $M_V$ for its trace on a model $V$.
A generalized pair $(T,D,\mathbf M)$ has effective rational boundary
$D$ and rational Cartier $K_T+D+M_T$. It is generalized log canonical
if, on every projective birational model $e:V\to T$, the divisor
$D_V$ defined by
\[
 K_V+D_V+M_V=e^*(K_T+D+M_T)
\]
has all coefficients at most one.
\end{definition}

The standard subboundary in the following theorem specifies the
meridian weights in \eqref{eq:meridian-group}.

\begin{theorem}\label{thm:antinef-group}
Let $(T,D,\mathbf M)$ be a connected projective generalized log canonical
pair with rational nef b-divisor as in
Definition~\ref{def:generalized-pair}. Assume
$-(K_T+D+M_T)$ is nef. For every effective standard rational
subboundary $\Theta\le D$, the group $\Orbgrp T\Theta$ is
finitely generated and virtually solvable, and every FAb quotient of
this group is finite.
\end{theorem}

\begin{proof}
We induct on $n=\dim T$. For $n=1$,
$T$ is smooth, the nef b-divisor descends to $T$, and
\[
 \deg(K_T+\Theta)\le\deg(K_T+D+M_T)\le0.
\]
By the orbicurve classification, including its punctured and
nonuniformizable cases, $\Orbgrp T\Theta$ is finitely generated
and virtually abelian
\cite[Propositions~5.3--5.6 and~5.11]{behrend2005}.
For positive genus the degree inequality forces genus one and
empty boundary, giving $\ZZ^2$. In genus zero it gives exactly
the spherical or Euclidean signatures. The spherical groups are
finite, including a two-point football of weights $a,b$ with
group $\ZZ/\gcd(a,b)\ZZ$; the Euclidean groups have a lattice
of finite index. Infinite weights are punctures in these same
presentations.
Suppose $n\ge2$ and the theorem holds in all smaller dimensions.
Finite generation follows from a finite triangulation of a projective
compactification compatible with the deleted singular and boundary
loci~\cite[Theorems~3.4.4 and~9.2.1]{bochnak1998}.
As in Theorem~\ref{thm:canonical-albanese}(i), the complement has
a finite CW model; taking the meridian quotient preserves finite
generation.

\emph{A dlt model with trivial generalized adjoint.}
Put $A_0=-(K_T+D+M_T)$ and
$\mathbf M^+=\mathbf M+\overline{A_0}$, where the bar denotes
pullback closure. The sum is nef on a common determining model.
Adding $A_0$ changes no generalized discrepancy, and
$K_T+D+M_T^+=0$. Take a crepant generalized dlt modification
\cite[Theorem~2.9 and Remark~2.31]{filipazzi2023},
\[
 e:V\to T,\qquad
 K_V+B+M_V^+=e^*(K_T+D+M_T^+)\sim_{\QQ}0,
\]
with $V$ $\QQ$-factorial and $B\ge0$.
On a resolution
$f:W\to V$ where the nef part descends, put
$E=f^*M_V^+-M_W^+$. This divisor is exceptional and $-E$ is
$f$-nef, so negativity gives $E\ge0$.
The generalized crepant boundary is
$B_W^{\rm gen}=B_W^{\rm ord}+E$. Hence the ordinary pair is lc,
and every ordinary lc place is also a generalized lc place.
The generalized dlt condition implies that $(V,B)$ is log smooth
near the generic point of each such center. Resolving the complement of
that open
introduces only discrepancies greater than $-1$, which is the
ordinary dlt condition. Thus $(V,B)$ is dlt.
Define $\Gamma\le B$ by giving every floor component coefficient
one and otherwise retaining the coefficient of the strict transform
of $\Theta$. Thus $(V,\Gamma)$ is dlt and standard.

There is an epimorphism
\begin{equation}\label{eq:antinef-dlt-group-quotient}
 \Orbgrp V\Gamma\twoheadrightarrow\Orbgrp T\Theta.
\end{equation}
Indeed, remove the exceptional image, of codimension at least two,
from the downstairs regular locus and remove any additional floor
primes. The resulting opens agree. Every extracted divisor is deleted
upstairs. Restoring an additional floor prime downstairs imposes its
first meridian or the finite power prescribed by its original weight.

\emph{The pseudoeffective case.}
The trace $M_V^+$ is pseudoeffective: approximate its nef determining
divisor by ample rational classes with effective representatives and
push them to $V$. Numerical pushforward is well defined because $V$
is $\QQ$-factorial; the difference between a numerically trivial
divisor upstairs and the pullback of its pushforward is exceptional
and relatively numerically trivial, hence zero by applying the negativity
lemma~\cite[Lemma~2.3.26(i)]{fujino2014a} to both signs.
It follows that the pushforward of a numerically trivial class
is numerically trivial, so pushforward is a continuous linear map
on the finite-dimensional numerical divisor spaces.
Consequently
\[
 (K_V+\Gamma)+(B-\Gamma)+M_V^+\sim_{\QQ}0
\]
shows that if $K_V+\Gamma$ is pseudoeffective, both it and its
negative are pseudoeffective, and hence $K_V+\Gamma\equiv0$.
Here the pseudoeffective cone is pointed: on a smooth resolution this
follows from duality with movable curve classes, which span the numerical
curve space~\cite[Theorems~0.2 and~2.2]{boucksom2004}.
Indeed pullback sends both signs to the pseudoeffective cone on
the resolution. They pair to zero with all movable curves;
the complete-intersection classes $H^{n-2}A$, for fixed ample
$H$ and varying ample $A$, already span the numerical curve
space by the numerical Hodge index pairing. The pulled-back
class is zero, and pullback on numerical Cartier classes is
injective.
Theorem~\ref{thm:numerically-trivial-group} now gives virtual
solvability, which passes through
\eqref{eq:antinef-dlt-group-quotient}.

\emph{The Mori fiber space.}
Otherwise choose an ample Cartier class $A$ and a rational $t>0$
large enough that $K_V+\Gamma+tA$ is nef. For sufficiently
divisible $m\gg t$, take a general $A_m\in|mA|$ and put
$H=(t/m)A_m$. This divisor has class $tA$, and the pair
$(V,\Gamma+H)$ is lc because the coefficient is small and the divisor
is transverse to the dlt strata.
The nonpseudoeffective dlt program
with scaling terminates in a Mori fiber space
$p:(V',\Gamma')\to Z$~\cite[Section~2.20]{fujino2014b}.
Its models remain $\QQ$-factorial dlt, no prime is extracted, and
surviving weights are unchanged. Meridian comparison
(Proposition~\ref{prop:ordinary-meridian-comparison}) gives
\begin{equation}\label{eq:antinef-endpoint-group-quotient}
 \Orbgrp{V'}{\Gamma'}\twoheadrightarrow\Orbgrp V\Gamma
       \twoheadrightarrow\Orbgrp T\Theta.
\end{equation}

The generalized adjoint remains crepant throughout this program.
Fix one rational top form for the canonical divisors and write
\[
 q(K_V+B+M_V^+)=\operatorname{div}(\varphi)
\]
for some integer $q>0$. Push this equality through divisorial steps;
flips preserve it in codimension one. On a common resolution all
pullbacks equal the same principal divisor. Their generalized
boundaries therefore coincide. The endpoint is generalized lc and
\[
 K_{V'}+B'+M_{V'}^+\sim_{\QQ}0,\qquad B'\ge\Gamma'.
\]
\emph{A Mori base of dimension zero.}
Suppose $\dim Z=0$. The endpoint then has Picard number one.
If $\lfloor\Gamma'\rfloor=0$, the pair is klt with finite weights
and ample negative adjoint. Its ordinary topological meridian group
is therefore finite by Braun's theorem~\cite[Theorem~2]{braun2020}.
If the floor is nonempty, choose a component $S$ and apply
Theorem~\ref{thm:positive-standard-collar}:
its collar image $J_S$ has finite index, its positive meridian
$h_S$ is central, and
\[
 \Orbgrp S{\Theta_S}\twoheadrightarrow J_S/\langle h_S\rangle,
 \qquad \Theta_S=\operatorname{Diff}_S(\Gamma'-S).
\]
The adjunction pair is normal projective dlt and standard, with
$-(K_S+\Theta_S)$ ample. Use the pullback closure of this ample
divisor as the nef part of a generalized pair. Since its dimension is
$n-1$, the induction hypothesis implies that $\Orbgrp S{\Theta_S}$
is virtually solvable. The inverse image in
$J_S$ of a solvable finite-index subgroup of
$J_S/\langle h_S\rangle$ is solvable, since its kernel is central
cyclic. It has finite index in the endpoint group, which is therefore
virtually solvable.

\emph{A positive-dimensional base.}
Suppose $0<\dim Z<n$. A sufficiently general fiber $F$ is normal
and connected, and generic flatness, generic smoothness, and
restriction of a log resolution give
\[
 F_{\mathrm{reg}}=F\cap V'_{\mathrm{reg}},\qquad
 K_F+\Gamma'_F=(K_{V'}+\Gamma')|_F,
\]
with $(F,\Gamma'_F)$ lc and $-(K_F+\Gamma'_F)$ ample.
Horizontal coefficients restrict with their original weights. Use
the pullback closure of the ample divisor $-(K_F+\Gamma'_F)$ as
the nef part. By induction, $\Orbgrp F{\Gamma'_F}$ is virtually solvable.
Let $N$ be its actual image in $\Orbgrp{V'}{\Gamma'}$.
The nef part just chosen descends to $F$, so its generalized
discrepancies equal the ordinary discrepancies.

The generalized canonical bundle formula
\cite[Definition~2.19 and Theorem~2.20]{filipazzi2023} associates
a discriminant and a nef moduli b-divisor to a rational generalized
lc-trivial fibration. Apply it to
$p:(V',B',\mathbf M^+)\to Z$, for which the divisor pulled back
from the base is zero.
To verify the rank condition, take a resolution where the nef part
descends. The rounded modified discrepancy divisor is effective
and exceptional: nonexceptional coefficients in $[0,1)$ round
to zero, while coefficient-one discrepancies are omitted in
$\mathbf A^*$. Its associated sheaf pushes to $\mathcal O_{V'}$
by normality, then to $\mathcal O_Z$ by connectedness of fibers.
The direct image therefore has rank one, so the formula gives
\[
 K_Z+B_Z+N_Z\sim_{\QQ}0
\]
with rational nef b-divisor $\mathbf N$. For a prime $P\subset Z$,
let $t_P$ be the generic generalized lc threshold of $p^*P$, so
$\operatorname{coeff}_P B_Z=1-t_P$.

For every prime $E\subset V'$ dominating $P$, let $e_E$ be its
pullback multiplicity and $w_E$ its weight in $\Gamma'$, including
$w_E=1$ when unmarked. Put
\[
 a_P=\gcd\{e_Ew_E:w_E<\infty,\ p(E)=P\},\qquad
 \Xi=\sum_P(1-1/a_P)P,
\]
with empty gcd infinity. Every prime $P$ has a divisorial
inverse-image component by properness and dominance. Outside the finitely many divisorial
components of the discriminant and the vertical marked locus,
the generic fibers over $P$ are reduced and an unmarked component
has multiplicity one. Thus $a_P=1$ there, and $\Xi$ has finite
support.
Global generalized lc gives $t_P\ge0$. If
$b_E=\operatorname{coeff}_E B'$, the discrepancy after adding
$t p^*P$ is $1-b_E-te_E$. Since $B'\ge\Gamma'$, this yields
\[
 t_P\le\frac{1}{e_Ew_E}\quad(w_E<\infty),\qquad
 t_P=0\quad(w_E=\infty).
\]
Hence $B_Z\ge0$ and $\Xi\le B_Z$, also when the set defining
the gcd is empty. On each higher base model the same computation uses the
crepant generalized total pair, whose thresholds are nonnegative.
All higher discriminant coefficients are at most one; consequently
$(Z,B_Z,\mathbf N)$ is generalized lc. By induction
$\Orbgrp Z\Xi$ is virtually solvable.

\emph{The group extension.}
Choose a dense regular base open over which a relative log resolution
and its boundary and exceptional strata are smooth. The marked proper
family is a differentiable bundle on its regular boundary complement.
Horizontal finite meridian powers are fiber relations. The group of
this open surjects onto the endpoint group. The bundle sequence has
normal kernel equal to the image of the whole fiber-complement group.
Quotienting by
the horizontal finite powers keeps their normal closure inside
that kernel; its image is a quotient of
$\Orbgrp F{\Gamma'_F}$. The dense-open surjection to the
endpoint maps this normal kernel onto $N$. Thus $N$ is normal
and virtually solvable, and the group of the chosen base open surjects onto the
endpoint group modulo $N$.
At a prime $E$ over $P$, a base parameter is a unit times
a normal parameter to the power $e_E$. Restoring an unmarked prime
imposes $\delta_P^{e_E}$, a finite weight imposes
$\delta_P^{e_Ew_E}$, and an infinite weight imposes no relation.
B\'ezout's identity therefore gives an epimorphism
\[
 \Orbgrp Z\Xi\twoheadrightarrow
              (\Orbgrp{V'}{\Gamma'})/N.
\]
Divisors over subsets of the base of codimension at least two can
impose further relations. Base singularities and divisors of infinite
weight remain deleted.

The normal subgroup $N$ and the quotient
$(\Orbgrp{V'}{\Gamma'})/N$ are virtually solvable.
Apply \cite[Lemma~5.5]{dinh2012tits} with
$A=\Orbgrp{V'}{\Gamma'}$ and $B=N$ to conclude that
$\Orbgrp{V'}{\Gamma'}$ is virtually solvable. The epimorphisms in
\eqref{eq:antinef-endpoint-group-quotient} finish the induction.

Let $Q$ be a FAb quotient and choose a solvable finite-index
subgroup $H\le Q$. Since every finite-index subgroup of $Q$ has
finite abelianization, $H'$ still has finite index in $Q$. Repeating
this argument down the finite derived series shows that the trivial
subgroup has finite index. Thus $Q$ is finite.
\end{proof}

\section{Minimal quotients and the proof of the main theorem}\label{sec:13-ordinary-minimum}

Assuming that the canonical kernel is infinite, we construct a quotient
of minimal positive relative dimension whose weighted adjoint is big
on very general projective fibers of every prepared model. We use the
Iitaka leaves to derive a contradiction with specialness, then finish the proof of
Theorem~\ref{thm:introduction-main} by applying the finite-kernel criterion.

\subsection{Construction of the minimal quotient}\label{subsec:minimal-quotient-construction}

\begin{definition}\label{def:prepared-ordinary-model}
Let $Q$ be an algebraic principal torus torsor over a complex abelian
variety. Let $Y$ be smooth connected quasi-projective and
$f:Y\to Q$ a dominant morphism inducing an epimorphism of
fundamental groups. Let $c:Y\dashrightarrow Z$ be a dominant rational
map with geometrically connected general fiber, with $Z$ smooth and
$g:Z\to Q$ projective surjective with connected general fibers and
$gc=f$. Put $G=\pi_1(Y)$, let $B_c\triangleleft G$ be the image in $G$
of the fundamental group of a smooth resolved general $c$-fiber, and put
$R=G/B_c$. We call a resolution over $Y$ \emph{retained} if every
exceptional divisor over the open variety $Y$ remains in the resolution.

A \emph{prepared model} consists of a smooth $Z_m$ projective
over $Q$ and birational over $Q$ to $Z$, a smooth retained proper
birational $\mu_m:M_m\to Y$, and a morphism $c_m:M_m\to Z_m$
representing $c$, with $g_mc_m=f\mu_m$. Choose a dense open
$Z_m^\circ=Z_m\setminus(D_m\cup\Sigma_m)$ over which $c_m$ is a smooth
surjective topological fiber bundle with connected fibers, with $D_m$
a reduced SNC divisor and $\operatorname{codim}\Sigma_m\ge2$.
Require that the quotient $G\to R$ induce an epimorphism
$\eta_m:\pi_1(Z_m^\circ)\to R$, compatible on common dense opens and
common prepared refinements. For a prime $D\subset D_m$, let $m_D$
be the actual order of $\eta_m(\mu_D)$ in this fixed $R$, allowing
$m_D=\infty$. Set
$\Delta_m=\sum_{m_D>1}(1-1/m_D)D$.
Other interior primes have weight one. All relations from retained
source divisors, including those mapped into base codimension at least
two, are imposed before these orders are defined.
\end{definition}

The inclusion of a dense algebraic open in a smooth variety induces a
surjection on fundamental groups
\cite[Proposition~2.9.1]{kollar1993}. If the complement has codimension
at least two, loops and homotopy disks can be moved off it by general
position, so the induced map is an isomorphism. By transverse disks and
van Kampen, filling divisors
kills the normal closure of their generic meridians.
A proper birational morphism between smooth varieties is an isomorphism
over an open whose complement downstairs has codimension at least two.
It therefore induces an isomorphism on fundamental groups. By
\cite[Corollary~2.9.3]{kollar1993}, for a rational fibration with
geometrically connected general fiber, the whole resolved general-fiber
image is normal; common retained resolutions preserve that image by
the same birational comparison. All resolved fibers below are whole
fibers, with every divisor over the original open retained.

\begin{theorem}\label{thm:ordinary-source-minimum}
Let $X$ be a smooth connected special complex quasi-projective
variety, and let $f_0:Y_0\to Q_0$ be its canonical second Albanese map
from Theorem~\ref{thm:canonical-albanese}, with kernel $P_0$.
Assume $P_0$ is infinite. Then there are a connected finite etale cover
$\rho:Q\to Q_0$, a connected pullback $Y=Y_0\times_{Q_0}Q$, the lift
$f:Y\to Q$, and a connected rational fibration $c:Y\dashrightarrow Z$
over $Q$, with $g:Z\to Q$ projective surjective and connected fibers,
having the following properties.
\begin{enumerate}[label=(\roman*)]
\item Identify $G=\pi_1(Y)$ with its subgroup in $\pi_1(Y_0)$.
Then $P_0\subseteq G$. Put $P=P_0$ and let $B_c$ be the image in $G$ of
the fundamental group of a resolved general $c$-fiber. One has
$B_c\triangleleft G$, $B_c\le P$,
and an exact sequence
\[
 1\longrightarrow\Pi=P/B_c\longrightarrow R=G/B_c
 \longrightarrow\pi_1(Q)\longrightarrow1,
\]
where $\Pi$ is infinite and FAb. Set $\ell=\dim Z-\dim Q$.
Then $1\le\ell\le\dim Y_0-\dim Q_0$.
\item Prepared models and common prepared refinements exist.
For every such model, a very general point $u\in Q$
gives a smooth connected open fiber $S=f^{-1}(u)$ and a smooth
connected projective $\ell$-fold $V=g_m^{-1}(u)$ with SNC standard
boundary $\Delta_V=\Delta_m|_V$. The induced rational map
$a_m:S\dashrightarrow V$ is dominant with geometrically connected
general fiber. The surjection
$\nu:\pi_1(S)\to P\to\Pi$ is compatible on a dense source open
with a surjection $\Orbgrp V{\Delta_V}\to\Pi$ induced by
$S\dashrightarrow V$. Each meridian has exactly its specified order
in $\Pi$, including infinite order.
\item For each such $S$, every smooth projective $W$ with
$1\le\dim W<\ell$, and every connected rational fibration
$h:S\dashrightarrow W$, the image in $\Pi$ of the fundamental group of
a smooth connected sufficiently general resolved $h$-fiber has finite
index. Its map is through a retained resolution to $S$ and then $\nu$.
\end{enumerate}
The same cover, $c$, $B_c$, $R$, $\Pi$ and $\ell$ work for every
prepared model; the very general point may depend on the
model. The target $Q$ is again an algebraic
principal torus torsor over an abelian variety.
\end{theorem}

\begin{proof}
\emph{Good fibers.}
Put $G_0=\pi_1(Y_0)$ and $E_0=\pi_1(Q_0)$. By
Theorem~\ref{thm:canonical-albanese} and
Proposition~\ref{prop:canonical-kernel-covers}, $P_0$ is FAb and is the
image of the fundamental group of a general $f_0$-fiber.
Choose a nonempty connected open $Q_0^*\subset Q_0$ on which the
whole fibers are smooth and geometrically connected and have this
same image $P_0$. Applying Theorem~2.2 of \cite{achinger2019}
to $f_0$, we obtain a bundle locus on which transport preserves $P_0$.

Take a smooth compactification of $Y_0$ projective over $Q_0$, with
full SNC boundary.
On each of the countably many relative Hilbert schemes, stratify the
locus of geometrically integral cycles meeting $Y_0$. Resolve the
generic universal cycle and its boundary, spread this resolution over
a nonempty parameter open, and shrink until its fibers are smooth
geometrically connected proper birational resolutions of the cycles.
Noetherian induction on the remaining parameter subsets gives
countably many diagrams
\[
 p_j:U_j\to T_j,\qquad e_j:U_j\to Y_0,\qquad
 f_0e_j=t_jp_j,
\]
whose fibers are smooth connected resolutions of whole closed integral
subvarieties of actual $f_0$-fibers. Every such subvariety occurs by
closing it in the compactification. Delete exactly the inverse image of
its boundary; each fiber evaluation remains proper. Apply
\cite[Theorem~2.2]{achinger2019} separately to each $p_j$, and refine
to smooth connected parameter strata on which the $p_j$ are
topological bundles. Under $e_j$, a transported loop cylinder becomes
a homotopy in $Y_0$, with base point moving along the evaluated path.
The subgroup images are therefore conjugate by that path. By the
bundle sequence, the fiber image is normal in $\pi_1(U_j)$, so its
image under evaluation is normal in $e_{j*}\pi_1(U_j)$.

If $A_j=\overline{e_j(U_j)}\ne Y_0$, constructibility of
$f_0(Y_0\setminus A_j)$ gives a nonempty open
$O_j\subset f_0(Y_0\setminus A_j)$. Put
\[
 Q_0^{\mathrm{good}}=Q_0^*\cap\bigcap_{A_j\ne Y_0}O_j.
\]
Its complement is a countable union of proper closed subsets, so it
meets every nonempty open of $Q_0$. A whole fiber over this set can
be contained in $A_j$ only if $A_j=Y_0$. Call these fibers \emph{good}.

\emph{Minimizing the base dimension.}
Since $P_0$ is infinite, $s=\dim Y_0-\dim Q_0>0$.
For every good $f_0$-fiber $S_0$ and every
connected rational fibration $h:S_0\dashrightarrow W$ to a smooth
projective variety, let
$C_h\triangleleft P_0$ be its actual resolved fiber image; normality
follows from \cite[Corollary~2.9.3]{kollar1993} and the
retained-resolution comparison above. Define
\[
 \ell=\min\{\dim W:P_0/C_h\text{ is infinite}\}.
\]
A compactification of $S_0$ is a candidate of dimension $s$, so
$1\le\ell\le s$.

Choose a minimizing fibration on a good fiber $S_0$ and put $C=C_h$.
Its leaves have dimension $s-\ell$, and a dense flat family of their
closures occurs in one of the resolved Hilbert strata above. Its
evaluation image contains the selected whole $S_0$ in its closure.
Since this fiber is good, the evaluation map is dominant onto $Y_0$.
Common retained resolutions identify its fiber images with $C$ up to
transport. By the bundle sequence, $C$ is normal in
$L=e_{j*}\pi_1(U_j)$, which has finite index in $G_0$ by
\cite[Proposition~2.9.1]{kollar1993}. Therefore
$G=N_{G_0}(C)$ has finite index in $G_0$. Since $C\triangleleft P_0$,
one has $P_0\subseteq G$.

Put $J=f_{0*}(G)\le E_0$ and take the connected finite etale cover
$\rho:Q\to Q_0$ corresponding to $J$
\cite[Expos\'e~XII, Theorem~5.1]{grothendieck2004}. Since
$f_{0*}^{-1}(J)=G$, the pullback $Y=Y_0\times_{Q_0}Q$ is connected
with fundamental group $G$. Its projection $f:Y\to Q$ induces
$G\twoheadrightarrow J=\pi_1(Q)$ with kernel $P=P_0$. Moreover,
\[
 f^{-1}(u)\simeq f_0^{-1}(\rho(u))\qquad(u\in Q),
\]
compatibly with the fiber group images in $P_0$. Thus every fiber over
$\rho^{-1}(Q_0^{\mathrm{good}})$ participates in the same minimum.
The evaluation lifts to $Y$ because $L\subseteq G$, and remains
dominant. Since $p_{j*}$ is surjective, $t_{j*}\pi_1(T_j)=f_{0*}(L)$
lies in $J$; hence the parameter map also lifts to $Q$. The lifted
family has the same fiber image $C\triangleleft G$.
By Proposition~\ref{prop:finite-cover-torus-torsor}, $Q$ is a torus
torsor over an abelian variety.

\emph{Globalizing the minimizing family.}
Apply the Shafarevich theorem
\cite[Definition~3.2 and Corollary~3.5.1]{kollar1993} to
$C\triangleleft G$. We obtain a rational map
$c_0:Y\dashrightarrow Z_0$ with connected fibers and countably many
proper closed subsets $D_i\subset Y$. A closed integral subvariety not contained in their
union is contracted precisely when its normalization has finite
image in $G/C$. The same criterion holds for its whole resolution,
whose map to the normalization induces a surjection on fundamental
groups by \cite[Proposition~2.9.1]{kollar1993}.
For each $D_i$, smoothness of $p_j$ and dominance of $e_j$ imply that
the parameters whose whole member lies in $D_i$ form a proper closed
subset. Very general members therefore satisfy this
criterion. Their image in $G/C$ is trivial, so $c_0e_j$ is constant
on them. It factors rationally through $p_j$: the image of
$(p_j,c_0e_j)$ is generically finite over $T_j$, and its field is
$\CC(T_j)$ because the geometric generic $p_j$-fiber is integral.

Take a retained resolution over $Y$ of the graph of $c_0$ with
projective target, and let $B_c$ be its whole general-fiber image.
A very general member maps to a sufficiently general point of this
target, by dominance.
Resolve its proper graph into the retained source. Its image lies
in the whole $c_0$-fiber, while proper birational invariance preserves
its image $C$ in $G$. Thus $C\le B_c$. Separately, a very general
whole $c_0$-fiber has a closed integral image in $Y$, is proper
birational over that image, and avoids being contained in any $D_i$.
The last assertion follows by applying constructibility to the image
of the inverse image of $Y\setminus D_i$ in the target. The
Shafarevich criterion on this image gives
$B_c/(B_c\cap C)$ finite. Consequently $B_c/C$ is finite.

A morphism to $Q$ with finite fundamental-group image is constant.
Indeed, a nonconstant projection to its abelian base would be nonconstant
on a smooth curve. Completing that curve gives a map with trivial
image in the torsion-free abelian lattice, hence a lift to a vector
space, which is constant. Once the image lies in a torus fiber, a
nonconstant character on a curve has a zero or pole on its completion
and a loop of nonzero winding. The torus lattice injects into
$\pi_1(Q)$, excluding this case as well.

Since $C\le P$, the general $c_0$-fiber has finite image in $\pi_1(Q)$,
so $f$ factors through $c_0$. A smooth projective relative
compactification of the target graph in $\overline Z_0\times Q$,
followed by resolution, gives $c:Y\dashrightarrow Z$ and
$g:Z\to Q$. Its source graph is closed in $Y\times_Q Z$, hence
proper over $Y$. Common retained resolutions preserve $B_c$.
Since $\CC(Q)$ is relatively algebraically closed in $\CC(Y)$, it
is also relatively algebraically closed in $\CC(Z)$. The proper
Stein factor over the normal base $Q$ is therefore trivial. Thus $g$ is
surjective with connected fibers and $B_c\le P$.

The infinitude of $P_0/C$ and finiteness of $B_c/C$ imply that
$P_0/B_c$ is infinite. Put $k=\dim Z-\dim Q$.
Since the leaves are contracted, $k\le\ell$.
On a sufficiently general $S=f^{-1}(u)$ with
$\rho(u)\in Q_0^{\mathrm{good}}$, the induced map to the smooth
projective $g$-fiber has connected general fiber and actual image
$B_c$. Under $S\simeq f_0^{-1}(\rho(u))$, this is a candidate for the
minimum, since $P_0/B_c$ is infinite. Hence $k\ge\ell$, so $k=\ell$.
The infinite group $\Pi=P_0/B_c$ is FAb as a quotient of $P_0$.
This proves (i).

\emph{Preparation and the fixed quotient.}
Fix a retained resolved morphism to a smooth target projective
over $Q$ and apply \cite[Theorem~2.2]{achinger2019} on its dense
smooth connected-fiber locus. Resolve the complement of this bundle open
on the target, then resolve the new source graph over $Y$, leaving
the bundle open unchanged. The target stays projective over $Q$,
so the graph is proper over $Y$. This constructs a prepared model.
Resolve the closure of the target correspondence in the product over
$Q$, and its retained source graph, to obtain a common refinement;
the same construction handles any finite collection.
Birational maps between retained resolutions preserve $G$ and $B_c$.
On the bundle locus, the sequence
\[
 \pi_1(F_m)\longrightarrow\pi_1(M_m^\circ)
 \longrightarrow\pi_1(Z_m^\circ)\longrightarrow1
\]
shows that $G\to R$ factors uniquely through an epimorphism
$\eta_m:\pi_1(Z_m^\circ)\to R$. On a common dense base open, the
two representations pull back on a common retained source to the
same $G\to R$. Apply the bundle theorem to this morphism over a
smaller open. Since the fibers are connected, the fundamental group
of the total space surjects onto that of the base open. The two
representations therefore agree. Their
composite to $\pi_1(Q)$ is $g_{m*}$.
If a retained source prime dominates a base prime with multiplicity
$e$, its transverse disk gives $\eta_m(\mu_D)^e=1$. Its true order
therefore divides $e$. Thus a prime of infinite order has no such
source prime above it.
Relations from divisors over base codimension at least two already
hold in $G/B_c$.

For each prepared model choose $u\in\rho^{-1}(Q_0^{\mathrm{good}})$ so that $S$,
$M_{m,u}$ and $V$ are smooth, $M_{m,u}\to S$ is proper birational,
$V$ is transverse to the dominating SNC strata and avoids the others,
and $\Sigma_m\cap V$ has codimension at least two. Generic smoothness,
the dimension formula and properness give these conditions on a dense
open, which meets $\rho^{-1}(Q_0^{\mathrm{good}})$. By the Stein
argument, $V$ is connected. A dominating exceptional divisor has
general fiber of dimension at most $\dim S-1$, so no component of
the smooth $M_{m,u}$ lies entirely in the exceptional locus. The unique component
meeting the isomorphism open is therefore all of $M_{m,u}$.
Restricting this bundle to $V^\circ$ shows that $a_m$ has geometrically
connected general fiber.
Its bundle sequence and dense-open surjectivity show that
\[
 \im(\pi_1(V^\circ)\to R)=\im(\pi_1(S)\to R)=\Pi.
\]
A transverse meridian has the same order in $\Pi\subset R$. Fill
order-one primes and codimension-two omissions and impose the specified
finite powers. This proves (ii).

Let $h$ be as in (iii). Under its isomorphism with $f_0^{-1}(\rho(u))$, the good fiber $S$
is among those used to define the minimum. Thus its
normal fiber image $C_h\le P_0$ has $P_0/C_h$ finite, and therefore
$\Pi/(C_hB_c/B_c)=P_0/(C_hB_c)$
is finite. This proves (iii).
\end{proof}

\subsection{Bigness of the weighted adjoint}\label{subsec:minimal-quotient-bigness}

\begin{theorem}\label{thm:source-minimal-bigness}
For the fixed data of Theorem~\ref{thm:ordinary-source-minimum}, the
weighted adjoint $K_V+\Delta_m|_V$ is big on a very general smooth
projective fiber $V$ of every prepared model $g_m:Z_m\to Q$.
\end{theorem}

\begin{proof}
Fix a prepared model and take $S,V,\Delta,\Pi,\nu$ as in
Theorem~\ref{thm:ordinary-source-minimum}(ii), with $n=\ell$.
Suppose a sufficiently general whole resolved fiber $T$ of a fibration
$S\dashrightarrow W$, $0<\dim W<n$, maps on a dense open through a
group $A_F$ all of whose FAb quotients are finite. By minimality in
part (iii),
\[
 I=\im(\pi_1(T)\to\Pi)\ \le\ J=\im(A_F\to\Pi)
\]
has finite index in $\Pi$. The inclusion follows from the surjection
on fundamental groups induced by a dense open of $T$.
Thus $J$ is FAb and finite, contradicting the infinitude of $\Pi$.

\emph{Minimal models and Mori fiber spaces.}
Put $L=K_V+\Delta$. The smooth rational SNC pair $(V,\Delta)$ is
$\QQ$-factorial dlt. If $L$ is
not pseudoeffective, choose an effective ample scaling divisor $A$
with $(V,\Delta+A)$ lc and $L+A$ nef. The non-pseudoeffective
case of \cite[Section~2.20]{fujino2014b} gives an $L$-negative
ample-scaling program ending at a Mori fiber space
$p:(M,C_M)\to Z$.
If $L$ is pseudoeffective, \cite[Theorem~11.1]{openai2026b} gives a
good log minimal model of this projective lc pair. The log-model
convention permits coefficient-one extracted primes
\cite[Definitions~2.7--2.8]{liu2025}. By
\cite[Lemma~2.9, (1)$\Rightarrow$(3)]{liu2025}, the existence of this
model gives a terminating ample-scaling MMP for the original $L$,
ending in an ordinary good minimal model $(M,C_M)$.
Here real semiampleness means real linear equivalence to the pullback
of an ample real divisor by a contraction.
By \cite[Lemmas~2.1.22--2.1.23]{fujino2014a}, real semiampleness
of the rational adjoint implies that a positive Cartier multiple is
globally generated after clearing denominators. Let
$p:M\to Z$ be the resulting connected semiample contraction.

In either program the pairs remain $\QQ$-factorial dlt, no prime is
extracted, and surviving coefficients are unchanged
\cite[Lemma~4.8.13 and Propositions~4.8.14 and~4.8.16]{fujino2014a}.
The scaling divisor selects the extremal rays; at each stage the
boundary is the transform of $\Delta$. Applying negativity to the
$L$-negative steps
\cite[Lemmas~2.3.26--2.3.27]{fujino2014a} and adding the resulting
divisor comparisons on a common resolution $a:T\to V$, $b:T\to M$ gives
\[
 a^*L=b^*(K_M+C_M)+E_{\rm exc},\qquad
 E_{\rm exc}\ge0\quad\text{$b$-exceptional}.
\]
The divisor $E_{\rm exc}$ is exceptional over $M$: surviving primes
keep their coefficients, and each divisorial contraction contributes
only primes contracted by the composite to $M$. Since $b$ is proper
birational, $M$ is normal, and $E_{\rm exc}$ is effective and
$b$-exceptional, the pushforward argument in the proof of
Lemma~\ref{lem:canonical-kernel-logarithmic-comparison} gives
$b_*\mathcal O_T(mE_{\rm exc})=\mathcal O_M$ in every divisible degree.
The projection formula then gives
\[
 H^0(V,mL)\simeq H^0(T,ma^*L)
 \simeq H^0(M,m(K_M+C_M)).
\]
These are the natural identifications in every divisible degree,
compatible with multiplication and all equal-degree section ratios.
Proposition~\ref{prop:ordinary-meridian-comparison} gives
compatible epimorphisms
\[
 \Orbgrp M{C_M}\twoheadrightarrow\Orbgrp V\Delta
 \twoheadrightarrow\Pi.
\]
Write $d=\dim Z$. The Mori case has $d<n$; in the semiample case,
the section comparison gives $d=\kappa(V,L)$.

\emph{Minimality excludes a lower-dimensional base.}
If $d=0$, the endpoint adjoint is anti-ample or rationally trivial.
By Theorem~\ref{thm:antinef-group} with zero nef part, its FAb
quotient $\Pi$ is finite, a contradiction. Suppose $0<d<n$. A very general
fiber $(F,C_F)$ is projective lc with standard coefficients and with
anti-ample or rationally trivial adjoint. To obtain its group
map, shrink the smooth base locus until $p$ is flat, its restriction to
$M_{\rm reg}$ is smooth, and a log resolution with its boundary strata
restricts to the fibers. Restricting discrepancies shows that the fiber pair is lc. Flatness
over the regular base and smoothness on $M_{\rm reg}$ imply
$F_{\rm reg}=F\cap M_{\rm reg}$.
Horizontal boundary components restrict reducedly, with their original
coefficients and meridians. Hence
$\Orbgrp F{C_F}\to\Orbgrp M{C_M}\to\Pi$ is defined, and
by Theorem~\ref{thm:antinef-group}, every FAb quotient of the first
group is finite.

Resolve $Z$ to a smooth projective $W$. The composite
$S\dashrightarrow W$ has geometrically connected general fiber, since
its function-field extension is a composite of regular extensions. A
whole sufficiently general resolved fiber has a dense open mapping to
$F_{\rm reg}\setminus\Supp C_F$, compatibly with $\nu$. Apply the
consequence of minimality above with $A_F=\Orbgrp F{C_F}$ to obtain a
contradiction. Thus the Mori case is impossible, and in the semiample
case $d=n$, so $L$ is big.
\end{proof}

\subsection{Iitaka leaves and the specialness obstruction}\label{subsec:iitaka-specialness}

\begin{definition}\label{def:ordinary-minimum-log-field}
Write $\kappa_{\log}(U,C)=\kappa(U,C)$ for the logarithmic dimension
of Definition~\ref{def:canonical-kernel-logarithmic-pair}. When it is
nonnegative, the \emph{field of section ratios} is the subfield of
$\CC(U)$ generated by ratios of nonzero sections of equal positive
divisible degree of the full logarithmic adjoint. Its independence of
compactification and the comparison of section rings for torus torsors
follow from
Lemma~\ref{lem:canonical-kernel-logarithmic-comparison}.
\end{definition}

For the rational logarithmic adjoints of nonnegative Iitaka dimension
used below, the full section-ratio field is the function field of their
connected Iitaka fibration. On the smooth projective compactification,
put $\Gamma=D_\infty+\overline C$, including every deleted prime.
This is a smooth geometric orbifold in the sense of
\cite[Definition~2.1]{campana2009}: every rational coefficient
$\gamma\in[0,1]$ has multiplicity $1/(1-\gamma)\in\QQ_{\ge1}\cup\{\infty\}$.
By \cite[Section~3.3, p.~836]{campana2009}, the sufficiently large
divisible complete systems have connected general fibers and define
the same rational fibration. Write $K_0$ for the function field of
its normal base, of transcendence degree $\kappa_{\log}$. Over $\CC$,
the general fibers are connected, so the finite Stein factor of a
resolved system has degree one. Thus $K_0$ is relatively algebraically
closed in the ambient function field. Each equal-degree ratio belongs to one of these
systems after a Veronese replacement, since
$s/t=st^{v-1}/t^v$ for arbitrarily divisible $v$; conversely the
ratios of any such complete system generate $K_0$. Thus the full
field equals $K_0$. Consequently a rational function constant on very general
whole Iitaka fibers belongs to this full field, and nested full fields
of equal transcendence degree inside the same ambient function field
coincide.

\begin{lemma}\label{lem:ordinary-minimum-finite-leaves}
Let $g:Z\to Q$ be projective surjective with connected general
fibers, where $Z$ is smooth connected quasi-projective and $Q$ is a
principal torus torsor over an abelian variety. Put
$e=\dim Z-\dim Q>0$. Let $C$ be an effective rational divisor on
$Z$ with coefficients less than one, extending to an SNC pair with
the full boundary on a smooth projective compactification. Assume
that, for a very general smooth fiber $F$ of $g$, the divisor
$K_F+C|_F$ is big. Then the
following hold.
\begin{enumerate}[label=(\roman*)]
\item The dimension $k=\kappa_{\log}(Z,C)$ satisfies $k\ge e$.
\item There is a retained resolution of the full Iitaka map with
finite coefficients giving an effective exceptional adjoint error.
For a very general actual open leaf $H$, let $S_H$ be its reduced
closed image in $Q$. There are a smooth variety $V_H$, a proper
birational morphism $a_H:H\to V_H$, and a finite etale morphism
$v_H:V_H\to S_H^\nu$ factoring the specified map as
\[
 H\xrightarrow{a_H}V_H\xrightarrow{v_H}S_H^\nu\longrightarrow Q.
\]
Both $V_H$ and $S_H^\nu$ are principal torus torsors over abelian varieties.
Every positive finite boundary component on $H$ is $a_H$-exceptional,
and the specified map $\pi_1(H)\to\pi_1(Q)$ is injective.
\item After shrinking the parameter space, these normal models form
finite-type families, and one connected finite etale cover of $Q$
makes every connected lifted very general leaf proper birational
onto its whole image.
The connected stabilizer subtori and abelian image directions of
these images are constant on the family.
\end{enumerate}
\end{lemma}

\begin{proof}
\emph{Iitaka dimension and the leaves.}
Compactify the proper map over a toric compactification of $Q$. Deleted
primes receive coefficient one; retained exceptional primes of discrepancy
$a>-1$ receive rational coefficients in $(\max(0,-a),1)$. The adjoint
error is effective exceptional, so the section comparison in the
proof of Lemma~\ref{lem:canonical-kernel-logarithmic-comparison}
preserves all complete section spaces and their multiplication.
Since $Z\to Q$ is proper, the open graph is closed over $Q$ and the
inverse image of $Q$ is exactly the retained open.
Thus the reduced inverse image of the deleted base boundary has
coefficient one. The compactified map has connected fibers by its
regular generic field extension and normal Stein base, and its
general fiber adjoint remains big by the comparison with effective
exceptional error. By
Lemma~\ref{lem:canonical-kernel-logarithmic-comparison}, $Q$ has logarithmic
dimension zero. The compact pairs are log smooth and log canonical,
and the floor of the source boundary contains the reduced inverse
image of the deleted base boundary. The generic fiber adjoint is big:
on the relatively log-smooth locus, choose a very general fiber
outside the cohomology jumping loci of every divisible adjoint power,
so that its section dimensions equal the generic ones.
The first inequality of
\cite[Theorem~12, arXiv:1711.05854v1]{wei2020}
therefore gives $k\ge e$.
Resolve the full Iitaka map to a morphism
$\bar h:\overline Z'\to W$ with connected fibers and $W$ a smooth
projective $k$-fold, using the full-field observation after
Definition~\ref{def:ordinary-minimum-log-field} and the same coefficient
rule. Write $Z'$ for the retained open, $C'$ for its finite boundary,
and $h=\bar h|_{Z'}$. After clearing the denominators of the adjoint and its effective
exceptional error, we apply \cite[Notation~4.1 and Lemma~4.2]{taji2021}
to this Iitaka resolution. The restricted adjoint on a very general
compact fiber has Iitaka dimension zero.
Adjunction and transversality to the full deleted boundary identify
this restriction with the logarithmic adjoint of
$H=Z'\cap\bar h^{-1}(w)$. Thus $\kappa_{\log}(H,C'|_H)=0$.
This is a whole connected open leaf, closed in $Z'$ and proper over $Q$.

To prove that $H\to Q$ is generically finite, use the single proper
Stein diagram
\[
 Z'\xrightarrow{v}J\longrightarrow I\subset W\times Q
\]
of $(h,g)$, where $J$ is normal, $v$ has connected fibers, and
$J\to I$ is finite. The map $(h,g)$ is proper because $Z'\to Q$ is
projective and $W$ is separated. By the regularity of
$\CC(W)\subset\CC(Z')$ and generic flatness, $J_w$ is normal and
integral for general $w$. Hence $B_H=J_w$ is the actual normal Stein
base of $H\to Q$. Its resolution $\widehat B_H$ is already
generically finite over its image in $Q$, so
\[
 \overline\kappa(\widehat B_H)\ge0
\]
by Lemma~\ref{lem:canonical-kernel-logarithmic-comparison}.

Relative bigness and generic base change supply one divisible degree
whose complete relative section map has vertical differential rank
$e$ on a dense open of $Z'$. Choose $j\in J_w$ generally so that
$q=g(j)$ is smooth for $g$, $J\to Q$ is smooth at $j$, and the whole
fiber $T=v^{-1}(j)$ is smooth, connected, transverse to $C'$, and
meets this rank open. Such a choice is possible by generic smoothness and properness of
$v$ on the boundary strata. Put $F=g^{-1}(q)$.
The whole projective $T$ is a smooth fiber of $F\to J_q$, so its
normal bundle in $F$ is trivial and
\[
 (K_F+C'|_F)|_T=K_T+C'|_T.
\]
The map defined by the chosen complete sections on $F$ restricts to
$T$ with differential rank $\dim T$, so this fiber adjoint is big.
Take a retained resolution of the main base change over $\widehat B_H$.
Compatible compactifications of its two proper maps to $H$ and $\widehat B_H$
have exactly this open as both inverse images. The coefficient rule
preserves the zero logarithmic dimension on the source. The compact pairs
are again log smooth and log canonical,
and the reduced inverse image of the deleted base boundary lies in
the floor of the source boundary. The fiber adjoint is generically
big by the same very-general-fiber argument. Thus the first
inequality of \cite[Theorem~12, arXiv:1711.05854v1]{wei2020} gives
\[
 0\ge\dim T+\overline\kappa(\widehat B_H)\ge\dim T.
\]
Consequently $H\to Q$ is generically finite. Its unweighted logarithmic
dimension is at most zero by effectiveness of $C'|_H$, and at least
zero by the same generically finite comparison. Thus
$\overline\kappa(H)=0$.

\emph{Normal models and finite boundary.}
Apply Lemma~\ref{lem:canonical-kernel-zero-images} to this proper,
generically finite map. It identifies $S_H^\nu$, gives the proper
birational factor $H\to V_H$ followed by a finite etale map
$V_H\to S_H^\nu$, and proves injectivity of $\pi_1(H)\to\pi_1(Q)$.
By Proposition~\ref{prop:finite-cover-torus-torsor}, $V_H$ is a torus torsor.
If a positive finite boundary component on $H$ dominated a prime $J_H$
on $V_H$, choose its adapted compactification by
Lemma~\ref{lem:canonical-kernel-divisor-positivity}; there
$\kappa(\overline J_H)>0$. Compactify and resolve $a_H$ without deleting
interior primes, retaining the effective-exceptional comparison.
The adapted $\overline J_H$ contains no boundary stratum. Pull back
a logarithmic volume form of the torsor compactification. Its
contribution to the leaf adjoint is nonnegative everywhere. A retained
exceptional prime over the smooth interior has strictly positive
ordinary discrepancy. For a deleted exceptional prime over
$\overline J_H$, let its center have codimension $c$ and lie on $k$
boundary branches. Then $c>k$. With center parameters $x_1,\ldots,x_c$
and boundary parameters $x_1,\ldots,x_k$, the Jacobian wedge estimate
gives log discrepancy at least
$\sum_{i=k+1}^c\operatorname{ord}_E(x_i)>0$.
The strict transform of $J_H$ has the positive prescribed finite
coefficient. Taking the minimum of the finitely many positive
coefficient ratios over the components of $a_H^*\overline J_H$ gives
\[
 K+D+C\ge\epsilon a_H^*\overline J_H\qquad(\epsilon\in\QQ_{>0})
\]
on this common compactification, where $K$ is its canonical divisor
and $D+C$ its full deleted and finite boundary. This contradicts the preserved zero
Iitaka dimension of the leaf adjoint and proves (ii).

\emph{One finite cover for the family.}
Over a connected smooth open $W^0$, take the reduced closed image of
$(h,g):Z'|_{W^0}\to W^0\times Q$ for the retained Iitaka resolution
$Z'$. This map is proper and generically finite. Normalize its image
in its own field and in $\CC(Z')$. After shrinking $W^0$, geometric
normality, integrality and flatness identify their fibers with
$S_H^\nu$ and $V_H$: the birational isomorphism opens meet every fiber.
The nonsmooth and non-etale loci have constructible images; their
absence on very general fibers permits a further shrinking on which
both families are smooth and the finite map is etale of constant degree
$d$. Denote the diagram by
$\mathcal V\xrightarrow{b}\mathcal S\to W^0\times Q$, and its
evaluations by $e_V=e_Sb$.

Apply the single-family theorem \cite[Theorem~2.2]{achinger2019} to
$\mathcal S\to W^0$ on a connected dense stratum. The finite-cover
transport argument in the proof of
Lemma~\ref{lem:canonical-kernel-special-proper-image} applies:
path lifting gives compatible charts for $b$. Since $e_V=e_Sb$, the
evaluated base-point paths agree.
Thus transport conjugates the actual subgroup pairs by one common
element of $E_Q$:
\[
 P_H=\im\pi_1(H)\ \le\ N_H=\im\pi_1(S_H^\nu)
 \ \le\ E_Q=\pi_1(Q),\qquad [N_H:P_H]=d.
\]
Here the index is the actual finite etale degree: proper birationality
identifies $\pi_1(H)$ with $\pi_1(V_H)$, and the specified
$\pi_1(S_H^\nu)\to E_Q$ is injective.
By the torsor lattice sequence, $E_Q$ is finitely generated and
nilpotent. The subgroup-separation result used in the proof of
Lemma~\ref{lem:canonical-kernel-special-proper-image} therefore applies:
each element outside $P_H$ can be separated from $P_H$ in a finite
quotient of $E_Q$. Apply it to the finitely many nontrivial cosets
of $P_H$ in $N_H$.
For each coset representative,
take the preimage of the image of $P_H$ in a finite quotient separating
it from $P_H$. Their intersection is a finite-index subgroup $L$
containing $P_H$ with $L\cap N_H=P_H$.
Take its normal core $J_Q\triangleleft E_Q$.
Then $J_Q\cap N_H\subset P_H$, also for all transported conjugate
pairs. Normality gives, for every $n\in N_H$,
\[
 J_Q\cap N_H=n(J_Q\cap N_H)n^{-1}\subset nP_Hn^{-1}.
\]
On the associated algebraic finite etale cover of $Q$, each
connected lifted component, including every double-coset choice, has
finite normal degree
$[J_Q\cap N_H:J_Q\cap nP_Hn^{-1}]=1$.
Thus its leaf is proper birational onto its whole image. Finally,
$N_H$ meets the central torus lattice in exactly the stabilizer
cocharacter lattice, by Lemma~\ref{lem:canonical-kernel-zero-images}.
This intersection and the rational span of the image of $N_H$ in the
abelian lattice are invariant under conjugacy. They determine the fixed
stabilizer subtorus and abelian image direction, proving (iii).
\end{proof}

\begin{lemma}\label{lem:ordinary-minimum-mixed-leaves}
Let $f,c,g,Q,Z,B_c,R$ satisfy
Definition~\ref{def:prepared-ordinary-model}, and assume
$\dim Z-\dim Q>0$. Put $G=\pi_1(Y)$, and let $(Z_m,\Delta_m)$
be a prepared model whose adjoint on a very general projective
$g_m$-fiber is big. Then the following hold.
\begin{enumerate}[label=(\roman*)]
\item There is a rational $\epsilon\in(0,1)$ such that the full
weighted logarithmic Iitaka field equals that obtained by decreasing
every interior coefficient one by $\epsilon$. Its dimension $k$
is at least $\dim Z-\dim Q$.
\item There is a retained Iitaka resolution $U$ projective over $Q$ such
that the divisor defined by its meridian orders has no component dominating the
Iitaka base. If $H$ is a very general actual open Iitaka leaf, its
image $N$ in $R$ is normal and maps injectively to $\pi_1(Q)$.
If $K$ is the actual resolved general-fiber image of the composite
map from $Y$ to this Iitaka base, then
\[
 B_c\le K,\qquad K/B_c=N,\qquad G/K=R/N.
\]
Every element of $\ker(R\to\pi_1(Q))$ retains its exact order in
$R/N$, including infinite order.
\end{enumerate}
\end{lemma}

\begin{proof}
\emph{The full Iitaka field.}
Let $C$ be the closure of the original interior divisor $\Delta_m$ on
a compactification, and let $C_\infty$ be its coefficient-one part.
Choose rational $\epsilon\in(0,1)$ sufficiently small that
$K_F+(C-\epsilon C_\infty)|_F$ stays big. This works on very general
fibers: choose one fiber satisfying base change for every degree and
rational $\epsilon$, and use openness of its big cone. Fix a divisible
degree whose complete system on that fiber has full differential
rank. Generic base change and openness of that rank condition give
the same $\epsilon$ and degree on sufficiently general fibers.
Lemma~\ref{lem:ordinary-minimum-finite-leaves} applies to the
pair whose interior coefficients are less than one and gives Iitaka
dimension at least $\dim Z-\dim Q>0$.

Resolve its full Iitaka map. Give every retained exceptional prime of
discrepancy $a>-1$ a rational coefficient in $(\max(0,-a),1)$; deleted
primes keep coefficient one. The effective exceptional adjoint error
preserves sections, also after further retained resolutions. The adjoint
remains relatively big. By Lemma~\ref{lem:ordinary-minimum-finite-leaves},
every component of this auxiliary boundary on a very general
leaf $H$ is exceptional for the proper birational map $H\to V_H$ to
its smooth torsor model. Further resolutions preserve that model,
since image and function field are unchanged.

Any effective boundary supported on these exceptional primes has log
adjoint dimension zero on $H$: compactify the proper map to the toric
compactification of $V_H$. Its logarithmic volume form pulls back with an effective
exceptional divisor, so exceptional pushforward and
$H^0(\overline V_H,\mathcal O)=\CC$ imply that every divisible section
space is one-dimensional. Compare the full original adjoint upstairs,
assigning coefficient one to retained exceptional primes. Its restriction to $H$
has this form: the original boundary and the boundary with decreased
coefficients have the same strict support, and every exceptional prime
belongs to the auxiliary support. The complete sections of the original
adjoint identify with those
of this larger adjoint upstairs: log canonicity of the original full
pair ensures that assigning coefficient one to exceptional primes
introduces only effective exceptional error. Adjunction on the compact Iitaka fiber
identifies their restrictions with the stated one-dimensional
whole-leaf section space. In each divisible degree with nonzero
sections, choose a denominator and a basis containing it. On a very
general whole leaf the denominator does not vanish identically and
all basis ratios are constant. By the observation after
Definition~\ref{def:ordinary-minimum-log-field}, every such ratio
belongs to the full field for the decreased coefficients. Since this
holds in every degree, that field contains the original full field.
Conversely, multiplication by the effective boundary difference
preserves ratios and gives the reverse containment. The two full fields
are equal in $\CC(Z)$, proving (i).

\emph{Meridian orders and the quotient by the leaf group.}
Prepare this Iitaka model $U$. The support of $\Delta_U$ is contained in the
strict original support and the exceptional primes, hence in the
auxiliary support: strict-transform meridians keep their orders in
the same $R$, and all new retained exceptional primes have positive
auxiliary coefficients. First fill exactly the order-one primes in
the prepared complement and restore its codimension-two omissions.
Their meridians are trivial in $R$, so the prepared representation
factors as
\[
 \eta_U:\pi_1(U\setminus\Supp\Delta_U)\twoheadrightarrow R.
\]
Its composite to $\pi_1(Q)$ is the specified geometric map, by
compatibility on the bundle open and dense-open surjectivity.
On a leaf, the exceptional images and the
nonisomorphism locus in $V_H$ have codimension at least two.
Remove their union $A_H\subset V_H$. Its inverse image is a dense
open $O_H\subset H\setminus\Supp\Delta_U$, isomorphic to
$V_H\setminus A_H$. In the chain from $\pi_1(O_H)$ through
$\pi_1(H\setminus\Supp\Delta_U)$ and $\pi_1(H)$ to $\pi_1(V_H)$,
the dense-open maps are surjective and the composite is an isomorphism.
Thus
\[
 \pi_1(H\setminus\Supp\Delta_U)\simeq\pi_1(H)
 \simeq\pi_1(V_H).
\]
A horizontal prime in $\Supp\Delta_U$ would have a transverse meridian
in this leaf which becomes trivial when filled in $H$, hence is
already trivial in the displayed complement group. Under $\eta_U$
it would be trivial in the same $R$, contradicting its winding-one
meridian's actual order greater than one, including infinite order.
Thus $\Supp\Delta_U$ is vertical.

Choose a dense connected smooth base open avoiding the closures of
these finitely many vertical images, and apply
\cite[Theorem~2.2]{achinger2019} to the single open Iitaka morphism.
The whole bundle total space is dense in
$U\setminus\Supp\Delta_U$, so its group surjects onto $R$.
By the bundle sequence, the actual leaf image $N$ is normal. The
specified map $\pi_1(H)\to\pi_1(Q)$ is injective by
Lemma~\ref{lem:ordinary-minimum-finite-leaves}, hence so is $N\to\pi_1(Q)$.
Write $r:G\twoheadrightarrow R$ for the fixed quotient and take one
common retained resolution $M\to Y$ with $a:M\to U$. A whole
sufficiently general fiber $M_w$ of $ha$ is smooth connected, since
the two regular generic field extensions compose. It dominates
$H=h^{-1}(w)$ with geometrically connected general fiber, as seen
on the connected $a$-bundle locus. Hence
$\pi_1(M_w)\twoheadrightarrow\pi_1(H)$ by
\cite[Proposition~2.9.1]{kollar1993}. Prepared compatibility extends
from the dense bundle open to $a^{-1}(U\setminus\Supp\Delta_U)$
by dense-open surjectivity. It therefore holds on the whole $M_w$,
and gives $r(K)=N$.

Separately, choose a sufficiently general point $z\in H$ in the
global $a$-bundle open. The entire fiber $a^{-1}(z)$ lies in $M_w$
and has image exactly $B_c$ in $G$, by common retained resolution
of the original $c$-fiber. Thus $B_c\le K$. Together with
$r(K)=N$ and $\ker r=B_c$, this gives $K=r^{-1}(N)$,
$K/B_c=N$ and $G/K=R/N$.
Since $N\cap\ker(R\to\pi_1(Q))=1$, quotienting by $N$ changes
no finite or infinite order in the relative kernel. This proves (ii).
\end{proof}

\begin{theorem}\label{thm:cofinal-specialness-obstruction}
For the fixed special source and positive relative dimension in
Theorem~\ref{thm:ordinary-source-minimum}, bigness of the weighted
adjoint on a very general smooth projective fiber cannot hold on every
prepared model.
\end{theorem}

\begin{proof}
\emph{Choosing a model of minimal Iitaka dimension.}
Let $\mathcal C$ be the collection of all prepared models and assume
bigness on each. Under a projective
birational morphism of smooth prepared bases, weighted log sections in
divisible degree descend: strict transforms keep their meridian orders,
the boundary at infinity has coefficient one, and bounds at each prime
extend by normality. Common graphs give the same comparison between
compactifications.

Lemma~\ref{lem:ordinary-minimum-mixed-leaves} gives
$\kappa_{\log}(Z_m,\Delta_m)\ge\ell$. Let $b\ge\ell>0$ be the
minimum on $\mathcal C$, and choose a model attaining it. The same
lemma identifies its full Iitaka field with that obtained by decreasing
the interior coefficients equal to one, with a single decrease for
which the general fiber adjoint remains big. Take a common retained
Iitaka resolution with the finite-coefficient model of
Lemma~\ref{lem:ordinary-minimum-finite-leaves}, and prepare it for the
same quotient $R$, obtaining $U\in\mathcal C$ projective over $Q$.
The whole-leaf properties of
Lemma~\ref{lem:ordinary-minimum-mixed-leaves} persist under these
retained refinements. In particular $\Supp\Delta_U$ is vertical:
a new exceptional prime dominating the Iitaka base has center
generically outside the old vertical support, where its meridian is
null-homotopic in the old smooth complement.

Write $\mathcal F(V,D)$ for the full field in
Definition~\ref{def:ordinary-minimum-log-field}. The new actual
divisor $\Delta_U$ satisfies, by the strict-prime descent above,
\[
 \mathcal F(U,\Delta_U)\subseteq
 \mathcal F(Z_m,\Delta_m)\subset\CC(Z).
\]
The minimum on $\mathcal C$ and this descent give
$\kappa_{\log}(U,\Delta_U)=b$. Since the smaller full field is relatively
algebraically closed, as
proved after that definition, the two fields are equal. Thus the resolved
map $h:U\to B$ to the smooth
projective compact Stein base is also the full Iitaka map for
the actual $\Delta_U$. Fill its order-one primes and codimension-two
omissions as in Lemma~\ref{lem:ordinary-minimum-mixed-leaves}.
The image $N\triangleleft R$ of the leaf group injects into $\pi_1(Q)$, and
$N\cap\ker(R\to\pi_1(Q))=1$.

\emph{A common cover and torus torsor quotient.}
By Lemma~\ref{lem:ordinary-minimum-finite-leaves}, there is one finite
etale target cover on which every connected lifted very general leaf
maps with degree one onto its normal image. Their stabilizer subtori and abelian image
directions are constant. The pullbacks of $Y$ and $U$ are connected: their group maps
to $\pi_1(Q)$ are surjective, the latter by connectedness of its general
fiber and \cite[Proposition~2.9.1]{kollar1993}. General $c$-fibers lift
unchanged, so their subgroup remains $B_c$.
The whole images are closed because $(h,g):U\to B\times Q$ is
proper: $U\to Q$ is projective and $B$ is separated. After the
cover, use the compact Stein parameter of the connected lifted
leaves and the full inverse image of a dense parameter open. Their
evaluation onto the target cover stays dominant.
Their torsor normalizations give
$\dim S_H=\dim T_H+\dim A_{\rm dir}$, where $T_H$ is the fixed
connected stabilizer and $A_{\rm dir}$ is the fixed abelian image
direction. These are the whole-image hypotheses of
Lemma~\ref{lem:canonical-kernel-family-splitting}, which gives
a further abelian isogeny and a finite
etale cover $j:Q'\to Q$ with
\[
 Q'\longrightarrow A_{\rm dir}\times Q_1,\qquad r:Q'\to Q_1,
\]
where $A_{\rm dir}$ is the fixed abelian image direction,
the first map is a principal torsor under a connected torus, and every
very general connected lifted leaf is proper birational onto a whole
$r$-fiber. Proposition~\ref{prop:finite-cover-torus-torsor} identifies
$Q'$ as a torus torsor over an abelian variety.

For this cover, put
\[
 \begin{gathered}
 E_Q'=j_*\pi_1(Q'),\qquad Y'=Y\times_QQ',\qquad U'=U\times_QQ',\qquad
 G'=\pi_1(Y')=f_*^{-1}(E_Q'),\qquad P_0\le G',\\
 R'=G'/B_c=(R\to\pi_1(Q))^{-1}(E_Q')\le R.
 \end{gathered}
\]
Since the specified $\pi_1(H)\to R$ identifies a whole leaf group
with $N$, a connected lift corresponds to its inverse image of
$E_Q'$. Its actual image is therefore $N'=N\cap R'$. By normality, changing
the compatible base-point paths does not change this identity, and
$N'\cap\Pi=1$.
Interior meridians lie in the relative kernel:
they bound transverse disks on the smooth base mapping to $Q$.
Their lifts have winding one and unchanged orders in $R'$. The subgroup
$N'$ remains normal and injects into $\pi_1(Q')$.

For the full field under this cover, take a connected Galois closure
with deck group $A$ and compatible finite normal compactifications.
At a deleted prime of ramification index $e$, pluricanonical order
$v$ pulls back to $ev+m(e-1)$ in degree $m$; the coefficient-one
bound is therefore equivalent downstairs and upstairs. At retained
primes the cover is etale. The deleted-exceptional comparison in
the proof of Lemma~\ref{lem:canonical-kernel-logarithmic-comparison}
gives the remaining pullback bounds. Invariant sections descend by
their bounds at every target prime and normality. For an invariant
full ratio $t=s_1/s_0$, the equal-degree sections
$s=\prod_{\sigma\in A}\sigma(s_0)$ and
$ts=s_1\prod_{\sigma\ne1}\sigma(s_0)$ are invariant sections of the
same degree and descend. Hence the invariant full field is the
original one. Taking intermediate invariants handles the possibly
non-Galois $Q'\to Q$.

The full field on $U'$ is consequently finite over the original
one; its relative algebraic closedness identifies it with the field
of the compact Stein model $B'$ of the lifted Iitaka map. Thus
$\dim B'=b$, and its very general whole fibers are the connected lifted leaves.
The map $rg$ is constant on them, so it descends rationally to $B'$.
Resolve its graph in $B'\times Q_1$ to obtain a smooth $b$-fold $W$
projective over $Q_1$, with a birational morphism $W\to B'$.
The smooth model
$E=Q'\times_{Q_1}W$
is projective over $Q'$. The induced rational map
$U'\dashrightarrow E$ is birational: on a very general
$W$-fiber it is the proper birational map onto the whole $r$-fiber.
For the actual resolved general-fiber image $K'$ of
$Y'\dashrightarrow W$, the whole composite fiber maps dominantly
with geometrically connected general fiber to the whole lifted leaf.
Its group surjects onto that leaf group, giving
$(G'\to R')(K')=N'$. Separately it contains an entire sufficiently
general original $c$-fiber, whose connected lifts are unchanged.
As in Lemma~\ref{lem:ordinary-minimum-mixed-leaves},
$B_c\le K'$, $K'/B_c=N'$ and $G'/K'=R'/N'$.
Every element of $\ker(R'\to\pi_1(Q'))$ retains its exact order
in $R'/N'$.

\emph{Meridian orders on the quotient models.}
The source graph over $Y'$ is proper, since $E$ is projective over
$Q'$. Its bundle sequence factors the quotient of $\pi_1(Y')$ through
the group of the base open, compatibly on common graphs. Compute the
meridian orders on $E$ in $R'$. A prime horizontal over $W$ has order one.
Indeed its DVR has a center on $U'$ by properness of $U'\to Q'$.
The two maps from this DVR to $B'$, through $U'$ and through $W$,
agree generically and hence everywhere by separatedness. Its center
therefore dominates the same compact Stein parameter. On a very
general leaf the prime maps to a divisor of the whole $r$-fiber.
Its center on $U'$ cannot have codimension at
least two: the center's dimension on that leaf would be at most the
leaf dimension minus two, smaller than its divisorial image in $Q'$.
It is therefore a transform of a horizontal prime of $U'$, with
the same DVR and meridian. Its
order is one by Lemma~\ref{lem:ordinary-minimum-mixed-leaves}.

The projection $p:E\to W$ is smooth with geometrically integral
fibers. By the dimension bound, every vertical prime dominates a base prime.
The inverse image of that base prime is reduced and irreducible by
smoothness and geometric integrality. Thus $\Delta_E=p^*\Theta_W$.
After filling order-one primes, the image of a whole $p$-fiber is
exactly $N'$: its proper birational map from a whole lifted leaf
induces a fundamental-group isomorphism, compatible with the
representations on their common dense open. Modulo $N'$, the bundle
sequence gives a representation of the base complement onto
$R'/N'=G'/K'$. A local section identifies a base meridian with a
meridian of its reduced inverse image, modulo a fiber element. The
order is unchanged in $R'/N'$, so the weights in $\Theta_W$ are the
orders in $G'/K'$, including infinity.
For every smooth $W'$ with a projective birational morphism $W'\to W$,
apply this argument to $E'=Q'\times_{Q_1}W'$, keeping $U'$ and the
maps to $B'$ fixed. This gives
$\Delta_{E'}=p'^*\Theta_{W'}$ and an epimorphism
$\beta_{W'}:\pi_1((W')^\circ)\twoheadrightarrow G'/K'$
onto the fixed quotient $G'/K'$. The whole generic composite fiber
is unchanged, so $K'$ is unchanged as well. The weight support is finite, contained in a
finite bundle complement; further resolutions of this complement
make it SNC. New primes outside its total transform have order one.
In the following comparisons first take such an SNC higher model.

\emph{Bounding the Iitaka dimension on each quotient model.}
Fix $W'$. Choose a connected finite Galois target cover
$\widetilde Q\to Q$ dominating $Q'$, and write
$\mathcal A=\operatorname{Gal}(\widetilde Q/Q)$ and
$\mathcal H=\operatorname{Gal}(\widetilde Q/Q')$. Resolve equivariantly the graph of
all deck translates of the rational map from $U\times_Q\widetilde Q$
to the corresponding pullback of $E'$. All factors are projective over
$\widetilde Q$. The resulting $J$ has a free deck action, since it maps
equivariantly to $\widetilde Q$. Its quotient $J_0$ is smooth,
projective over $Q$, and birational to $U$; projectivity follows by
descending a power of the tensor product of the translates of a relative
ample bundle. Etale descent gives $J=J_0\times_Q\widetilde Q$.
The map to the factor indexed by the identity was defined over $Q'$,
so it is $\mathcal H$-equivariant. By separatedness, its extension
is also equivariant.
Taking this intermediate quotient gives a proper birational map
$J_0\times_QQ'\to E'$.

Prepare $J_0$ for the original $c$, obtaining a dominating
$T_0\in\mathcal C$. Its connected pullback $T'_0=T_0\times_QQ'$
maps properly birationally to $E'$. Its divisor of meridian weights is
the etale pullback of the original one, since relative meridians lift with
winding one and unchanged order. The ramification and
deleted-exceptional comparison for $U'\to U$ applies to
$T'_0\to T_0$: all deleted coefficients are one, and the map is
etale at retained primes. Every divisible weighted log section pulls back to a section, and
pullback commutes with multiplication and taking ratios. Push them from
$T'_0$ to $E'$ using the comparison at each prime from the first
paragraph. This gives
\[
 \kappa_{\log}(E',p'^*\Theta_{W'})\ge b.
\]
Since $E'\to A_{\rm dir}\times W'$ is a principal torus torsor,
Lemma~\ref{lem:canonical-kernel-logarithmic-comparison} identifies its
section ring, up to a constant volume factor, with that of
$(W',\Theta_{W'})$. Since $\dim W'=b$, the full weighted log adjoint
of $(W',\Theta_{W'})$ is big for every $W'$.
In particular the argument applies directly to every higher smooth
model on which the actual weight support and full deleted boundary
are SNC. For an initially non-SNC model, descend the sections from
an SNC refinement at its strict primes to obtain the same lower
bound for its adjoint. The comparison graphs may depend on $W'$, whereas
$Q',Y',B_c,R',N'$ and $K'$ remain fixed.

\emph{Pullback to the special source.}
The finite etale cover $Y'$ is special
\cite[Proposition~10.11]{campana2009}.
Fix $X_*=\overline Y'$, with full reduced SNC boundary $D_Y$,
and apply the fixed-source preparation in the
proof of Lemma~\ref{lem:canonical-kernel-ordinary-exactness} to the
dominant rational map $X_*\dashrightarrow\overline W$, for a
projective compactification $\overline W$ of $W$.
After this construction, every original prime of $X_*$, including
every deleted prime, dominates the base or a base prime.
Its restricted valuation is respectively trivial or
$e\,\operatorname{ord}_D$ with $e>0$; this persists on every higher
proper normal base. Resolve the base, its infinity and the finite
recomputed meridian support. Over $W$ we obtain a smooth $W'$ with a
projective birational morphism $W'\to W$, whose weighted log
adjoint is big.
Any source prime still mapping into base codimension at least two
is exceptional over $X_*$.

Every retained source prime maps into $W'$: its generic DVR maps to
$Q_1$ through $rf'$, and properness of $W'\to Q_1$ extends its rational
lift into $W'$. Separatedness identifies this lift with the prepared
map. The quotient $G'\to G'/K'$ agrees with the base
representation on every source open avoiding the boundary: the two
agree on a common bundle open, whose inclusion induces a surjection on
fundamental groups. In particular, all retained source divisors
supply their relations.
If such a prime dominates a base prime of meridian order $m_D$ with
multiplicity $e$, its transverse disk gives
$\beta_{W'}(\mu_D)^e=1$ in the fixed $G'/K'$. Thus $m_D$
is finite and divides $e$; an infinite-order prime has no retained
source prime above its generic point.

Write $\varphi:X_*\dashrightarrow\overline W'$ for the rational map,
and let $\mathcal L$ be the saturation of its generic differential
line in $\Omega^b_{X_*}(\log D_Y)$. It is rank-one reflexive,
because its quotient in this locally free sheaf is torsion-free;
it is therefore invertible on the smooth $X_*$. At each divisorial
DVR choose its primitive generator in the ambient logarithmic
differential module. Its tensor powers remain primitive, so the
zeros of the pulled-back differential remain as zeros of its
coefficient in $\mathcal L$.

Since $Y'$ is special and $1\le b\le\dim Z\le\dim Y'$, the
logarithmic differential inequality for the saturated line on
$(X_*,D_Y)$ gives $\kappa(X_*,\mathcal L)<b$.

Pull back divisible weighted log sections and check every original
prime of $X_*$. A prime dominating the base has no pole
from a base divisor. At a retained prime over a base prime of finite weight,
$t=x^eu$ gives differential vanishing
\[
 e-1\ge e(1-1/m_D).
\]
At a deleted prime,
$d(x^eu)=x^e(eu\,dx/x+du)$ vanishes to order at least $e$ in a
logarithmic frame. This suffices for every coefficient at most one,
including the boundary at infinity. More generally, if several
target parameters have positive orders $e_i$ at a deleted prime,
write their pullbacks as $x^{e_i}u_i$. The whole $b$-fold differential
wedge is divisible by $x^{\sum_i e_i}$ in the logarithmic frame.
After subtracting the weighted poles its order is at least
$\sum_i(1-\gamma_i)e_i\ge0$, since every coefficient $\gamma_i\le1$;
unweighted parameters have $\gamma_i=0$.
Normality therefore gives the natural injections
\[
 H^0\!\left(\overline W',m(K_{\overline W'}+D_{W'}+
 \overline\Theta_{W'})\right)
 \hookrightarrow H^0(\overline Y',\mathcal L^{\otimes m})
\]
in every sufficiently divisible degree, all from the same generic
differential pullback. They preserve multiplication, and for equal
degree sections $s,t$ the quotient of their images is
$\varphi^*(s/t)$. Since the base adjoint is big, a sufficiently
divisible complete base system has image dimension $b$. Dominance
and this equality of ratios give an image subsystem of
$\mathcal L^{\otimes m}$ with the same image dimension. Hence
$\kappa(\overline Y',\mathcal L)\ge b$, contradicting specialness.
\end{proof}

\begin{proof}[Proof of Theorem~\ref{thm:introduction-main}]
Choose the canonical second Albanese map $f_0:Y_0\to Q_0$ on a connected finite
etale cover of $X$ by Theorem~\ref{thm:canonical-albanese}. Fix
$G_0=\pi_1(Y_0)$, $P_0=\ker(f_{0*})$ and
$E_0=G_0/P_0=\pi_1(Q_0)$.
If $P_0$ is infinite,
Theorem~\ref{thm:ordinary-source-minimum} gives a fixed relative quotient
of positive relative dimension on a finite etale cover $Y\to Y_0$.
The source $Y$ is special
\cite[Proposition~10.11]{campana2009}.
By Theorem~\ref{thm:source-minimal-bigness}, the weighted adjoint is big
on a very general projective fiber of every model in this prepared
collection. This contradicts
Theorem~\ref{thm:cofinal-specialness-obstruction}. Both statements
allow the very general point to depend on the model.
Thus $P_0$ is finite. Applying
Proposition~\ref{prop:finite-kernel-criterion} to $f_0$ gives a
finite-index subgroup $H\le G_0$ of class at most two. Finally,
\[
 [\pi_1(X):H]=[\pi_1(X):G_0]\,[G_0:H]<\infty.
\]
\end{proof}

\section*{Acknowledgments}
\addcontentsline{toc}{section}{Acknowledgments}

We thank Richard Hain for his careful reading of our manuscript
and for pointing out inaccuracies in our citations.

\section*{AI disclosure}
\addcontentsline{toc}{section}{AI disclosure}

OpenAI's proof of the
abelianity conjecture for special compact K\"ahler manifolds
\cite{openai2026d} inspired us to develop an argument
for the quasi-projective case. Both proofs reduce to showing that a
finitely generated FAb kernel is finite. They choose a base of minimal
dimension retaining an infinite quotient and define boundary
coefficients by actual meridian orders. Descent along the Iitaka
fibration then yields a contradiction with specialness.

We realize the second Albanese map as a morphism to an algebraic torus
torsor and identify the kernel of the induced homomorphism on
fundamental groups. We treat logarithmic boundaries with
possibly infinite meridian orders and establish the generalized log
canonical antinef group theorem (Theorem~\ref{thm:antinef-group}),
including its non-klt case. We combine this group theorem with
OpenAI's existence theorem for good log minimal models
\cite{openai2026b} and the Mori fiber space case of the minimal model
program to obtain the bigness needed in our proof.

The human authors conceived the new mathematical ideas underlying
these additional arguments and directed the mathematical work.
We thank ChatGpt-6 Astra for help with discussion and revision of the
arguments, and take full responsibility for the reasoning, results,
and proofs.

\bigskip
\noindent\begin{minipage}{\textwidth}
\small
\textsc{Tianshu Jiang}\\
Nanjing University\\
\textit{Email:} \href{mailto:jts2026@nju.edu.cn}{\nolinkurl{jts2026@nju.edu.cn}}

\medskip
\textsc{Jiayu Li}\\
Nanjing University\\
\textit{Email:} \href{mailto:jiayuli@nju.edu.cn}{\nolinkurl{jiayuli@nju.edu.cn}}
\end{minipage}
\end{document}